\documentclass[a4paper,11pt,leqno]{article}
\usepackage{amsfonts,amsmath,amscd}
\usepackage{amssymb}
\usepackage{amssymb,amsmath,amsthm,amscd,amsfonts}
\usepackage{amsmath,amssymb,amsthm,amscd}
\usepackage{graphicx}

\numberwithin{equation}{section}
\theoremstyle{remark}

\begin{document}
\title{\textbf{Classification of isoparametric submanifolds\\
admitting a reflective focal submanifold\\
in symmetric spaces of non-compact type}}
\author{Naoyuki Koike}
\date{\today}
\maketitle

%
%
%

\begin{abstract}
In this paper, we assume that all isoparametric submanifolds have flat section.  
The main purpose of this paper is to prove that, if a full irreducible complete isoparametric submanifold 
of codimension greater than one in a symmetric space of non-compact type admits a reflective focal submanifold 
and if it of real analytic, then it is a principal orbit of a Hermann type action on the symmetric space.  
A hyperpolar action on a symmetric space of non-compact type admits a reflective singular orbit if and only if 
it is a Hermann type action.  
Hence is not extra the assumption that the isoparametric submanifold admits a reflective focal submanifold.  
Also, we prove that, if a full irreducible complete isoparametric submanifold of codimension greater than one 
in a symmetric space of non-compact type satisfies some additional conditions, then it is a principal orbit of 
the isotropy action of the symmetric space, where we need not impose that the submanifold is of real analytic.  
We use the building theory in the proof.  
\end{abstract}

\vspace{0.1truecm}


\hspace{2truecm}{\small geneity, anti-Kaehler isoparametric submanifold, Hermann type action}



\vspace{0.3truecm}

\section{Introduction} 
In 1985, C. L. Terng (\cite{Te1}) introduced the notion of an isoparametric submanifold (of general codimension) 
in a Euclidean space and, in 1995, C. L. Terng and G. Thorbergsson (\cite{TT}) introduced the notion of 
an {\it equifocal submanifold} in a (Riemannian) symmetric space as its generalized notion.  This notion is 
defined as a compact submanifold with flat section, trivial normal holonomy group and parallel focal structure.  
Here "with flat section" means that the images (whcih is called the normal umbrella) of the normal spaces of 
the submanifold by the normal exponential map are flat totally geodesic submanifolds and ``the parallelity of 
the focal structure'' means that, for any parallel normal vector field $\widetilde v$ of the submanifold, 
the focal radii along the normal geodesic $\gamma_{\widetilde v_x}$ with 
$\gamma_{\widetilde v_x}'(0)=\widetilde v_x$ 
are independent of the choice of $x$ (with considering the multiplicities), where $\gamma_{\widetilde v_x}'(0)$ 
is the velocity vector of $\gamma_{\widetilde v_x}$ at $0$.  
Compact isoparametric hypersurfaces in a sphere or a hyperbolic space are equifocal.  
E. Heintze, X. Liu and C. Olmos (\cite{HLO}) defined the notion of an {\it isoparametric submanifold} in a general 
complete Riemannian manifold as a (properly embedded) complete submanifold with section and trivial normal holonomy 
group whose sufficiently close parallel submanifolds are of constant mean curvature with respect to the radial 
direction.  Here ``with section`` means that the normal umbrellas of the submanifold are totally geodesic 
(the normal umbrellas are called {\it sections}).  

\vspace{0.5truecm}

\noindent
{\bf Assumption.} In this paper, we assume that all isoparametric submanifolds have flat section, that is, 
the induced metric on the sections are flat.  

\vspace{0.5truecm}

For a compact submanifold in a symmetric space of compact type, they (\cite{HLO}) proved that it is equifocal 
if and only if it is an isoparametric submanifold (with flat section).  
In 1989, C. L. Terng (\cite{Te2}) introduced the notion of an isoparametric submanifold in a (separable) Hilbert 
space and intiated its research.  In 1995, C. L. Terng and G. Thorbergsson (\cite{TT}) proved that the research of 
an equifocal submanifold in a symmetric space $G/K$ of compact type is reduced to that of an isoparametric 
submanifold in the Hilbert space $H^0([0,1],\mathfrak g)$ through the composition of the parallel 
transport map $\phi:H^0([0,1],\mathfrak g)\to G$ for $G$ and the natural projection $\pi:G\to G/K$, where 
$\mathfrak g$ denotes the Lie algebra of $G$ and $H^0([0,1],\mathfrak g)$ denotes the path space of 
all $L^2$-integrable curves(=pathes) in $\mathfrak g$.  
Denote by $I(V)$ the group of all isometries of a (separable) Hilbert space $V$, where we note that $I(V)$ is not 
a Banach Lie group (see \cite{Ha} and \cite{KW} (Appendix of \cite{Koi3} also)).  
Let $\widetilde M$ be a full irreducible complete isoparametric submanifold of codimension greater 
than one in $V$.  Here ``completeness'' means ``metric completeness'', where we note that, for 
a Riemannian Hilbert manifold, if it is metrically complete, then it is also geodesically complete, but 
the convese does not necessarily hold (see \cite{A}).  In main theorems of \cite{Koi12}, \cite{Koi13} and 
\cite{Koi14}, we assumed that the submanifolds are metrically complete as anti-Kaehler Hilbert manifolds 
without mentioned.  In Section 3, we shall state the definition of metric completeness of an anti-Kaehler Hilbert 
manifold.  Throughout this paper, ``completeness'' means ``metric completeness'' and we shall write it 
without abbreviated.  Set $H:=\{F\in I(V)\,|\,F(\widetilde M)=\widetilde M\}$.  
In 1999, E. Heintze and X. Liu (\cite{HL}) proved that $\widetilde M$ is extrinsically homogeneous in the sense 
that $Hu=\widetilde M$ holds for any $u\in M$.  
This result is the infinite dimensional version of the extrinsic homogeneity theorem for 
a finite dimensional compact isoparametric submanifold in a Euclidean space by G. Thorbergsson 
(\cite{Th}).  The extrinsic homogeneity theorem in \cite{Th} states that full irreducible compact isoparametric 
submanifolds of codimension greater than two in a Euclidean space are extrinsically homogeneous.  
G. Thorbergsson proved this statement by constructing the topological Tits building of spherical type associated 
to the isoparametric submanifold (in more general, he defined this topological Tits building for full irreducible 
isoparametric submanifolds of rank greater than one in a Euclidean space) and using it, where we note that, 
if the isoparametric submanifold is a principal orbit of the s-representation of an irreducible symmetric space 
$G/K$ of non-compact type and rank greater than one, then its associated topological Tits building coincides 
with the Tits building of the semi-simple Lie group $G$ (which is defined as the building having parabolic subgroups 
as vetices) as Tits building.  
Later, C. Olmos (\cite{O2}) proved this result by Thorbergsson in simpler method, that is, 
by constructing the normal homogeneous structure for the isoparametric submanifold and using the result in \cite{OS} 
(without use of the above topological Tits building), where the normal homogeneous structure means a certain kind of 
connection on the Whitney sum of the tangent bundle and the normal bundle of the submanifold.  E. Heintze and X. Liu 
(\cite{HL}) proved the above extrinsic homogeneity theorem in the method similar to the proof in \cite{O2}.  
In 2002, by using the extrinsic homogeneity theorem of Heintze-Liu, U. Christ (\cite{Ch}) proved that a full 
irreducible equifocal submanifold of codimension greater than one in a simply connected symmetric space of compact 
type is extrinsically homogeneous.  However, there was a gap in his proof because the above group $H$ in the theorem 
of Heintze-Liu is not Banach Lie group but he interpreted it as a Banach Lie group.  
Let $I_b(V)$ be the subgroup of $I(V)$ generated by one-parameter transformation groups induced by 
the Killing vector fields defined entirely on $V$.  It is easy to show that $I_b(V)$ is a Banach Lie group.  
Set $H_b:=H\cap I_b(V)$, which is a Banach Lie subgroup of $I_b(V)$.  
Recently, C. Gorodski and E. Heintze (\cite{GH}) proved that $\widetilde M$ is extrinsically homogeneous in the sense 
that $H_bu=\widetilde M$ holds for any $u\in \widetilde M$.  This improved extrinsic homogeneity theorem closed the 
gap in the proof of Christ.  According to the extrinsic homogeneity theorem by Christ and Theorem 2.3 in \cite{HPTT}, 
we can derive that, if $M$ is an irreducible equifocal(=isoparametric) submanifold of codimension greater than one 
in a simply connected symmetric space of compact type, then it is a principal orbit of a hyperpolar action.  
On the other hand, according to the classification of the hyperpolar actions by A. Kollross (\cite{Kol}), 
all hyperpolar actions of cohomogeneity greater than one on the irreducible symmetric space of compact type 
are Hermann actions.  Also, O. Goertsches and G. Thorbergsson (\cite{GT}) proved that principal orbits of 
Hermann actions are curvature-adapted, where ``curvature-adaptedness'' means that, for any unit normal vector $v$ 
of $M$, the normal Jacobi operator $R(v)$ preserves $T_xM$ ($x:$the base point $v$) invariantly and that 
$R(v)$ commutes with the shape operator $A_v$, where $R$ is the curvature tensor of the ambient symmetric space 
and $R(v):=R(\cdot,v)v$.  From these facts, we can derive the following fact:

\vspace{0.3truecm}

{\sl All complete equifocal(=isoparametric) submanifolds of codimension greater than one in simply}

{\sl connected irreducible symmetric spaces of compact type are principal orbits of Hermann actions}

{\sl and they are curvature-adapted.}

\vspace{0.3truecm}

In 2000, by the discussion with G. Thorbergsson at Nagoya University (The 47-th Geometry Symposium), 
the author was very interesting in the following open problem:

\vspace{0.2truecm}

{\it Is there a similar theory for equifocal submanifolds in simply connected non-compact symmetric 

spaces?}

\vspace{0.2truecm}

\noindent
This is one of seven open problems proposed in \cite{TT}.  
The author interpreted that this open problem means the following:

\vspace{0.2truecm}

{\it Can we reduce the study of an equifocal submanifold in a simply connected non-compact symme-

tric space to the study of the lift of the submanifold to a Hilbert space through a Riemannian 

submersion (of the Hilbert space onto the symmetric space) or the study of the lift of some 

extended submanifold of the original submanifold (which is a submanifold in some extended 

symmetric space of the original symmetric space) to some pseudo-Hilbert space through a pseudo-

Riemannian submersion (of the pseudo-Hilbert space onto the extended symmetric space)?}

\vspace{0.2truecm}

\noindent
Under this motivation, the author introduced the notion of a complex equifocal submanifold in a symmetric space of 
non-compact type and started its study.  
We shall explain why we introduced the notion of a complex equifocal submanifold in a symmetric 
space of non-compact type.  When a non-compact submanifold $M$ in a symmetric space $G/K$ of non-compact type is 
deformed as its principal curvatures approach to zero, its focal set vanishes beyond the ideal boundary 
$(G/K)(\infty)$ of $G/K$.  
For example, when an open portion of a totally umbilic sphere (whose only principal curvature is greater than 
$\sqrt{-c}$) in a hyperbolic space of constant curvature $c(<0)$ is deformed as its principal curvatures 
approach to $\sqrt{-c}$ (remaining to be totally umbilic), its focal point approach to $(G/K)(\infty)$ 
and, when it furthermore is deformed as its principal curvatures approach to 
a positive value smaller than $\sqrt{-c}$ (remaining to be totally umbilic), the focal point vanishes beyond 
$(G/K)(\infty)$.  According to these facts, we recognized that, for a non-compact submanifold in a symmetric space 
of non-compact type, the parallelity of the focal structure is not an essential condition.  
So, we (\cite{Koi2}) introduced the notion of a {\it complex focal radius} of the 
submanifold along the normal geodesic $\gamma_{v}$.  
See Section 2 about the definition of this notion.  
Furthermore, we (\cite{Koi2}) defined the notion of a {\it complex equifocal submanifold} 
as a (properly embedded) complete submanifold with flat section, trivial normal holonomy group and 
parallel complex focal structure, where we note that this submanifold should be called 
an equi-complex focal submanifold but we called it a complex equifocal submanifold for the simplicity.  
We proved that all isoparametric submanifolds (in the sense of \cite{HLO}) are complex equifocal 
and that, conversely all curvature-adapted complex equifocal submanifolds are isoparametric 
(see Theorem 15 of \cite{Koi3}).  Thus, for a complete submanifold in $G/K$, it is cuvature-adapted complex 
equifocal if and only if it is curvature-adapted isoparametric.  Hence, throuhout this paper, we shall use the 
terminology ``curvature-adapted isoparametric'' more familiar than ``curvature-adapted complex equifocal''.  

We consider the case where $M$ is of class $C^{\omega}$ (i.e., real analytic).  
Then we (\cite{Koi3}) defined the complexification $M^{\mathbb C}$ of $M$ as an anti-Kaehler submanifold 
in the anti-Kaehler symmetric space $G^{\mathbb C}/K^{\mathbb C}$.  Here we note that $G^{\mathbb C}/K^{\mathbb C}$ 
is a space including both $G/K$ and its compact dual $G_{\kappa}/K$ as submanifolds transversal to each other 
and that it is interpreted as the complexification of both $G/K$ and $G_{\kappa}/K$.  Also we note that the induced 
metric on $G/K$ coincides with the original metric of $G/K$ and that the induced metric on $G_{\kappa}/K$ is 
the $(-1)$-multiple of the metric of the symmetric space $G_{\kappa}/K$ of compact type.  
We (\cite{Koi3}) showed that 
$z$ is a complex focal radius of $M$ along the normal geodesic $\gamma_{v}$ if and only if 
$\gamma_{v}^{\mathbb C}(z)$ is a focal point of $M^{\mathbb C}$ along the complexified geodesic 
$\gamma_{v}^{\mathbb C}$.  Here $\gamma_{v}^{\mathbb C}$ is defined by 
$\gamma_{v}^{\mathbb C}(z):=\gamma_{av+bJv}(1)$ ($z=a+b\sqrt{-1}\in{\mathbb C}$), where 
$J$ denotes the complex structure of $G^{\mathbb C}/K^{\mathbb C}$ and $\gamma_{av+bJv}$ is the geodesic in 
$G^{\mathbb C}/K^{\mathbb C}$ with $\gamma_{av+bJv}'(0)=av+bJv$.  
Thus the complex focal radii of $M$ are the quantities indicating the positions of focal points of $M^{\mathbb C}$.  

We (\cite{Koi3}) introduced the notion of an {\it anti-Kaehler isoparametric submanifold} in the infinite 
dimensional anti-Kaehler space and furthermore, defined the {\it parallel transport map} for $G^{\mathbb C}$ as 
an anti-Kaehler submersion of an infinite dimensional anti-Kaehler space (which is denoted by 
$H^0([0,1],\mathfrak g^{\mathbb C})$) consisting of certain kind of paths in the Lie algebra 
$\mathfrak g^{\mathbb C}$ of $G^{\mathbb C}$ onto $G^{\mathbb C}$.  Denote by $\phi$ the parallel transport map 
for $G^{\mathbb C}$.  
We (\cite{Koi3}) proved that the research of a complex equifocal $C^{\omega}$-submanifold in a symmetric space 
$G/K$ of non-compact type is reduced to that of an anti-Kaehler isoparametric submanifold in the infinite dimensional 
anti-Kaehler space $H^0([0,1],\mathfrak g^{\mathbb C})$ by lifting the complexification of the original 
submanifold through the composition of the parallel transport map 
$\phi:H^0([0,1],\mathfrak g^{\mathbb C})\to G^{\mathbb C}$ for $G^{\mathbb C}$ and the natural projection 
$\pi:G^{\mathbb C}\to G^{\mathbb C}/K^{\mathbb C}$.  More precisely, we showed that a $C^{\omega}$-submanifold in 
$G/K$ is complex equifocal if and only if the lift of its complexification to $H^0([0,1],\mathfrak g^{\mathbb C})$ 
is anti-Kaehler isoparametric.  

We (\cite{Koi12}) proved that any full irreducible (metrically) complete anti-Kaehler isoparametric 
$C^{\omega}$-submanifold $\widetilde M$ with $J$-diagonalizable shape operators of codimension greater than one in an 
infinite dimensional anti-Kaehler space $V$ is extrinsically homogeneous in the sense that $Hu=\widetilde M$ holds 
for any $u\in\widetilde M$ as $H:=\{F\in I_h(V)\,|\,F(\widetilde M)=\widetilde M\}$, where $I_h(V)$ denotes 
the group of all holomorphic isometries of $V$ and ``with $J$-diagonalizable shape operators'' means that 
the complexifications of the shape operators are diagonalized with respect to $J$-orthonormal bases.  
Note that we assumed that, in main theorem, an anti-Kaehler isoparametric submanifold is (metrically) complete 
without mentioned.  Recently we (\cite{Koi13}) improved this extrinsic homogeneity 
theorem as follows.  

\vspace{0.3truecm}

\noindent
{\bf Fact 1.1.} {\sl Let $\widetilde M$ be a full irreducible complete anti-Kaehler isoparametric 
$C^{\omega}$-submanifold with $J$-diagonalizable shape operators of codimension greater than one in an infinite 
dimensional anti-Kaehler space $V$.  Then $\widetilde M$ is extrinsically homogeneous in the sense that 
$H_bu=\widetilde M$ holds for any $u\in\widetilde M$ as $H_b:=\{F\in I^b_h(V)\,|\,F(\widetilde M)=\widetilde M\}$, 
where $I_h^b(V)$ denotes the subgroup of $I_h(V)$ generated by one-parameter transformation groups induced by 
holomorphic Killing vector fields defined entirely on $V$.}

\vspace{0.3truecm}

Here we note that $I_h^b(V)$ is a Banach Lie group and $H_b$ is a Banach Lie subgroup of $I_h^b(V)$.  
Let $M$ be a complete curvature-adapted submanifold with flat section in a symmetric space $G/K$.  
If $G/K$ is of compact type or Euclidean type, then the following fact $(\ast_{\mathbb R})$ holds:

\vspace{0.3truecm}

\hspace{0.1truecm}$(\ast_{\mathbb R})\,\,$
{\sl For any unit normal vector $v$ of $M$, the nullity spaces for focal radii along}

\hspace{1.15truecm}{\sl 
the normal geodesic $\gamma_v$ span $T_xM\ominus({\rm Ker}\,A_v\cap{\rm Ker}\,R(v))$.}

\vspace{0.3truecm}

\noindent
Here $T_xM\ominus({\rm Ker}\,A_v\cap{\rm Ker}\,R(v))$ denotes $T_xM\cap({\rm Ker}\,A_v\cap{\rm Ker}\,R(v))^{\perp}$.  
However, if $G/K$ is of non-compact type, then this fact $(\ast_{\mathbb R})$ does not necessarily hold.  
For example, in the case where $G/K$ is a hyperbolic space of constant curvature $c(<0)$ and where $M$ is 
a hypersurface, $(\ast_{\mathbb R})$ holds if and only if 
all the absolute values of the principal curvatures of $M$ at each point are 
greater than $\sqrt{-c}$.  
So, in this paper, we consider the following condition:

\vspace{0.3truecm}

\hspace{0.1truecm}$(\ast_{\mathbb C})\,\,$
{\sl For any unit normal vector $v$ of $M$, the nullity spaces for 
complex focal radii }

\hspace{1truecm}{\sl along the normal geodesic $\gamma_v$ span 
$(T_xM)^{\mathbb C}\ominus({\rm Ker}\,A_v\cap{\rm Ker}\,R(v))^{\mathbb C}$.}

\vspace{0.3truecm}

\noindent
This condition $(\ast_{\mathbb C})$ is the condition weaker than $(\ast_{\mathbb R})$.  
In the case where $G/K$ is of non-comact type, 
$(\ast_{\mathbb C})$ also does not necessarily hold.  
For example, in the case where $G/K$ is a hyperbolic space of constant 
curvature $c(<0)$ and where $M$ is a hypersurface, 
$M$ satisfies $(\ast_{\mathbb C})$ if and only if all the principal curvatures of 
$M$ at each point of $M$ are not equal to $\pm\sqrt{-c}$.  

In this paper, we first prove the following result.  

\vspace{0.5truecm}

\noindent{\bf Theorem A.} {\sl Let $M$ be a complete isoparametric submanifold in a symmetric space $G/K$ of 
non-compact type or compact type.  
If $M$ admits a reflective focal submanifold, then it is curvature-adapted.}

\vspace{0.5truecm}

\noindent
{\it Remark 1.1.} In Theorem A, the condition of the existence of a reflective focal submanifold is indispensable.  
In fact, we have the following examples.  
Let $G=KAN$ be the Iwasawa's decomposition of $G$.  We can find many (complex) hyperpolar actions 
as subgroup actions of the solvable group $S:=AN$.  Since such hyperpolar actions admits no singular orbit, 
the principal orbits of the actions admits no focal submanifold.  Among such actions, we can find ones 
whose principal orbits are not curvature-adapted (see \cite{Koi9}).  

\vspace{0.5truecm}

Next we prove the following extrinsic homogeneity theorem.  

\vspace{0.5truecm}

\noindent{\bf Theorem B.} {\sl Let $M$ be a full irreducible complete curvature-adapted isoparametric 
$C^{\omega}$-submanifold of codimension greater than one in a symmetric space $G/K$ of non-compact type.  
If $M$ satisfies the above condition $(\ast_{\mathbb C})$, then $M$ is extrinsically homogeneous.}

\vspace{0.5truecm}

The proof of this theorem is performed by showing the extrinsic homogeneity of $M$ (see Theorem 7.1) 
by using Fact 1.1 through the anti-Kaehler submersion 
$\widetilde{\phi}:=\pi\circ\phi:H^0([0,1],\mathfrak g^{\mathbb C})\to G^{\mathbb C}/K^{\mathbb C}$ 
(see Figure 1).  

\vspace{10pt}

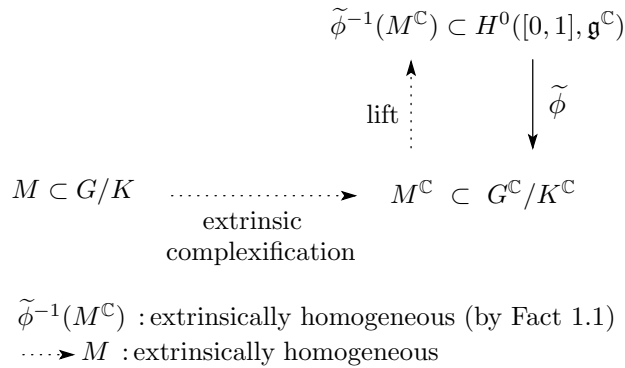
\begin{figure}[h]
\centerline{
\unitlength 0.1in
\begin{picture}( 76.1000, 18.1500)(-34.4000,-27.9500)
\put(20.3000,-18.7000){\makebox(0,0)[rt]{{\small $M\subset G/K$}}}%
%
\special{pn 8}%
\special{pa 2220 1960}%
\special{pa 3170 1960}%
\special{dt 0.045}%
\special{sh 1}%
\special{pa 3170 1960}%
\special{pa 3104 1940}%
\special{pa 3118 1960}%
\special{pa 3104 1980}%
\special{pa 3170 1960}%
\special{fp}%
\put(23.8000,-20.4000){\makebox(0,0)[lt]{{\small extrinsic}}}%
\put(21.9000,-22.0000){\makebox(0,0)[lt]{{\small complexification}}}%
\put(33.5000,-18.8000){\makebox(0,0)[lt]{{\small $M^{\mathbb C}\,\,\subset\,\,G^{\mathbb C}/K^{\mathbb C}$}}}%
%
\special{pn 8}%
\special{pa 4090 1250}%
\special{pa 4090 1710}%
\special{fp}%
\special{sh 1}%
\special{pa 4090 1710}%
\special{pa 4110 1644}%
\special{pa 4090 1658}%
\special{pa 4070 1644}%
\special{pa 4090 1710}%
\special{fp}%
\put(41.7000,-13.8000){\makebox(0,0)[lt]{$\widetilde{\phi}$}}%
\put(45.7000,-11.5000){\makebox(0,0)[rb]{{\small $\widetilde{\phi}^{-1}(M^{\mathbb C})\subset H^0([0,1],\mathfrak g^{\mathbb C})$}}}%
%
\special{pn 8}%
\special{pa 3460 1710}%
\special{pa 3460 1260}%
\special{dt 0.045}%
\special{sh 1}%
\special{pa 3460 1260}%
\special{pa 3440 1328}%
\special{pa 3460 1314}%
\special{pa 3480 1328}%
\special{pa 3460 1260}%
\special{fp}%
\put(34.0000,-15.8000){\makebox(0,0)[rb]{{\small lift}}}%
\put(14.2000,-25.0000){\makebox(0,0)[lt]{{\small $\widetilde{\phi}^{-1}(M^{\mathbb C})\,:\,$extrinsically homogeneous (by Fact 1.1)}}}%
%
\special{pn 8}%
\special{pa 1450 2790}%
\special{pa 1710 2790}%
\special{dt 0.045}%
\special{sh 1}%
\special{pa 1710 2790}%
\special{pa 1644 2770}%
\special{pa 1658 2790}%
\special{pa 1644 2810}%
\special{pa 1710 2790}%
\special{fp}%
\put(17.5000,-27.2000){\makebox(0,0)[lt]{{\small $M\,:\,$extrinsically homogeneous}}}%
\put(28.7000,-29.4000){\makebox(0,0)[lb]{$\,$}}%
\end{picture}%
\hspace{13.25truecm}}
\caption{The method of the proof of the extrinsic homogeneity}
\label{figure1}
\end{figure}

\vspace{10pt}

Let $G/K$ be a symmetric space of non-compact type and $H$ a closed subgroup of $G$.  
If there exists an involution $\sigma$ of $G$ with 
$({\rm Fix}\,\sigma)_0\subset H\subset{\rm Fix}\,\sigma$, then 
we (\cite{Koi4}) called the $H$-action on $G/K$ a {\it action of Hermann type}, where 
${\rm Fix}\,\sigma$ is the fixed point group of $\sigma$ and 
$({\rm Fix}\,\sigma)_0$ is the identity component of ${\rm Fix}\,\sigma$.  
In \cite{Koi10}, we called this kind of actions on semi-simple pseudo-Riemannian symmetric spaces (in more general) 
a {\it Hermann type action}.  In this paper, we shall use this terminology.  
According to the result in \cite{Koi4}, it follows that principal orbits of a Hermann type action are 
curvature-adapted complex equifocal (hence isoparametric) $C^{\omega}$-submanifolds and that they satisfy 
the condition $(\ast_{\mathbb C})$.  Also, a Hermann type action admits a reflective singular orbit and hence 
the principal orbits of the action admit a reflective focal submanifold.  

The main result of this paper is as follows.  

\vspace{0.5truecm}

\noindent{\bf Theorem C.} {\sl Let $M$ be a full irreducible complete isoparametric $C^{\omega}$-submanifold of 
codimension greater than one in a symmetric space $G/K$ of non-compact type.  
If $M$ admits a reflective focal submanifold, then it is a principal orbit of a Hermann type action 
on $G/K$.}

\vspace{0.5truecm}

If $M$ is a principal orbit of the isotropy action of a symmetric space $G/K$ of non-compact type, 
then it satisfies the following condition:

\vspace{0.3truecm}

\hspace{0.1truecm}$(\ast'_{\mathbb R})\,\,$
{\sl For any unit normal vector $v$ of $M$, the nullity spaces for focal radii along}

\hspace{1.15truecm}{\sl 
the normal geodesic $\gamma_v$ span $T_xM$.}

\vspace{0.5truecm}

By using the building theory, we prove that the following fact holds 
without the assumption of the real analyticity of the submanifold.  

\vspace{0.4truecm}

\noindent{\bf Theorem D.} {\sl Let $M$ be a full irreducible complete curvature-adapted 
isoparametric submanifold of codimension greater than two in an irreducible symmetric space $G/K$ of 
non-compact type.  
If $M$ satisfies the above condition $(\ast'_{\mathbb R})$, then $M$ is a principal orbit 
of the isotropy action of $G/K$.}

\vspace{0.5truecm}

In Section $2$-$5$, we shall recall the basic notions and facts.  
In Section 6, we shall prove Theorem A by using the basic facts stated in Section 5.  
In Section 7, we shall prove Theorem B by using Fact 1.1.  
In Section 8, we shall prove Theorems C (Main theorem) by using Theorems A and B.  
In Section 9, we shall classify isoparametric submanifolds as in Theorem C.  
In Section 10, we prove Theorem D.  

\section{Complex focal radius} 
In this section, we shall recall the notion of a complex focal radius and some facts related to it, 
which will be used in Sections 7 and 8.  
Let $M$ be a submanifold in a complete Riemannian manifold $N$, 
$\psi\,:\,T^{\perp}M\to M$ the normal bundle of $M$ and $\exp^{\perp}$ the 
normal exponential map of $M$.  Denote by ${\mathcal V}$ the vertical distribution 
on $T^{\perp}M$ and 
${\mathcal H}$ the horizontal distribution on $T^{\perp}M$ with respect to the 
normal connection of $M$.  
Let $v$ be a unit normal vector of $M$ at $x(\in M)$ and $r$ 
a real number.  Denote by $\gamma_v$ the normal geodesic of $M$ of direction 
$v$ (i.e., $\gamma_v(s)=\exp^{\perp}(sv)$).  
If $\psi_{\ast}({\rm Ker}\,\exp^{\perp}_{\ast rv})\not=\{0\}$, then 
$\exp^{\perp}(rv)$ (resp. $r$) is called a {\it focal point} 
(resp. a {\it focal radius}) {\it of} $M$ {\it along} $\gamma_v$.  
For a focal radius $r$ of $M$ along $\gamma_v$, 
$\psi_{\ast}({\rm Ker}\,\exp^{\perp}_{\ast rv})$ is called the 
{\it nullity space} for $r$ and its dimension is called 
the {\it multiplicity} of $r$.  
Denote by ${\mathcal FR}^{\mathbb R}_{M,v}$ the set of all focal radii of $M$ along 
$\gamma_v$.  Set 
$${\mathcal F}^{\mathbb R}_{M,x}:=\mathop{\bigcup}_{v\in T^{\perp}_xM\,\,{\rm s.t.}\,\,
|| v||=1}\{rv\,|\,r\in{\mathcal FR}^{\mathbb R}_{M,v}\},$$
which is called the {\it tangential focal set of} $M$ {\it at} $x$.  Note that 
$\exp^{\perp}({\mathcal F}^{\mathbb R}_{M,x})$ is the focal set of $M$ at $x$.  
If, for any $y\in M$, the normal umbrella $\Sigma_y:=\exp^{\perp}(T^{\perp}_yM)$ is totally 
geodesic in $G/K$ and the induced metric on $\Sigma_y$ is flat, then $M$ is 
called a {\it submanifold with flat section}.  
Assume that $N$ is a symmetric space $G/K$ and that $M$ is a submanifold with flat section.  Then we can show that, 
for any $rv\in T^{\perp}M$, ${\rm Ker}\,\exp^{\perp}_{\ast rv}\subset{\mathcal H}_{rv}$ holds and that 
$$\exp^{\perp}_{\ast rv}(X_{rv}^L)=P_{\gamma_{rv}\vert_{[0,1]}}
\left(\left(\cos(r\sqrt{R(v)})-\frac{\sin(r\sqrt{R(v)})}{\sqrt{R(v)}}\circ 
A_v\right)(X)\right)\,\,\,\,(X\in T_xM)\leqno{(2.1)}$$
holds, where $X_{rv}^L$ is the horizontal lift of $X$ to $rv$, $P_{\gamma_{rv}\vert_{[0,1]}}$ is 
the parallel translation along the normal geodesic $\gamma_{rv}\vert_{[0,1]}$, $R(v)$ is 
the normal Jacobi operator $R(\bullet,v)v$ and $A$ is the shape tensor of $M$ (see Figure 2).  
Hence ${\mathcal FR}^{\mathbb R}_{M,v}$ coincide with the set of all zero points of 
the real-valued function 
$$F_v(s)={\rm det}\left(\cos(s\sqrt{R(v)})
-\frac{\sin(s\sqrt{R(v)})}{\sqrt{R(v)}}\circ A_v\right)\quad\,\,
(s\in{\mathbb R}).$$

\vspace{10pt}

\begin{figure}[h]
\centerline{
\unitlength 0.1in
\begin{picture}( 34.2000, 26.5000)( -2.4000,-34.4000)
%
\special{pn 8}%
\special{pa 1130 2554}%
\special{pa 1138 2522}%
\special{pa 1148 2492}%
\special{pa 1160 2462}%
\special{pa 1174 2434}%
\special{pa 1190 2406}%
\special{pa 1208 2380}%
\special{pa 1228 2354}%
\special{pa 1246 2328}%
\special{pa 1270 2306}%
\special{pa 1292 2284}%
\special{pa 1316 2264}%
\special{pa 1342 2244}%
\special{pa 1368 2226}%
\special{pa 1396 2210}%
\special{pa 1424 2194}%
\special{pa 1454 2180}%
\special{pa 1482 2166}%
\special{pa 1512 2156}%
\special{pa 1542 2146}%
\special{pa 1574 2136}%
\special{pa 1604 2128}%
\special{pa 1636 2122}%
\special{pa 1654 2120}%
\special{sp}%
%
\special{pn 8}%
\special{pa 2144 2702}%
\special{pa 2152 2672}%
\special{pa 2162 2642}%
\special{pa 2172 2612}%
\special{pa 2188 2582}%
\special{pa 2202 2554}%
\special{pa 2220 2528}%
\special{pa 2240 2502}%
\special{pa 2260 2478}%
\special{pa 2282 2454}%
\special{pa 2306 2432}%
\special{pa 2330 2412}%
\special{pa 2356 2392}%
\special{pa 2382 2376}%
\special{pa 2410 2358}%
\special{pa 2438 2344}%
\special{pa 2466 2328}%
\special{pa 2496 2316}%
\special{pa 2526 2304}%
\special{pa 2556 2294}%
\special{pa 2586 2286}%
\special{pa 2618 2278}%
\special{pa 2648 2270}%
\special{pa 2668 2268}%
\special{sp}%
%
\special{pn 8}%
\special{pa 1138 2562}%
\special{pa 1170 2562}%
\special{pa 1202 2566}%
\special{pa 1234 2568}%
\special{pa 1266 2570}%
\special{pa 1298 2572}%
\special{pa 1328 2576}%
\special{pa 1360 2578}%
\special{pa 1392 2580}%
\special{pa 1424 2584}%
\special{pa 1456 2586}%
\special{pa 1488 2590}%
\special{pa 1520 2594}%
\special{pa 1552 2596}%
\special{pa 1584 2600}%
\special{pa 1616 2604}%
\special{pa 1648 2608}%
\special{pa 1678 2612}%
\special{pa 1710 2616}%
\special{pa 1742 2622}%
\special{pa 1774 2626}%
\special{pa 1806 2630}%
\special{pa 1836 2636}%
\special{pa 1868 2640}%
\special{pa 1900 2648}%
\special{pa 1932 2652}%
\special{pa 1964 2658}%
\special{pa 1994 2664}%
\special{pa 2026 2670}%
\special{pa 2058 2676}%
\special{pa 2088 2684}%
\special{pa 2120 2690}%
\special{pa 2144 2696}%
\special{sp}%
%
\special{pn 8}%
\special{pa 1654 2130}%
\special{pa 1686 2132}%
\special{pa 1718 2134}%
\special{pa 1750 2136}%
\special{pa 1782 2138}%
\special{pa 1814 2140}%
\special{pa 1844 2144}%
\special{pa 1876 2146}%
\special{pa 1908 2150}%
\special{pa 1940 2152}%
\special{pa 1972 2156}%
\special{pa 2004 2160}%
\special{pa 2036 2162}%
\special{pa 2068 2166}%
\special{pa 2100 2170}%
\special{pa 2132 2174}%
\special{pa 2164 2176}%
\special{pa 2194 2182}%
\special{pa 2226 2186}%
\special{pa 2258 2190}%
\special{pa 2290 2196}%
\special{pa 2322 2200}%
\special{pa 2352 2206}%
\special{pa 2384 2210}%
\special{pa 2416 2216}%
\special{pa 2448 2222}%
\special{pa 2480 2228}%
\special{pa 2510 2232}%
\special{pa 2542 2240}%
\special{pa 2574 2246}%
\special{pa 2604 2252}%
\special{pa 2636 2260}%
\special{pa 2660 2266}%
\special{sp}%
%
\special{pn 8}%
\special{pa 1776 2628}%
\special{pa 1748 2612}%
\special{pa 1722 2596}%
\special{pa 1694 2578}%
\special{pa 1668 2558}%
\special{pa 1644 2538}%
\special{pa 1620 2518}%
\special{pa 1598 2494}%
\special{pa 1580 2468}%
\special{pa 1564 2440}%
\special{pa 1552 2410}%
\special{pa 1546 2378}%
\special{pa 1550 2346}%
\special{pa 1558 2316}%
\special{pa 1574 2288}%
\special{pa 1594 2264}%
\special{pa 1616 2240}%
\special{pa 1640 2220}%
\special{pa 1666 2200}%
\special{pa 1694 2184}%
\special{pa 1720 2168}%
\special{pa 1750 2154}%
\special{pa 1778 2140}%
\special{pa 1782 2138}%
\special{sp}%
%
\special{pn 20}%
\special{sh 1}%
\special{ar 1546 2374 10 10 0  6.28318530717959E+0000}%
\special{sh 1}%
\special{ar 1546 2374 10 10 0  6.28318530717959E+0000}%
%
\special{pn 8}%
\special{pa 1628 2520}%
\special{pa 2068 2374}%
\special{dt 0.045}%
%
\special{pn 20}%
\special{sh 1}%
\special{ar 2038 2372 10 10 0  6.28318530717959E+0000}%
\special{sh 1}%
\special{ar 2038 2372 10 10 0  6.28318530717959E+0000}%
%
\special{pn 8}%
\special{pa 1540 1390}%
\special{pa 1536 2332}%
\special{dt 0.045}%
%
\special{pn 8}%
\special{pa 1250 1650}%
\special{pa 1760 1180}%
\special{fp}%
%
\special{pn 8}%
\special{pa 1390 1190}%
\special{pa 850 1650}%
\special{pa 2290 1650}%
\special{pa 2780 1190}%
\special{pa 2780 1190}%
\special{pa 1390 1190}%
\special{fp}%
%
\special{pn 20}%
\special{sh 1}%
\special{ar 1540 1380 10 10 0  6.28318530717959E+0000}%
\special{sh 1}%
\special{ar 1540 1380 10 10 0  6.28318530717959E+0000}%
%
\special{pn 13}%
\special{pa 1538 2354}%
\special{pa 1638 2180}%
\special{fp}%
\special{sh 1}%
\special{pa 1638 2180}%
\special{pa 1588 2228}%
\special{pa 1612 2226}%
\special{pa 1622 2248}%
\special{pa 1638 2180}%
\special{fp}%
%
\special{pn 8}%
\special{pa 1170 1380}%
\special{pa 2580 1380}%
\special{fp}%
%
\special{pn 8}%
\special{pa 2030 2360}%
\special{pa 2030 1380}%
\special{dt 0.045}%
%
\special{pn 20}%
\special{sh 1}%
\special{ar 2030 1380 10 10 0  6.28318530717959E+0000}%
\special{sh 1}%
\special{ar 2030 1380 10 10 0  6.28318530717959E+0000}%
%
\special{pn 13}%
\special{pa 2030 1380}%
\special{pa 2170 1220}%
\special{fp}%
\special{sh 1}%
\special{pa 2170 1220}%
\special{pa 2112 1258}%
\special{pa 2136 1260}%
\special{pa 2142 1284}%
\special{pa 2170 1220}%
\special{fp}%
%
\special{pn 8}%
\special{pa 1730 1220}%
\special{pa 2320 1220}%
\special{dt 0.045}%
%
\special{pn 8}%
\special{pa 1660 1290}%
\special{pa 2100 1290}%
\special{dt 0.045}%
%
\special{pn 8}%
\special{pa 2080 1330}%
\special{pa 1890 1530}%
\special{dt 0.045}%
%
\special{pn 8}%
\special{pa 1410 1520}%
\special{pa 1890 1520}%
\special{dt 0.045}%
%
\special{pn 8}%
\special{pa 1490 1450}%
\special{pa 1970 1450}%
\special{dt 0.045}%
%
\special{pn 8}%
\special{pa 1580 2460}%
\special{pa 2030 2370}%
\special{dt 0.045}%
%
\special{pn 8}%
\special{pa 1460 1020}%
\special{pa 1680 1250}%
\special{dt 0.045}%
\special{sh 1}%
\special{pa 1680 1250}%
\special{pa 1648 1188}%
\special{pa 1644 1212}%
\special{pa 1620 1216}%
\special{pa 1680 1250}%
\special{fp}%
%
\special{pn 8}%
\special{pa 1900 1020}%
\special{pa 2100 1300}%
\special{dt 0.045}%
\special{sh 1}%
\special{pa 2100 1300}%
\special{pa 2078 1234}%
\special{pa 2070 1258}%
\special{pa 2046 1258}%
\special{pa 2100 1300}%
\special{fp}%
%
\special{pn 8}%
\special{pa 1210 2050}%
\special{pa 1600 2260}%
\special{dt 0.045}%
\special{sh 1}%
\special{pa 1600 2260}%
\special{pa 1552 2212}%
\special{pa 1554 2236}%
\special{pa 1532 2246}%
\special{pa 1600 2260}%
\special{fp}%
%
\special{pn 8}%
\special{pa 1240 2780}%
\special{pa 1600 2500}%
\special{dt 0.045}%
\special{sh 1}%
\special{pa 1600 2500}%
\special{pa 1536 2526}%
\special{pa 1558 2534}%
\special{pa 1560 2558}%
\special{pa 1600 2500}%
\special{fp}%
%
\special{pn 8}%
\special{pa 2690 1010}%
\special{pa 2380 1380}%
\special{dt 0.045}%
\special{sh 1}%
\special{pa 2380 1380}%
\special{pa 2438 1342}%
\special{pa 2414 1340}%
\special{pa 2408 1316}%
\special{pa 2380 1380}%
\special{fp}%
\put(12.2000,-28.3000){\makebox(0,0)[rt]{$M$}}%
\put(11.6000,-20.6000){\makebox(0,0)[rb]{$X$}}%
\put(21.1000,-19.2000){\makebox(0,0)[lb]{$\gamma_{rv}(1)$}}%
\put(31.8000,-22.2000){\makebox(0,0)[lt]{$\gamma_{rv}$}}%
\put(15.6000,-9.7000){\makebox(0,0)[rb]{{\small $0$-section}}}%
\put(19.4000,-9.7000){\makebox(0,0)[rb]{$X_{rv}^L$}}%
%
\special{pn 8}%
\special{pa 2550 1590}%
\special{pa 2040 1390}%
\special{dt 0.045}%
\special{sh 1}%
\special{pa 2040 1390}%
\special{pa 2096 1434}%
\special{pa 2090 1410}%
\special{pa 2110 1396}%
\special{pa 2040 1390}%
\special{fp}%
\put(26.1000,-15.7000){\makebox(0,0)[lt]{$rv$}}%
\put(27.0000,-9.6000){\makebox(0,0)[lb]{$T_x^{\perp}M$}}%
\put(14.9000,-23.2000){\makebox(0,0)[rt]{$x$}}%
\put(26.5000,-24.4000){\makebox(0,0)[lt]{$N$}}%
\put(29.7000,-13.4000){\makebox(0,0)[lt]{$T^{\perp}M$}}%
%
\special{pn 13}%
\special{pa 1606 1790}%
\special{pa 1606 2000}%
\special{fp}%
\special{sh 1}%
\special{pa 1606 2000}%
\special{pa 1626 1934}%
\special{pa 1606 1948}%
\special{pa 1586 1934}%
\special{pa 1606 2000}%
\special{fp}%
\put(16.6500,-18.3000){\makebox(0,0)[lt]{$\exp^{\perp}$}}%
\put(14.1000,-34.4000){\makebox(0,0)[lt]{$\,$}}%
%
\special{pn 13}%
\special{pa 2030 1380}%
\special{pa 2320 1220}%
\special{fp}%
\special{sh 1}%
\special{pa 2320 1220}%
\special{pa 2252 1236}%
\special{pa 2274 1246}%
\special{pa 2272 1270}%
\special{pa 2320 1220}%
\special{fp}%
%
\special{pn 8}%
\special{pa 1680 2200}%
\special{pa 2030 2370}%
\special{dt 0.045}%
%
\special{pn 8}%
\special{pa 1630 2260}%
\special{pa 2030 2370}%
\special{dt 0.045}%
%
\special{pn 13}%
\special{pa 2050 2370}%
\special{pa 2190 2370}%
\special{fp}%
\special{sh 1}%
\special{pa 2190 2370}%
\special{pa 2124 2350}%
\special{pa 2138 2370}%
\special{pa 2124 2390}%
\special{pa 2190 2370}%
\special{fp}%
%
\special{pn 8}%
\special{pa 2220 1930}%
\special{pa 2040 2360}%
\special{dt 0.045}%
\special{sh 1}%
\special{pa 2040 2360}%
\special{pa 2084 2306}%
\special{pa 2062 2312}%
\special{pa 2048 2292}%
\special{pa 2040 2360}%
\special{fp}%
%
\special{pn 8}%
\special{pa 2170 1000}%
\special{pa 2230 1260}%
\special{dt 0.045}%
\special{sh 1}%
\special{pa 2230 1260}%
\special{pa 2234 1192}%
\special{pa 2218 1208}%
\special{pa 2196 1200}%
\special{pa 2230 1260}%
\special{fp}%
\put(22.0000,-9.6000){\makebox(0,0)[rb]{$Y$}}%
%
\special{pn 8}%
\special{ar 1780 2700 980 340  4.4728615 5.2974740}%
\put(24.4000,-21.6000){\makebox(0,0)[lb]{$\exp^{\perp}_{\ast\,rv}(Y)$}}%
%
\special{pn 8}%
\special{ar 2420 2380 320 290  3.2394930 3.2788372}%
\special{ar 2420 2380 320 290  3.3968700 3.4362143}%
\special{ar 2420 2380 320 290  3.5542471 3.5935913}%
\special{ar 2420 2380 320 290  3.7116241 3.7509684}%
\special{ar 2420 2380 320 290  3.8690012 3.9083454}%
\special{ar 2420 2380 320 290  4.0263782 4.0657225}%
\special{ar 2420 2380 320 290  4.1837553 4.2230995}%
\special{ar 2420 2380 320 290  4.3411323 4.3804766}%
\special{ar 2420 2380 320 290  4.4985094 4.5378536}%
\special{ar 2420 2380 320 290  4.6558864 4.6946289}%
%
\special{pn 8}%
\special{pa 2110 2320}%
\special{pa 2090 2360}%
\special{fp}%
\special{sh 1}%
\special{pa 2090 2360}%
\special{pa 2138 2310}%
\special{pa 2114 2312}%
\special{pa 2102 2292}%
\special{pa 2090 2360}%
\special{fp}%
\put(8.8000,-32.7000){\makebox(0,0)[lt]{$\exp^{\perp}_{\ast\,rv}(X_{rv}^L)=0,\,\,\exp^{\perp}_{\ast\,rv}(Y)\not=0$}}%
\put(11.2000,-32.4000){\makebox(0,0)[lb]{$\psi_{\ast\,rv}(X_{rv}^L)=\psi_{\ast\,rv}(Y)=X$}}%
%
\special{pn 8}%
\special{ar 2790 2420 540 200  3.2477708 3.2802032}%
\special{ar 2790 2420 540 200  3.3775005 3.4099329}%
\special{ar 2790 2420 540 200  3.5072302 3.5396627}%
\special{ar 2790 2420 540 200  3.6369600 3.6693924}%
\special{ar 2790 2420 540 200  3.7666897 3.7991221}%
\special{ar 2790 2420 540 200  3.8964194 3.9288519}%
\special{ar 2790 2420 540 200  4.0261492 4.0585816}%
\special{ar 2790 2420 540 200  4.1558789 4.1883113}%
\special{ar 2790 2420 540 200  4.2856086 4.3180411}%
\special{ar 2790 2420 540 200  4.4153384 4.4477708}%
\special{ar 2790 2420 540 200  4.5450681 4.5775005}%
\special{ar 2790 2420 540 200  4.6747978 4.7072302}%
\special{ar 2790 2420 540 200  4.8045275 4.8369600}%
\special{ar 2790 2420 540 200  4.9342573 4.9666897}%
\special{ar 2790 2420 540 200  5.0639870 5.0964194}%
\special{ar 2790 2420 540 200  5.1937167 5.2261492}%
\special{ar 2790 2420 540 200  5.3234465 5.3558789}%
\special{ar 2790 2420 540 200  5.4531762 5.4748600}%
%
\special{pn 8}%
\special{pa 2270 2370}%
\special{pa 2260 2400}%
\special{fp}%
\special{sh 1}%
\special{pa 2260 2400}%
\special{pa 2300 2344}%
\special{pa 2278 2350}%
\special{pa 2262 2330}%
\special{pa 2260 2400}%
\special{fp}%
\end{picture}%
\hspace{2.5truecm}}
\caption{Focal points of a submanifold with flat section}
\label{figure2}
\end{figure}
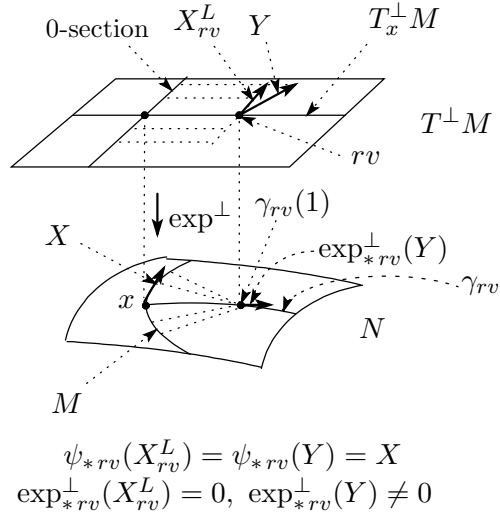

\vspace{10pt}

In particular, in the case where $G/K$ is a Euclidean space, we have 
$F_v(s)={\rm det}({\rm id}-sA_v)$ (${\rm id}\,:\,$the identity transformation 
of $T_xM$).  Hence 
${\mathcal FR}^{\mathbb R}_{M,v}$ is equal to the set of all the inverse numbers of 
the eigenvalues of $A_v$ and the nullity space for 
$r\in{\mathcal FR}^{\mathbb R}_{M,v}$ is equal 
to ${\rm Ker}(A_v-\frac 1r{\rm id})$.  Therefore the nullity spaces for focal 
radii of $M$ along $\gamma_v$ span $T_xM\ominus{\rm Ker}\,A_v$.  
In the case where $G/K$ is a sphere of constant curvature $c(>0)$, 
we have 
$$F_v(s)={\rm det}\left(\cos(s\sqrt c){\rm id}
-\frac{\sin(s\sqrt c)}{\sqrt c}A_v\right).$$
Hence we have 
$${\mathcal FR}^{\mathbb R}_{M,v}=\left\{\left.\frac{1}{\sqrt c}
\left(\arctan\frac{\sqrt c}{\lambda}+j\pi\right)\,\right|\,\lambda\,:\,
{\rm the}\,\,{\rm eigenvalue}\,\,{\rm of}\,\,A_v,\,\,\,\,j\in{\mathbb Z}\right\}$$
and the nullity space for 
$\displaystyle{\frac{1}{\sqrt c}\left(\arctan\frac{\sqrt c}{\lambda}
+j\pi\right)}$ is equal to ${\rm Ker}(A_v-\lambda\,{\rm id})$, where 
we note that $\arctan\frac{\sqrt c}{\lambda}$ means $\frac{\pi}{2}$ when 
$\lambda=0$.  Therefore the nullity spaces for focal radii of $M$ along $\gamma_v$ span $T_xM$.  
Note that the focal set of $M$ at $x$ is given by 
$${\mathcal F}_{M,x}^{\mathbb R}=\mathop{\bigcup}_{v\in T^{\perp}_xM\,\,{\rm s.t.}\,\,||v||=1}
\left\{\gamma_v(r)\,\left|\,r=\frac{1}{\sqrt c}\arctan\frac{\sqrt c}{\lambda}\,\,{\rm or}\,\,
\frac{1}{\sqrt c}\left(\arctan\frac{\sqrt c}{\lambda};\pi\right)\right.\right\}.$$
In the case where $G/K$ is a hyperbolic space of constant curvature $c(<0)$, we have 
$$F_v(s)={\rm det}\left(\cosh(s\sqrt{-c}){\rm id}-\frac{\sinh(s\sqrt{-c})}
{\sqrt{-c}}A_v\right).$$
Hence we have 
$${\mathcal FR}^{\mathbb R}_{M,v}=\left\{\left.\frac{1}{\sqrt{-c}}{\rm arctanh}
\frac{\sqrt{-c}}{\lambda}\,\right|\,\lambda\,:\,{\rm the}\,\,
{\rm eigenvalue}\,\,{\rm of}\,\,A_v\,\,{\rm s.t.}\,\,
|\lambda|>\sqrt{-c}\right\}\leqno{(2.2)}$$
and the nullity space for $\displaystyle{\frac{1}{\sqrt{-c}}
{\rm arctanh}\frac{\sqrt{-c}}{\lambda}}$ is equal to 
${\rm Ker}(A_v-\lambda\,{\rm id})$.  
Therefore the nullity spaces for focal radii of $M$ along $\gamma_v$ span 
$T_xM$ if and only if all the absolute values of eigenvalues of $A_v$ is 
greater than $\sqrt{-c}$.  
As a non-compact submanifold $M$ with flat section in a symmetric space $G/K$ 
of non-compact type is deformed as its principal curvatures approach to zero, 
its focal set vanishes beyond the ideal boundary $(G/K)(\infty)$ of 
$G/K$.  This fact follows from $(2.2)$.  According to this fact, we (\cite{Koi2}) considered that a focal radius 
of $M$ along the normal geodesic $\gamma_v$ should be defined in the complex number field ${\mathbb C}$.  
We (\cite{Koi2}) introduced the notion of a complex focal radius as the zero points 
of the complex-valued function $F^{\mathbb C}_v$ over ${\mathbb C}$ defined by 
$$F^{\mathbb C}_v(z):={\rm det}\left(\cos(z{\sqrt{R(v)}}^{\mathbb C})
-\frac{\sin(z{\sqrt{R(v)}}^{\mathbb C})}
{{\sqrt{R(v)}}^{\mathbb C}}\circ A^{\mathbb C}_v\right)\quad\,\,(z\in{\mathbb C}),$$
where $A^{\mathbb C}_v$ and ${\sqrt{R(v)}}^{\mathbb C}$ are the complexifications of 
$A_v$ and $\sqrt{R(v)}$, respectively.  
For a complex focal radius $z$ of $M$ along $\gamma_v$, 
$\displaystyle{{\rm Ker}\left(\cos(z{\sqrt{R(v)}}^{\mathbb C})
-\frac{\sin(z{\sqrt{R(v)}}^{\mathbb C})}
{{\sqrt{R(v)}}^{\mathbb C}}\circ A^{\mathbb C}_v\right)\,\,\,(\subset(T_xM)^{\mathbb C})}$ 
is called the {\it nullity space} for $z$ and its complex dimension is called 
the {\it multiplicity} of $z$.  
Denote by ${\mathcal FR}^{\mathbb C}_{M,v}$ the set of all complex focal radii of $M$ along $\gamma_v$.  
Set 
$${\mathcal F}^{\mathbb C}_{M,x}:=\mathop{\bigcup}_{v\in T^{\perp}_xM\,\,{\rm s.t.}\,\,
||v||=1}\{rv\,|\,r\in{\mathcal FR}^{\mathbb C}_{M,v}\}\,\,\,\,(\subset(T^{\perp}_xM)^{\mathbb C}),$$
which is called the {\it tangential complex focal set of} $M$ {\it at} $x$.  
In the case where $G/K$ is a Euclidean space, we have 
$F_v^{\mathbb C}(z)={\rm det}({\rm id}-zA^{\mathbb C}_v)$ (${\rm id}\,:\,$ 
the identity transformation of $(T_xM)^{\mathbb C}$).  Hence 
we have ${\mathcal FR}^{\mathbb C}_{M,v}={\mathcal FR}^{\mathbb R}_{M,v}$ and the nullity 
space for $z\in{\mathcal FR}^{\mathbb C}_{M,v}$ ie equal to ${\rm Ker}
(A_v^{\mathbb C}-\frac 1z{\rm id})$.  Therefore the nullity spaces for 
complex focal radii of $M$ along $\gamma_v$ span $(T_xM)^{\mathbb C}\ominus
{\rm Ker}\,A_v^{\mathbb C}$.  
Also, in the case where $G/K$ is a sphere of constant curvature $c(>0)$, 
we have 
$$F^{\mathbb C}_v(z)={\rm det}\left(\cos(z\sqrt c){\rm id}
-\frac{\sin(z\sqrt c)}{\sqrt c}A^{\mathbb C}_v\right).$$
Hence ${\mathcal FR}^{\mathbb C}_{M,v}={\mathcal FR}^{\mathbb R}_{M,v}$ and the nullity space 
for $\displaystyle{\frac{1}{\sqrt c}\left(\arctan\frac{\sqrt c}{\lambda}
+j\pi\right)}$ is equal to ${\rm Ker}(A^{\mathbb C}_v-\lambda\,{\rm id})$.  
Therefore the nullity spaces for complex focal radii of $M$ along $\gamma_v$ 
span $(T_xM)^{\mathbb C}$.  
Also, in the case where $G/K$ is a hyperbolic space of constant curvature 
$c(<0)$, we have 
$$F_v^{\mathbb C}(z)={\rm det}\left(\cos(\sqrt{-1}z\sqrt{-c}){\rm id}
-\frac{\sin(\sqrt{-1}z\sqrt{-c})}{\sqrt{-1}\sqrt{-c}}A_v^{\mathbb C}\right),$$
where $\sqrt{-1}$ denotes the imaginary unit.  
Hence ${\mathcal FR}^{\mathbb C}_{M,v}$ is equal to 
$$
%
\hspace{8truecm}}
\caption{The geometrical meaning of the complex focal radius}
\label{figure3}
\end{figure}

\vspace{10pt}

\begin{figure}[h]
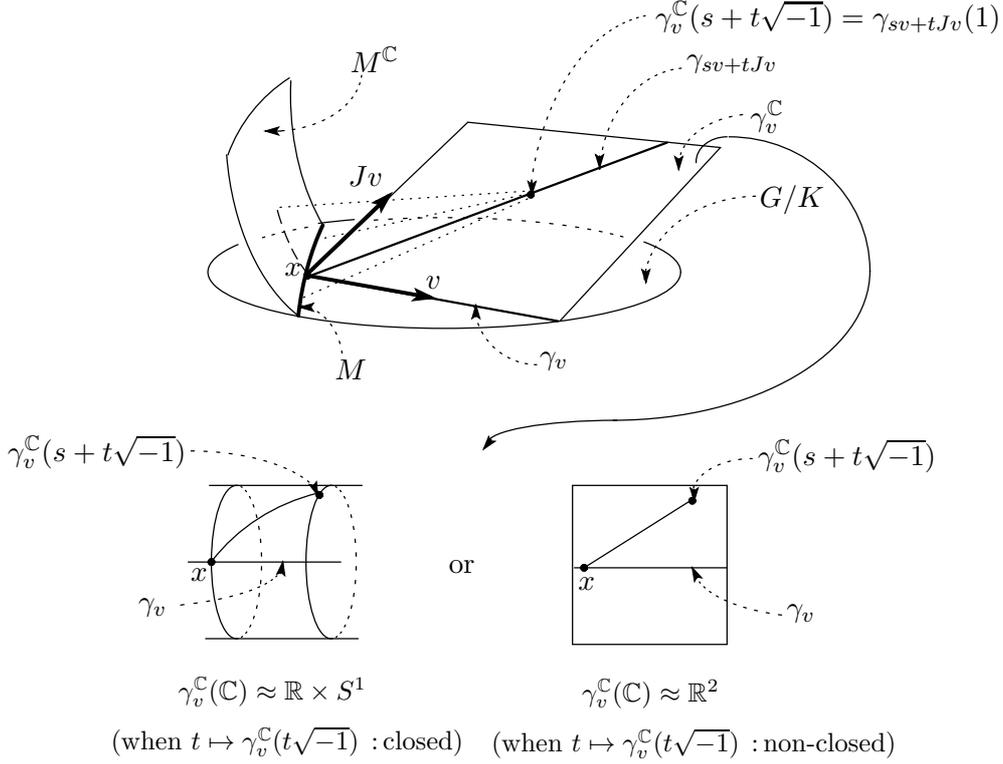

\centerline{
\unitlength 0.1in
%
%
\hspace{5.5truecm}}
\caption{The geometrical meaning of the complex focal radius (continued)}
\label{figure4}
\end{figure}

\vspace{10pt}


\noindent
the nullity space for 
$\displaystyle{\frac{1}{\sqrt{-c}}\left({\rm arctanh}\frac{\sqrt{-c}}{\lambda}
+j\pi\sqrt{-1}\right)}$ ($|\lambda|>\sqrt{-c}$) is equal to 
${\rm Ker}(A^{\mathbb C}_v-\lambda\,{\rm id})$ 
and the nullity space for 
$\displaystyle{\frac{1}{\sqrt{-c}}\left({\rm arctanh}\frac{\lambda}
{\sqrt{-c}}+(j+\frac12)\pi\sqrt{-1}\right)}$ ($|\lambda|>\sqrt{-c}$) is 
equal to ${\rm Ker}(A^{\mathbb C}_v-\lambda\,{\rm id})$.  
Therefore the nullity spaces for complex focal radii of $M$ along $\gamma_v$ 
span $(T_xM)^{\mathbb C}$ if and only if all the eigenvalues of $A_v$ are not 
equal to $\pm\sqrt{-c}$.  

Let $M$ be a $C^{\omega}$-submanifold with flat section in a symmetric space 
$G/K$.  Then we can define the complete complexification $M^{\mathbb C}$ of $M$ 
as a (metrically) complete anti-Kaehler submanifold in the anti-Kaehler symmetric space 
$G^{\mathbb C}/K^{\mathbb C}$ associated with $G/K$ (see the proof of Theorem B in \cite{Koi7}).  
Denote by $J$ and $\widehat R$ the complex structure and the curvature tensor 
of $G^{\mathbb C}/K^{\mathbb C}$, respectively, and 
$\widehat A$ and $\widehat{\exp}^{\perp}$ the shape tensor and 
the normal exponential map of $M^{\mathbb C}$, repsectively.  
Denote by $\widehat{\mathcal H}$ the horizontal distribution on the normal bundle 
$T^{\perp}(M^{\mathbb C})$ of $M^{\mathbb C}$ with respect to the normal connection of 
$M^{\mathbb C}$.  
Take $v\in T^{\perp}_xM(\subset T_x^{\perp}(M^{\mathbb C}))$ and $z=s+t\sqrt{-1}
\in{\mathbb C}$ ($s,t\in{\mathbb R}$).  
Then we can show ${\rm Ker}\,\widehat{\exp}^{\perp}_{\ast sv+tJv}\subset
\widehat{{\mathcal H}}_{sv+tJv}$ and 
$$\widehat{\exp}^{\perp}_{\ast sv+tJv}(X_{sv+tJv}^L)
=P_{\gamma_{sv+tJv}\vert_{[0,1]}}(Q_{v,z}(X))\,\,\,\,\,(X\in T_x(M^{\mathbb C})),$$
where $X_{s+tJv}^L$ is the horizontal lift of $X$ to $sv+tJv$, 
$P_{\gamma_{sv+tJv}}$ is the parallel translation along the normal 
geodesic $\gamma_{sv+tJv}$ of $M^{\mathbb C}$ and 
$$Q_{v,z}:=
\cos\left(s\sqrt{\widehat R(v)}+t\left(J\circ\sqrt{\widehat R(v)}\right)\right)
-\frac{\sin\left(s\sqrt{\widehat R(v)}+t\left(J\circ\sqrt{\widehat R(v)}\right)
\right)}{\sqrt{\widehat R(v)}}\circ\widehat A_v$$
($\widehat R(v):=\widehat R(\bullet,v)v$).  
Hence $\widehat{\exp}^{\perp}(sv+tJv)=\gamma_{sv+tJv}(1)=\gamma^{\mathbb C}_v(s+t\sqrt{-1})$ is a focal point of 
$M^{\mathbb C}$ along the geodesic $\gamma_{sv+tJv}$ if and only if 
$z=s+t\sqrt{-1}$ is a zero point of the complex-valued function $\widehat F_v$ 
over ${\mathbb C}$ defined by $\widehat F_v(z):={\rm det}\,Q_{v,z}$, 
where $Q_{v,z}$ is regarded as a ${\mathbb C}$-linear transformation of 
$T_x(M^{\mathbb C})$ regarded as a complex linear space by $J$.  
On the other hand, it is clear that 
the set of all zero points of $\widehat F_v$ is equal to that of 
$F_v^{\mathbb C}$.  Therefore $z=s+t\sqrt{-1}$ is a complex focal radius of $M$ 
along $\gamma_v$ if and only if $\gamma_v^{\mathbb C}(s+t\sqrt{-1})=\widehat{\exp}^{\perp}(sv+tJv)$ is 
a focal point of $M^{\mathbb C}$ along $\gamma_{sv+tJv}$ (see Figures 3 and 4).  Hence we see that 
$\widehat{\exp}^{\perp}({\mathcal F}^{\mathbb C}_{M,x})$ is the focal set of 
$M^{\mathbb C}\,(\subset G^{\mathbb C}/K^{\mathbb C})$ at $x$, where we identify $(T_xM)^{\mathbb C}$ with 
$T_x(M^{\mathbb C})$.  
Thus we can grasp the geometrical meaning of the complex focal radius.  

\section{Anti-Kaehler isoparametric submanifold}
In this section, we shall recall the notion of a proper anti-Kaehler isoparametric submanifold, which was 
introduced in \cite{Koi3}.  
Throughout this paper, we shall call this notion an {\it anti-Kaehler isoparametric submanifold with 
$J$-diagonalizable shape operators} because this terminology seems to be more familiar than 
``proper anti-Kaehler isoparametric submanifold''.  
Also, we shall state the defintion of the (metrically) completeness of 
an anti-Kaehler Hilbert manifold.   Furthermore we shall recall some facts related to this submanifold, 
which will be used in Sections 7 and 8.  

First we shall recall the notion of an infinite dimensional anti-Kaehler space, where we shall define this 
notion more smartly than the definition in \cite{Koi3}.  
Let $V$ be an infinite dimensional topological real vector space, $J$ be a continuous 
linear operator of $V$ such that $J^2=-{{\rm id}}$ and 
$\langle\,\,,\,\,\rangle$ be a continuous non-degenerate symmetric 
bilinear form of $V$ such that $\langle JX,JY\rangle=-\langle X,Y\rangle$ 
holds for every $X,Y\in V$.  
It is easy to show that there uniquely exists an orthogonal time-space decomposition 
$V=V_-\oplus V_+$ compatible with $J$, that is, 
the decomposition such that $\langle\,\,,\,\,\rangle\vert_{V_-\times V_-}$ is negative definite, 
$\langle\,\,,\,\,\rangle\vert_{V_+\times V_+}$ is positive definite, 
$\langle\,\,,\,\,\rangle\vert_{V_-\times V_+}=0$ and that $JV_{\pm}=V_{\mp}$ (see Figure 5).  
Define an inner product $\langle\,\,,\,\,\rangle^{\mathcal I}$ of $V$ by 
$$\langle\,\,,\,\,\rangle^{\mathcal I}:=-\pi_{V_-}^{\ast}\langle\,\,,\,\,\rangle
+\pi_{V_+}^{\ast}\langle\,\,,\,\,\rangle,$$
where $\pi_{V_{\pm}}$ denotes the projection of $V$ onto $V_{\pm}$.  
If $(V,\langle\,\,,\,\,\rangle^{\mathcal I})$ is a Hilbert space and that 
the distance topology associated with $\langle\,\,,\,\,\rangle^{\mathcal I}$ 
coincides with the original topology of $V$, then we (\cite{Koi3}) called 
$(V,\langle\,\,,\,\,\rangle,J)$ the {\it infinite dimensional anti-Kaehler space}.  
Here we state that, from each (infinite dimensional separable) Hilbert space, an infinite dimensional anti-Kaehler 
space is constructed in natural manner.  Let $(W,\langle\,\,,\,\,\rangle_W)$ be a (infinite dimensional separable) 
Hilbert space and $V:=W^{\mathbb C}$ be the complexification of $W$.  Also, let 
$\langle\,\,,\,\,\rangle_W^{\mathbb C}\,(:V\times V\to{\mathbb C})$ be the complexification of 
$\langle\,\,,\,\,\rangle_W$.  We regard $V$ as a topological real vector space.  
Define a continuous (${\mathbb R}$-)linear opeartor $J$ of $V$ by 
$Jv:=\sqrt{-1}v\,\,(v\in V)$ and a continuous non-degenerate symmetric (${\mathbb R}$-)bilinear form 
$\langle\,\,,\,\,\rangle$ of $V$ by 
$\langle v_1,v_2\rangle:=2{\rm Re}(\langle v_1,v_2\rangle)\,\,\,(v_1,v_2\in V)$, 
where ${\rm Re}(\cdot)$ is the real part of $(\cdot)$.  
Then $(V,\langle\,\,,\,\,\rangle,J)$ is an infinite dimensional anti-Kaehler space and 
$V=\sqrt{-1}W\oplus W$ is the orthogonal time-space decomposition compatible with $J$.  


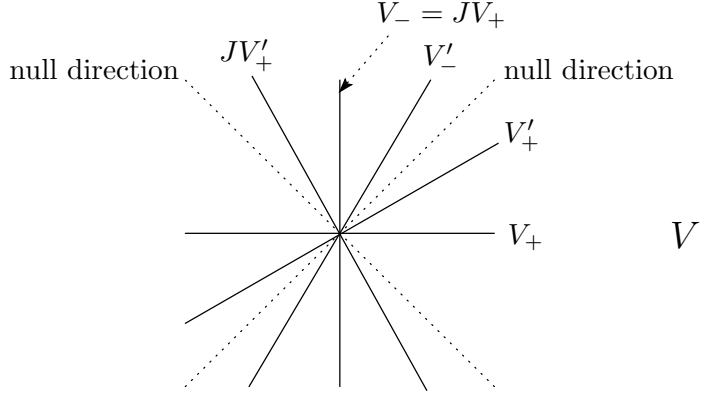
\begin{figure}[h]
\centerline{
\unitlength 0.1in
\begin{picture}( 38.0000, 26.7000)(  3.1000,-32.3000)
%
\special{pn 8}%
\special{pa 2400 2600}%
\special{pa 2400 1000}%
\special{fp}%
%
\special{pn 8}%
\special{pa 1600 1800}%
\special{pa 3200 1800}%
\special{fp}%
%
\special{pn 8}%
\special{pa 1600 2600}%
\special{pa 3200 1000}%
\special{dt 0.045}%
%
\special{pn 8}%
\special{pa 3200 2600}%
\special{pa 1600 1000}%
\special{dt 0.045}%
%
\special{pn 8}%
\special{pa 1930 2600}%
\special{pa 2870 1000}%
\special{fp}%
%
\special{pn 8}%
\special{pa 1600 2270}%
\special{pa 3220 1330}%
\special{fp}%
\put(32.7000,-17.5000){\makebox(0,0)[lt]{$V_+$}}%
\put(25.9000,-7.3000){\makebox(0,0)[lb]{$V_-=JV_+$}}%
%
\special{pn 8}%
\special{pa 2850 2620}%
\special{pa 1946 980}%
\special{fp}%
\put(28.3000,-9.7000){\makebox(0,0)[lb]{$V'_-$}}%
\put(32.5000,-13.7000){\makebox(0,0)[lb]{$V'_+$}}%
\put(20.6000,-9.5000){\makebox(0,0)[rb]{$JV'_+$}}%
\put(7.8000,-30.6000){\makebox(0,0)[lt]{$\langle\,\,,\,\,\rangle\vert_{V'_-\times V'_+}=0,\,\,\langle\,\,,\,\,\rangle\vert_{V'_+\times JV'_+}\not=0,\,\,\,\langle\,\,,\,\,\rangle^{\mathcal I}\vert_{V'_+\times JV'_+}=0$}}%
\put(41.1000,-17.4000){\makebox(0,0)[lt]{{\Large $V$}}}%
\put(14.5000,-32.3000){\makebox(0,0)[lt]{$\,$}}%
\put(14.2000,-28.5000){\makebox(0,0)[lt]{$\langle\,\,,\,\,\rangle\vert_{V_-\times V_+}=0,\,\,\langle\,\,,\,\,\rangle^{\mathcal I}\vert_{V_-\times V_+}=0$}}%
%
\special{pn 8}%
\special{pa 2650 770}%
\special{pa 2400 1060}%
\special{dt 0.045}%
\special{sh 1}%
\special{pa 2400 1060}%
\special{pa 2460 1024}%
\special{pa 2436 1020}%
\special{pa 2428 996}%
\special{pa 2400 1060}%
\special{fp}%
\put(32.5000,-9.9000){\makebox(0,0)[lb]{null direction}}%
\put(15.7000,-9.9000){\makebox(0,0)[rb]{null direction}}%
\end{picture}%
\hspace{2.5truecm}}
\caption{The uniqueness of the orthogonal time-space decomposition compatible with $J$}
\label{figure5}
\end{figure}


Next we recall the notion of an anti-Kaeher Hilbert manifold, where 
we shall define this notion more smartly than the definition in \cite{Koi3}).  
Also, we define the (metrically) completeness of an anti-Kaehler Hilbert manifold.  
Let $N$ be a Hilbert manifold modelled on a (separable) Hilbert space $(V,\langle\,\,,\,\,\rangle_V)$.  
Let $\langle\,\,,\,\,\rangle$ be a (smooth) section of the $(0,2)$-tensor bundle 
$T^{\ast}M\otimes T^{\ast}M$ such that $\langle\,\,,\,\,\rangle_x$ is 
a continuous non-degenerate symmetric bilinear form on $T_xM$ for each $x\in M$.  Also, let $J$ be 
a (smooth) section of the $(1,1)$-tensor bundle $T^{\ast}M\otimes TM$ such that 
$J_x$ is a continuous linear operator of $T_xM$ for each $x\in M$, $J^2=-{\rm id},\,\,\nabla J=0$ and that 
$\langle JX,JY\rangle=-\langle X,Y\rangle$ for every $X,Y\in TM$, where $\nabla$ denotes the Levi-Civita connection 
of $\langle\,\,,\,\,\rangle$.  
For each $x\in M$, there uniquely exists an orthogonal time-space decomposition 
$T_xM=W^x_-\oplus W^x_+$ (with respect to $\langle\,\,,\,\,\rangle_x$) compatible with $J_x$, that is, 
the decomposition such that $\langle\,\,,\,\,\rangle_x\vert_{W^x_-\times W^x_-}$ is negative definite, 
$\langle\,\,,\,\,\rangle_x\vert_{W^x_+\times W^x_+}$ is positive definite, 
$\langle\,\,,\,\,\rangle_x\vert_{W^x_-\times W^x_+}=0$ and that $J_xW^x_{\pm}=W^x_{\mp}$.  
Define an inner product $\langle\,\,,\,\,\rangle^{\mathcal I}_x$ of $T_xM$ by 
$$\langle\,\,,\,\,\rangle^{\mathcal I}_x:=-\pi_{W^x_-}^{\ast}\langle\,\,,\,\,\rangle_x
+\pi_{W^x_+}^{\ast}\langle\,\,,\,\,\rangle_x,$$
where $\pi_{W^x_{\pm}}$ denotes the projection of $T_xM$ onto $W^x_{\pm}$.  
Let $\langle\,\,,\,\,\rangle^{\mathcal I}$ be the section of $T^{\ast}M\otimes TM$ defined by 
assigning $\langle\,\,,\,\,\rangle^{\mathcal I}_x$ to each $x\in M$.  
If $(T_xM,\langle\,\,,\,\,\rangle^{\mathcal I}_x)$ is isometric to $(V,\langle\,\,,\,\,\rangle_V)$ for each $x\in M$, 
then we (\cite{Koi3}) called $(M,\langle\,\,,\,\,\rangle,J)$ an {\it anti-Kaehler Hilbert manifold}.  
If the Riemannian Hilbert manifold $(M,\langle\,\,,\,\,\rangle^{\mathcal I})$ is complete, then we say that 
the anti-Kaehler Hilbert manifold $(M,\langle\,\,,\,\,\rangle)$ is {\it (metrically) complete}.  
Note that the (metrically) completeness of a finite dimensional anti-Kaehler manifold also is defined similarly.  

Let $f$ be an isometric immersion of an anti-Kaehler Hilbert manifold 
$(M,\langle\,\,,\,\,\rangle,J)$ into an anti-Kaehler space $(V,\langle\,\,,\,\,\rangle,\widetilde J)$.  
If $\widetilde J\circ f_{\ast}=f_{\ast}\circ J$ holds, then we (\cite{Koi3}) called $M$ 
an {\it anti-Kaehler Hilbert submanifold in} $(V,\langle\,\,,\,\,\rangle,\widetilde J)$ {\it immersed by} $f$.  
If $M$ is of finite codimension and, for each $v\in T^{\perp}M$, the shape operator $A_v$ is a compact operator 
with respect to $f^{\ast}\langle\,\,,\,\,\rangle^{\mathcal I}$, then we (\cite{Koi3}) called $M$ 
an {\it anti-Kaehler Fredholm submanifold}.  
Assume that $M$ is an embedded anti-Kaehler Fredholm submanifold in $V$, where $M\subset V$ and $f$ is 
the inclusion map of $M$ into $V$.  For the simplicity, denote by the same symbol $J$ the complex structures of 
$M$ and $V$.  Also, denote by $A$ the shape tensor of $M$.  
Fix a unit normal vector $v$ of $M$.  If there exists $X(\not=0)\in TM$ with $A_vX=aX+bJX$, then we call 
the complex number $a+b\sqrt{-1}$ a $J$-{\it eigenvalue of} $A_v$ 
(or a $J$-{\it principal curvature of direction} $v$) and call $X$ a $J$-{\it eigenvector for} $a+b\sqrt{-1}$.  
Also, we call the space of all $J$-eigenvectors for $a+b\sqrt{-1}$ a $J$-{\it eigenspace for} $a+b\sqrt{-1}$.  
The $J$-eigenspaces are orthogonal to one another and they are $J$-invariant, respectively.  
We call the set of all $J$-eigenvalues of $A_v$ the $J$-{\it spectrum of} $A_v$ and denote it by 
${{\rm Spec}}_JA_v$.  Since $M$ is an anti-Kaehler Fredholm submanifold, ${{\rm Spec}}_JA_v\setminus\{0\}$ is 
described as follows:
$${{\rm Spec}}_JA_v\setminus\{0\}
=\{\mu_i\,|\,i=1,2,\cdots\}$$
$$\left(
\begin{array}{c}
\displaystyle{|\mu_i|>|\mu_{i+1}|\,\,\,{{\rm or}}
\,\,\,{\rm "}|\mu_i|=|\mu_{i+1}|\,\,\&\,\,
{{\rm Re}}\,\mu_i>{{\rm Re}}\,\mu_{i+1}{\rm "}}\\
\displaystyle{{{\rm or}}\,\,\,"|\mu_i|=|\mu_{i+1}| 
\,\,\&\,\,
{{\rm Re}}\,\mu_i={{\rm Re}}\,\mu_{i+1}\,\,\&\,\,
{{\rm Im}}\,\mu_i=-{{\rm Im}}\,\mu_{i+1}>0"}
\end{array}
\right).$$
Also, the $J$-eigenspace for each $J$-eigenvalue of $A_v$ other than $0$ is of finite dimension.  
We call the $J$-eigenvalue $\mu_i$ the $i$-{\it th} $J$-{\it principal curvature of direction} $v$.  
Assume that the normal holonomy group of $M$ is trivial.  
Fix a parallel normal vector field $\widetilde v$ of $M$.  
Assume that the number (which may be $\infty$) of distinct $J$-principal 
curvatures of direction $\widetilde v_x$ is independent of the choice of 
$x\in M$.  
Then we can define complex-valued functions $\widetilde{\mu}_i$ 
($i=1,2,\cdots$) on $M$ 
by assigning the $i$-th $J$-principal curvature of direction 
$\widetilde v_x$ to each $x\in M$.  We call this function 
$\widetilde{\mu}_i$ the $i$-{\it th} $J$-{\it principal curvature function of 
direction} $\widetilde v$.  
The submanifold $M$ is called an {\it anti-Kaehler isoparametric 
submanifold} if it satisfies the following condition:

\vspace{0.2truecm}

For each parallel normal vector field $\widetilde v$ of $M$, the 
number of distinct $J$-principal curvatures 

of direction $\widetilde v_x$ is independent of the choice of 
$x\in M$, each $J$-principal curvature function 

of direction $\widetilde v$ is constant 
on $M$ and it has constant multiplicity.  

\vspace{0.2truecm}

\noindent
Let $\{e_i\}_{i=1}^{\infty}$ be an orthonormal system of $T_xM$.  If 
$\{e_i\}_{i=1}^{\infty}\cup\{Je_i\}_{i=1}^{\infty}$ is an orthonormal base 
of $T_xM$, then we call $\{e_i\}_{i=1}^{\infty}$ (rather than 
$\{e_i\}_{i=1}^{\infty}\cup\{Je_i\}_{i=1}^{\infty}$) a 
$J$-{\it orthonormal base}.  
If there exists a $J$-orthonormal base consisting of $J$-eigenvectors of 
$A_v$, then $A_v$ is said to {\it be diagonalized with respect to the} 
$J$-{\it orthonormal base}.  If $M$ is anti-Kaehler isoparametric and, 
for each $v\in T^{\perp}M$, the shape operator $A_v$ is 
diagonalized with respect to a $J$-orthonormal base, then we (\cite{Koi3}) called 
$M$ a {\it proper anti-Kaehler isoparametric submanifold}.  
We named thus in similar to the terminology 
``{\it proper} isoparametric semi-Riemannian submanifold'' used in \cite{Koi1}.  
Throughout this paper, we shall call this submanifold 
an {\it anti-Kaehler isoparametric submanifold with $J$-diagonalizable shape operators} because this terminology 
seems to be more familiar than ``proper anti-Kaehler isoparametric submanifold''.  
Assume that $M$ is an anti-Kaehler isoparametric submanifold  with $J$-diagonalizable shape operators.  
Then, since the ambient space is flat and the normal holonomy group of $M$ is trivial, 
it follows from the Ricci equation that the shape operators $A_{v_1}$ 
and $A_{v_2}$ commute for arbitrary two unit normal vector $v_1$ and $v_2$ of $M$.  
Hence the shape operators $A_v$'s ($v\in T^{\perp}_xM$) 
are simultaneously diagonalized with respect to a $J$-orthonormal base.  
Let $\{E_i\,|\,i\in I\}$ be the family of distributions on $M$ such that, 
for each $x\in M$, 
$\{(E_i)_x\,|\,i\in I\}$ is the set of all common $J$-eigenspaces of 
$A_v$'s ($v\in T^{\perp}_xM$).  For each $x\in M$, we have 
$T_xM=\overline{\displaystyle{\mathop{\oplus}_{i\in I}(E_i)_x}}$, where 
$\overline{\displaystyle{\mathop{\oplus}_{i\in I}(E_i)_x}}$ denotes the closure of 
$\displaystyle{\mathop{\oplus}_{i\in I}(E_i)_x}$ with respect to $\langle\,\,,\,\,\rangle^{\mathcal I}_x$.  
We regard $T^{\perp}_xM$ ($x\in M$) as a complex vector space by 
$J_x|_{T^{\perp}_xM}$ and denote the dual space of the complex vector 
space $T^{\perp}_xM$ by $(T^{\perp}_xM)^{\ast_{\mathbb C}}$.  
Also, denote by $(T^{\perp}M)^{\ast_{\mathbb C}}$ the complex vector bundle over 
$M$ having $(T^{\perp}_xM)^{\ast_{\mathbb C}}$ as the fibre over $x$.  
Let $\lambda_i$ ($i\in I$) be the section of $(T^{\perp}M)^{\ast_{\mathbb C}}$ 
such that $A_v={{\rm Re}}(\lambda_i)_x(v){{\rm id}}+{{\rm Im}}(\lambda_i)_x(v)
J_x$ on $(E_i)_x$ for any $x\in M$ and any $v\in T^{\perp}_xM$.  
We call $\lambda_i$ ($i\in I$) $J$-{\it principal curvatures} of $M$ and 
$E_i$ ($i\in I$) $J$-{\it curvature distributions} of $M$.  
The distribution $E_i$ is integrable and each leaf of $E_i$ is a complex 
sphere.  Each leaf of $E_i$ is called a {\it complex curvature sphere}.  
It is shown that there uniquely exists a normal vector field $n_i$ of $M$ with 
$\lambda_i(\cdot)
=\langle n_i,\cdot\rangle-\sqrt{-1}\langle Jn_i,\cdot\rangle$ 
(see Lemma 5 of \cite{Koi3}).  We call $n_i$ ($i\in I$) the 
$J$-{\it curvature normals} of $M$.  Note that $n_i$ is parallel with 
respect to the complexification of the normal connection of $M$.  
Note that similarly are defined a (finite dimensional) proper anti-Kaehler 
isoparametric submanifold in a finite dimensional anti-Kaehler space, its 
$J$-principal curvatures, its $J$-curvature distributions and its 
$J$-curvature normals.  
Set ${\it l}^x_i:=(\lambda_i)_x^{-1}(1)$.  According to (i) of Theorem 2 in 
\cite{Koi3}, the tangential focal set of $M$ at $x$ is equal to 
$\displaystyle{\mathop{\cup}_{i\in I}{\it l}_i^x}$.  
We call each ${\it l}_i^x$ a {\it complex focal hyperplane of} $M$ {\it at} 
$x$.  Let $\widetilde v$ be a parallel normal vector field of $M$.  
If $\widetilde v_x$ belongs to at least one ${\it l}_i$, then it is called 
a {\it focal normal vector field} of $M$.  
For a focal normal vetor field $\widetilde v$, the focal map 
$f_{\widetilde v}$ is defined by $f_{\widetilde v}(x):=\exp^{\perp}
(\widetilde v_x)\,\,\,(x\in M)$.  The image $f_{\widetilde v}(M)$ is called a 
{\it focal submanifold} of $M$, where we denote by $F_{\widetilde v}$.  
For each $x\in F_{\widetilde v}$, the inverse image $f_{\widetilde v}^{-1}(x)$ 
is called a {focal leaf} of $M$.  
Denote by $T_i^x$ the complex reflection of order $2$ with respect to 
${\it l}_i^x$ (i.e., the rotation of angle $\pi$ having ${\it l}_i^x$ as the 
axis), which is an affine transformation of $T^{\perp}_xM$.  
Let ${\mathcal W}_x$ be the group generated by $T_i^x$'s ($i\in I$).  
According to Proposition 3.7 of \cite{Koi5}, ${\mathcal W}_x$ is discrete.  
Furthermore, it follows from this fact that ${\mathcal W}_x$ is isomorphic 
an affine Weyl group.  This group ${\mathcal W}_x$ is independent of the choice 
of $x\in M$ (up to group isomorphicness).  
Hence we simply denote it by ${\mathcal W}$.  We call this group ${\mathcal W}$ the 
{\it complex Coxeter group associated with} $M$.  
According to Lemma 3.8 of \cite{Koi5}, 
$W$ is decomposable (i.e., it is decomposed into a non-trivial product of 
two discrete complex reflection groups) if and only if there exist two 
$J$-invariant linear subspaces $P_1$ ($\not=\{0\}$) and $P_2$ ($\not=\{0\}$) 
of $T^{\perp}_xM$ such that $T^{\perp}_xM=P_1\oplus P_2$ (orthogonal 
direct sum), $P_1\cup P_2$ contains all $J$-curvature normals of 
$M$ at $x$ and that $P_i$ ($i=1,2$) contains at least one $J$-
curvature normal of $M$ at $x$.  
Also, according to Theorem 1 of \cite{Koi5}, $M$ is irreducible if and only if 
${\mathcal W}$ is not decomposable.  

\section{Parallel transport map} 
Y. Maeda, S. Rosenberg and P. Tondeur (\cite{MRT}) studied the minimality of the Gauge orbit in the space of 
the $H^0$-connections of the principal bundle $P$ having a compact semi-simple Lie group $G$ as the structure group 
over a compact Riemannian manifold $M$.  Let $c^{\ast}P$ be the pull-back bundle of $P$ by 
a $C^{\infty}$-path $c:[0,1]\to M$.  The space of $H^0$-connections on $c^{\ast}P$ is identified with the (seprable) 
Hilbert space $H^0([0,1],\mathfrak g)$ of the $H^0$-paths in the Lie algebra $\mathfrak g$ of $G$.  
The Hilbert Lie group $\Omega_e(G)(\subset H^1([0,1],G))$ of $H^1$-loops at $e$ in $G$ acts on 
$H^0([0,1],\mathfrak g)$ as the subaction of the Gauge group on the space of connections, where $e$ is the identity 
element of $G$.  The orbit map $\phi:H^0([0,1],\mathfrak g)\to H^0([0,1],\mathfrak g)/\Omega_e(G)\,(=G)$ is called 
the {\it parallel transport map} for $G$.  
See \cite{PiTh}, \cite{Te2}, \cite{Te3}, \cite{TT} and \cite{PaTe} about 
the study of the parallel transport map for a compact semi-simple Lie group.  
Now we shall consider the case where $[0,1]$ is replaced by the circle $S^1$ in the above 
definition.  Then we shall explain that $\phi$ should be called the {holonomy map} for $G$.  
Let $\gamma:S^1\to M$ be a $C^{\infty}$-loop.  
The space of $H^0$-connections on $\gamma^{\ast}P$ is identified with 
the (seprable) Hilbert space $H^0(S^1,\mathfrak g)$ of the $H^0$-loop in the Lie algebra $\mathfrak g$ of $G$.  
The Hilbert Lie group $\Omega_e(G)(\subset H^1(S^1,G))$ of $H^1$-loops at $e$ in $G$ acts on 
$H^0(S^1,\mathfrak g)$ as the subaction of the Gauge group on the space of connections.  
We consider the orbit map $\phi:H^0(S^1,\mathfrak g)\to H^0(S^1,\mathfrak g)/\Omega_e(G)\,(=G)$.  
Then, for each $C^{\infty}$-loop $\gamma:S^1\to M$, $\phi(\gamma)$ is equal to the generator of the holonomy group 
(which is a cyclic subgroup of $G$) of the connection $\omega_{\gamma}$ of the trivial $G$-bundle 
$S^1\times G\to S^1$ determined by $\gamma$.  In this sense, $\phi$ should be called the {holonomy map} for $G$.  

We (\cite{Koi3}) defined the notion of the parallel transport map for the 
complexification $G^{\mathbb C}$ of a semi-simple Lie group $G$.  In this section, we recall this notion and 
some facts related to this notion, which will be used in Sections 7 and 8.  
Let $K$ be a maximal compact subgroup of $G$, $\mathfrak g$ 
(resp. $\mathfrak k$) the Lie algebra of $G$ (resp. $K$) and 
$\mathfrak g=\mathfrak k\oplus\mathfrak p$ a Cartan decomposition of 
$\mathfrak g$.  Also, let $\langle\,\,,\,\,\rangle$ be the 
${\rm Ad}_G(G)$-invariant non-degenerate inner product of 
$\mathfrak g$.  The Cartan decomposition 
$\mathfrak g=\mathfrak k\oplus\mathfrak p$ is an orthogonal time-space 
decomposition of $\mathfrak g$ with respect to $\langle\,\,,\,\,\rangle$, 
that is, $\langle\,\,,\,\,\rangle|_{\mathfrak k\times\mathfrak k}$ is 
negative definite, 
$\langle\,\,,\,\,\rangle|_{\mathfrak p\times\mathfrak p}$ is positive 
definite and $\langle\,\,,\,\,\rangle|_{\mathfrak k\times\mathfrak p}$ 
vanishes.  
Set $\langle\,\,,\,\,\rangle^A:=2{\rm Re}\langle\,\,,\,\,\rangle^{\mathbb C}$, 
where $\langle\,\,,\,\,\rangle^{\mathbb C}$ is the complexification of 
$\langle\,\,,\,\,\rangle$ (which is a ${\mathbb C}$-bilinear form of 
$\mathfrak g^{\mathbb C}$).  
The ${\mathbb R}$-bilinear form $\langle\,\,,\,\,\rangle^A$ on 
$\mathfrak g^{\mathbb C}$ regarded as a real Lie algebra induces 
a bi-invariant pseudo-Riemannian metric on $G^{\mathbb C}$ and furthermore 
a $G^{\mathbb C}$-invariant anti-Kaehler metric on $G^{\mathbb C}/K^{\mathbb C}$.  
It is clear that $\mathfrak g^{\mathbb C}=(\mathfrak k\oplus\sqrt{-1}\mathfrak p)
\oplus(\sqrt{-1}\mathfrak k\oplus\mathfrak p)$ is an orthogonal time-space 
decomposition of $\mathfrak g^{\mathbb C}$ with respect to 
$\langle\,\,,\,\,\rangle^A$.  
For the simplicity, set $\mathfrak g^{\mathbb C}_-:=\mathfrak k\oplus\sqrt{-1}\mathfrak p$ 
and $\mathfrak g^{\mathbb C}_+:=\sqrt{-1}\mathfrak k\oplus\mathfrak p$.  
Note that $\mathfrak g^{\mathbb C}_-$ is the compact real form of 
$\mathfrak g^{\mathbb C}$.  
Set $\langle\,\,,\,\,\rangle^{\mathcal I}:=-
\pi^{\ast}_{\mathfrak g^{\mathbb C}_-}\langle\,\,,\,\,\rangle^A+
\pi^{\ast}_{\mathfrak g^{\mathbb C}_+}\langle\,\,,\,\,\rangle^A$, where 
$\pi_{\mathfrak g^{\mathbb C}_-}$ (resp. $\pi_{\mathfrak g^{\mathbb C}_+}$) is the 
projection of $\mathfrak g^{\mathbb C}$ onto $\mathfrak g^{\mathbb C}_-$ (resp. 
$\mathfrak g^{\mathbb C}_+$).  
Let $H^0([0,1],\mathfrak g^{\mathbb C})$ be the space of all $L^2$-integrable 
paths $u:[0,1]\to\mathfrak g^{\mathbb C}$ with respect to 
$\langle\,\,,\,\,\rangle^{\mathcal I}$ and 
$H^0([0,1],\mathfrak g^{\mathbb C}_-)$ (resp. 
$H^0([0,1],\mathfrak g^{\mathbb C}_+)$) the space of all 
$L^2$-integrable paths $u:[0,1]\to\mathfrak g^{\mathbb C}_-$ (resp. 
$u:[0,1]\to\mathfrak g^{\mathbb C}_+$) with respect to 
$-\langle\,\,,\,\,\rangle^A|_{\mathfrak g^{\mathbb C}_-\times
\mathfrak g^{\mathbb C}_-}$ (resp. 
$\langle\,\,,\,\,\rangle^A|_{\mathfrak g^{\mathbb C}_+
\times\mathfrak g^{\mathbb C}_+}$).  Clearly we have 
$H^0([0,1],\mathfrak g^{\mathbb C})=H^0([0,1],\mathfrak g^{\mathbb C}_-)
\oplus H^0([0,1],\mathfrak g^{\mathbb C}_+)$.  
Define a non-degenerate inner product $\langle\,\,,\,\,\rangle^A_0$ 
of $H^0([0,1],\mathfrak g^{\mathbb C})$ by 
$\langle u,v\rangle^A_0:=\int_0^1\langle u(t),v(t)\rangle^Adt$.  
It is easy to show that the decomposition 
$H^0([0,1],\mathfrak g^{\mathbb C})=
H^0([0,1],\mathfrak g^{\mathbb C}_-)\oplus H^0([0,1],\mathfrak g^{\mathbb C}_+)$ 
is an orthogonal time-space decomposition with respect to 
$\langle \,\,,\,\,\rangle^A_0$.  For the simplicity, set 
$H^{0,{\mathbb C}}_{\varepsilon}:=H^0([0,1],\mathfrak g^{\mathbb C}_{\varepsilon})$ 
($\varepsilon=-$ or $+$) and 
$\langle \,\,,\,\,\rangle^{\mathcal I}_0:=
-\pi^{\ast}_{H^{0,{\mathbb C}}_-}\langle\,\,,\,\,\rangle^A_0+
\pi^{\ast}_{H^{0,{\mathbb C}}_+}\langle\,\,,\,\,\rangle^A_0$, where 
$\pi_{H^{0,{\mathbb C}}_-}$ (resp. $\pi_{H^{0,{\mathbb C}}_+}$) 
is the projection of $H^0([0,1],\mathfrak g^{\mathbb C})$ onto 
$H^{0,{\mathbb C}}_-$ (resp. $H^{0,{\mathbb C}}_+$).  
It is clear that $\langle u,v\rangle^{\mathcal I}_0
=\int_0^1\langle u(t),v(t)\rangle^{\mathcal I}dt$ 
($u,\,v\in H^0([0,1],\mathfrak g^{{{\mathbb C}}})$).  
Hence $(H^0([0,1],\mathfrak g^{{{\mathbb C}}}),\,\langle\,\,,\,\,
\rangle^{\mathcal I}_0)$ is a Hilbert space, that is, 
$(H^0([0,1],\mathfrak g^{{{\mathbb C}}}),\,\langle\,\,,\,\,\rangle^A_0)$ 
is a pseudo-Hilbert space in the sense of \cite{Koi2}.  
Let $J$ be the endomorphism of $\mathfrak g^{\mathbb C}$ defined by 
$JX=\sqrt{-1}X$ ($X\in \mathfrak g^{{{\mathbb C}}}$).  Denote by the same 
symbol $J$ the bi-invariant almost complex structure of 
$G^{{{\mathbb C}}}$ induced from $J$.  
Define the endomorphism $\widetilde J$ of $H^0([0,1],\mathfrak g^{\mathbb C})$ by 
$\widetilde Ju=\sqrt{-1}u$ 
($u\in H^0([0,1],\mathfrak g^{{{\mathbb C}}})$).  From 
$\widetilde JH^{0,{{\mathbb C}}}_+=H^{0,{{\mathbb C}}}_-$, 
$\widetilde JH^{0,{{\mathbb C}}}_-=H^{0,{{\mathbb C}}}_+$ and 
$\langle\widetilde Ju,\widetilde Jv\rangle^A_0=-\langle u,v\rangle^A_0$ 
($u,v\in H^0([0,1],\mathfrak g^{{{\mathbb C}}})$), the space $(H^0([0,1],
\mathfrak g^{{{\mathbb C}}}),\langle\,\,,\,\,\rangle^A_0,\widetilde J)$ is 
an anti-Kaehler space.  Let 
$H^1([0,1],\mathfrak g^{\mathbb C})$ be a pseudo-Hilbert subspace of 
$H^0([0,1],\mathfrak g^{\mathbb C})$ consisting of 
all absolutely continuous paths $u:[0,1]\to \mathfrak g^{\mathbb C}$ 
such that the weak derivative $u'$ of $u$ is squared integrable 
(with respect to $\langle\,\,,\,\,\rangle^{\mathcal I}$).  
Also, let $H^1([0,1],G^{{{\mathbb C}}})$ be the Hilbert Lie group of 
all absolutely continuous paths $g:[0,1]\to G^{{{\mathbb C}}}$ 
such that the weak derivative $g'$ of $g$ is squared 
integrable (with respect to 
$\langle\,\,,\,\,\rangle^{\mathcal I}$), that is, 
$g_{\ast}^{-1}g'\in H^0([0,1],\mathfrak g^{\mathbb C})$.  
Define a map $\phi:H^0([0,1],\mathfrak g^{\mathbb C})\to G^{\mathbb C}$ by 
$\phi(u):=g_u(1)$ ($u\in H^0([0,1],\mathfrak g^{\mathbb C})$), 
where $g_u$ is the element of $H^1([0,1],G^{\mathbb C})$ with $g_u(0)=e$ and 
$g_{u\ast}^{-1}g'_u=u$.  
This map is called the {\it parallel transport map for} $G^{\mathbb C}$.  
This map is an anti-Kaehler submersion.  
Set $P(G^{{{\mathbb C}}},e\times G^{{{\mathbb C}}}):=
\{g\in H^1([0,1],G^{{{\mathbb C}}})\,|\,g(0)=e\}$ and 
$\Omega_e(G^{{{\mathbb C}}}):=
\{g\in H^1([0,1],G^{{{\mathbb C}}})\,|\,g(0)=g(1)=e\}$.  
The group $H^1([0,1],G^{{{\mathbb C}}})$ acts on $H^0([0,1],\mathfrak g^{{{\mathbb C}}})$ as the action of 
the gauge transformation group on the space of connections, that is, 
$$g\ast u:={\rm Ad}_{G^{\mathbb C}}(g)u-g'g_{\ast}^{-1}\qquad
(g\in H^1([0,1],G^{{{\mathbb C}}}),\,\,u\in H^0([0,1],\mathfrak g^{{{\mathbb C}}})).$$
It is shown that the following facts hold:

\vspace{0.2truecm}

(i) The above action of $H^1([0,1],G^{{{\mathbb C}}})$ on 
$H^0([0,1],\mathfrak g^{{{\mathbb C}}})$ is isometric,

(ii) The above action of $P(G^{{{\mathbb C}}},e\times G^{{{\mathbb C}}})$ 
in $H^0([0,1],\mathfrak g^{{{\mathbb C}}})$ is transitive and free, 

(iii) $\phi(g\ast u)=(L_{g(0)}\circ R_{g(1)}^{-1})
(\phi(u))$ for $g\in H^1([0,1],G^{{{\mathbb C}}})$ and $u\in$ $H^0([0,1],\mathfrak g^{{{\mathbb C}}})$,  

(iv) $\phi:H^0([0,1],\mathfrak g^{{{\mathbb C}}})\to 
G^{{{\mathbb C}}}$ is regarded as a $\Omega_e(G^{{{\mathbb C}}})$-bundle.  

(v) If $\phi(u)=(L_{x_0}\circ R_{x_1}^{-1})(\phi(v))$ 
($u,v\in H^0([0,1],\mathfrak g^{{{\mathbb C}}}),\,\,
x_0,x_1\in G^{{{\mathbb C}}})$, then there exists 

\hspace{0.4truecm} $g\in H^1([0,1],G^{{{\mathbb C}}})$ such that 
$g(0)=x_0,\,\,g(1)=x_1$ and $u=g\ast v$.  In particular, 

\hspace{0.4truecm} it follows that any $u\in H^0([0,1],
\mathfrak g^{{{\mathbb C}}})$ is described as $u=g\ast\hat 0$ in terms of some 

\hspace{0.4truecm} $g\in P(G^{{{\mathbb C}}},G^{{{\mathbb C}}}\times e)$.  


\section{Partial tubes} 
In this section, we recall some facts for partial tubes in a symmetric space, which will use to prove Theorem A in 
the next section.  
For a submanifold $F$ in a symmetric space $G/K$ of non-positive (or non-negative) curvature, 
M. Br$\ddot u$ck (\cite{Br}) introduced a certain kind of partial tube with flat section including the normal 
holonomy tube, where $F$ is assumed to admit the $\varepsilon$-tube for a sufficiently small positive number 
$\varepsilon$.  This notion is defined as follows.  
Let $\varepsilon_{\gamma}:={\rm inf}\{\vert r\vert\,\vert\,r:{\rm focal}\,\,
{\rm radius}\,\,{\rm of}\,\,M\,\,{\rm along}\,\,\gamma\}$, where $\gamma$ 
is a unit speed normal geodesic of $F$.  Set 
$$\varepsilon_F:={\rm inf}\{\varepsilon_{\gamma}\,\vert\,\gamma:{\rm unit}\,\,{\rm speed}\,\,{\rm normal}\,\,
{\rm geodesic}\,\,{\rm of}\,\,F\}.$$
Assume that $\varepsilon_F>0$.  Fix $x_0\in F$.  Let 
$\mathfrak C_{x_0}:=\{c:[0,1]\to F\,:\,{\rm a}\,\,{\rm piecewise}\,\,{\rm smooth}\,\,{\rm path}\,\,{\rm with}\,\,
c(0)=x_0\}$, $\Phi_{x_0}^0$ be the restricted normal holonomy group 
of $F$ at $x_0$ and $\mathfrak L_{x_0}$ be the Lie subalgebra of 
$\mathfrak{so}(T^{\perp}_{x_0}F)$ generated by $\{P_c^{-1}\circ
{\rm pr}_{T^{\perp}_{c(1)}F}\circ R_{c(1)}(P_cv_1,P_cv_2)\circ 
P_c\,\vert\,v_1,v_2\in T^{\perp}_{x_0}M,\,\,c\in\mathfrak C_{x_0}\}$, 
where $R$ denotes the curvature tensor of $G/K$ and $P_c$ is the parallel transport along $c$ with respect to 
the normal connection $\nabla^{\perp}$ of $F$ and ${\rm pr}_{T^{\perp}_{c(1)}F}$ is the orthogonal projection onto 
$T^{\perp}_{c(1)}F$.  Also, let $\widehat{\mathfrak L}_{x_0}$ be the Lie algebra generated by 
$\mathfrak L_{x_0}$ and ${\rm Lie}\,\Phi_{x_0}^0$.  Let $L_{x_0}:=\exp\,
\mathfrak L_{x_0}$ and $\widehat L_{x_0}:=\exp\,\widehat{\mathfrak L}_{x_0}$, 
where $\exp$ is the exponential map of $GL(T^{\perp}_{x_0}F)$.  Note that $L_{x_0}$ and $\widehat L_{x_0}$ are 
Lie subgroups of $SO(T^{\perp}_{x_0}F)$.  
For $v_0\in T^{\perp}_{x_0}F$, define a subbundle $B_{v_0}(F)$ of 
the normal bundle $T^{\perp}F$ of $F$ by 
$$B_{v_0}(F):=\{P_c(gv_0)\,\vert\,g\in\widehat L_{x_0},\,\,\,
c\in\mathfrak C_{x_0}\}$$
and $\widetilde B_{v_0}(F):=\exp^{\perp}(B_{v_0}(F))$, where $\exp^{\perp}$ denotes the normal exponential map of 
$F$.  For each vector $v_0$ with $\vert\vert v_0\vert\vert\,<\,\varepsilon_F$, 
$\widetilde B_{v_0}(F)$ is an immersed submanifold, that is, a partial tube over $F$ whose fibre over $x_0$ is 
$\exp^{\perp}(\widehat L_{x_0}v_0)$.  M. Br$\ddot{\rm u}$ck proved the following fact.  

\vspace{0.5truecm}

\noindent
{\bf Theorem 5.1(\cite{Br}).} 
{\sl Let $M$ be an equifocal submanifold in a symmetric space of non-positive 
(or non-negative) curvature having a focal submanifold $F$.  If the sections of $M$ are properly embedded, 
then $M$ is equal to the partial tube $\widetilde B_{v_0}(F)$, where $v_0$ is an element of $T^{\perp}F$ with 
$\exp^{\perp}(v_0)\in M$, and each fibre of $\widetilde B_{v_0}(F)(=M)$ is the image by $\exp^{\perp}$ of 
a principal orbit of an orthogonal representation on the normal space of $F$ which is 
equivalent to the direct sum representation of some s-represenations and a trivial representation.}

\vspace{0.5truecm}

According to the proof of this theorem in \cite{Br}, we can derive the following fact.  

\vspace{0.5truecm}

\noindent
{\bf Theorem 5.2.} 
{\sl Let $M$ be an isoparametric submanifold in a symmetric space of non-compact type 
having a focal submanifold $F$.  If the sections of $M$ are properly embedded, 
then $M$ is equal to the partial tube $\widetilde B_{v_0}(F)$, where $v_0$ is an element of $T^{\perp}F$ with 
$\exp^{\perp}(v_0)\in M$, and each fibre of $\widetilde B_{v_0}(F)(=M)$ is the image by $\exp^{\perp}$ of 
a principal orbit of an orthogonal representation on the normal space of $F$ which is 
equivalent to the direct sum representation of some s-represenations and a trivial representation.}

\vspace{0.5truecm}

We recall the notion of a (general) partial tube.  
Let $F$ be a submanifold in a Riemannian manifold $N$, $T^{\perp}F$ be the normal bundle of $F$, 
$\exp^{\perp}$ be the normal exponential map of $F$, $\nabla^{\perp}$ be the normal connection of $F$ and 
$\bar A$ be the shape tensor of $F$.  
Let $t(F)$ be a submanifold of $T^{\perp}F$ which is given as the sum of some normal holonomy subbundles of 
$T^{\perp}F$ and $\widetilde t(F):=\exp^{\perp}(t(F))$.  If $\exp^{\perp}|_{t(F)}$ is an immersion, then 
$\widetilde t(F)$ is called a {\it partial tube} over $F$.  
Denote by $A$ the shape tensor of $\widetilde t(F)$ and ${\mathcal V}$ (resp. ${\mathcal H}$) 
the vertical distribution (resp. the horizontal distribution) on $T^{\perp}F$, where ``horizontal distribution'' 
means that it is horizontal with respect to $\nabla^{\perp}$.  
Denote by $\widetilde X_{\xi}$ the horizontal lift of $X\in T_xF$ to $\xi\in\widetilde t(F)_x$.  
Denote by $A^x$ the shape tensor of the fibre $\widetilde t(F)_x:=\exp^{\perp}(t(F)\cap T^{\perp}_xF)$ of 
$\widetilde t(F)$ over $x(\in F)$ in the normal umbrella $\Sigma_x:=\exp^{\perp}(T^{\perp}_xF)$.  
For a $C^{\omega}$-function $\Psi$ and a linear transformation $Q$, we define a linear transformation 
$\Psi(Q)$ by 
$$\Psi(Q):=\sum_{k=0}^{\infty}\frac{\Psi^{(k)}(0)}{k!}Q^k.$$

According to the proof of Proposition 3.1 and Corollary 3.2 in \cite{Koi4}, we have the following facts for 
$A$.  

\vspace{0.5truecm}

\noindent
{\bf Proposition 5.3.} {\sl Assume that $N$ is a symmetric space $G/K$ of compact type or non-compact type 
and that $F$ is a submanifold with section, where ``with section'' means that the normal umbrella 
$\Sigma_y:=\exp^{\perp}(T^{\perp}_yF)$ at each $y\in F$ is totally geodesic in $G/K$ ($\Sigma_y$ is then called 
the section of $F$ through $y$).  
Let $v\in t(F)_x:=t(F)\cap T^{\perp}_xF$ and $w\in T^{\perp}_v\widetilde t(F)$.  
Then the following statements {\rm (i)} and {\rm(ii)} hold:

{\rm (i)} For $X\in{\mathcal V}_v$, we have $A_wX=A^x_wX$;

{\rm (ii)} Set $\bar w:=(P_{\gamma_v|_{[0,1]}})^{-1}(w)$, where $\gamma_v$ is the geodesic in $G/K$ with 
$\gamma'_v(0)=v$ and $P_{\gamma_v|_{[0,1]}}$ is the parallel transport map along $\gamma_v|_{[0,1]}$.  
Assume that the sectional curvature for the $2$-plane ${\rm Span}\{v,\bar w\}$ is equal to zero.  
For $Y\in{\mathcal H}_v$, we have 
$$\begin{array}{l}
\displaystyle{A_w\widetilde Y_v=P_{\gamma_v|_{[0,1]}}
\left(-({\rm ad}(\bar w)\circ\sinh({\rm ad}(v)))(Y)+\frac{\sinh({\rm ad}(v))}{{\rm ad}(v)}
(\bar A_{\bar w}Y)\right.}\\
\hspace{2.95truecm}\displaystyle{
\left.+\left(\left(\frac{\cosh({\rm ad}(v))-{\rm id}}{{\rm ad}(v)}-\frac{\sinh({\rm ad}(v))
-{\rm ad}(v)}{{\rm ad}(v)^2}\right)\circ{\rm ad}(\bar w)\right)(\bar A_vY)\right),}
\end{array}$$
where ${\rm ad}$ denotes the adjoint representation of the Lie algebra $\mathfrak g$ of $G$, and ${\rm ad}(v)^2$ and 
${\rm ad}(v)\circ{\rm ad}(\bar w)$ are regarded as a linear transformations of 
$\mathfrak p:={\rm Ker}(\theta_{\ast e}+{\rm id})(\approx T_{eK}(G/K))$ ($\theta\,:\,$ the Cartan involution of $G$ 
with $({\rm Fix}\,\theta)_0\subset K\subset{\rm Fix}\,\theta$).  
In particular, if $F$ is reflective, $R(v)Y=b_1^2Y$ and $R(\bar w)Y=b_2^2Y$, then we have 
$$A_w\widetilde Y_v=-b_2\tanh b_1\widetilde Y_v,$$
where $R$ denotes the curvature tensor of $G/K$, $R(\bullet)$ denotes the normal Jacobi operator for 
$(\bullet)$ and $b_i$ ($i=1,2$) are real numbers (resp. purely imaginary numbers) when 
$G/K$ is of non-compact type (resp. of compact type).}

\section{Proof of Theorem A} 
In this section, we shall prove Theorem A.  
Let $M$ be as in Theorem A and $F$ be a reflective focal submanifold of $M$.  
Denote by $A$ the shape tensor of $M$ and $R$ the curvature tensor of $G/K$.  
Without loss of generality, we may assume that $o:=eK\in F$.  
Let $\mathfrak g$ (resp. $\mathfrak k$) be the Lie algebra of $G$ (resp. $K$) and $\theta$ be an Cartan involution of 
$G$ with $({\rm Fix}\,\theta)_0\subset K\subset{\rm Fix}\,\theta$, where ${\rm Fix}\,\theta$ denotes the fixed 
point set of $\theta$ and $({\rm Fix}\,\theta)_0$ denotes the identity component of ${\rm Fix}\,\theta$.  
Denote by the same symbol $\theta$ the involution (i.e., $\theta_{\ast e}$) of $\mathfrak g$ induced form $\theta$ 
and set $\mathfrak p:={\rm Ker}(\theta+{\rm id}_{\mathfrak g})$, which is identified with the tangent space 
$T_o(G/K)$.  Denote by ${\rm Exp}$ the exponential map of $G/K$ at $o$.  

\vspace{0.5truecm}

\noindent
{\it Proof of Theorem A.} 
Take $Z_0\in\mathfrak p$ with ${\rm Exp}\,Z_0\in M$.  Set $x_0:={\rm Exp}\,Z_0,\,\,
\mathfrak t:=T_oF,\,\mathfrak t^{\perp}:=T^{\perp}_oF$ and 
$\mathfrak b:=(\exp\,Z_0)_{\ast o}^{-1}(T^{\perp}_{x_0}M)$.  Since $F$ is reflective, 
both $\mathfrak t$ and $\mathfrak t^{\perp}$ are Lie triple systems.  Also it is clear that $\mathfrak b$ is a 
maximal abelian subspace of $\mathfrak t^{\perp}$ containinig $Z_0$.  Take a maximal abelian subspace $\mathfrak a$ 
of $\mathfrak p$ including $\mathfrak b$.  Let $\triangle$ be a (restricted) root system with respect to 
$\mathfrak a$ and set $\triangle_{\mathfrak b}:=\{\alpha|_{\mathfrak b}\,\,|\,\,\alpha\in\triangle\}$.  
Let $(\triangle_{\mathfrak b})_+$ be the positive root system under a lexicographic ordering of 
$\mathfrak b^{\ast}$, $\mathfrak p_{\beta}$ be the root space for $\beta\in(\triangle_{\mathfrak b})_+$.  
Set $(\triangle_{\mathfrak b})^V_+:=\{\beta\in(\triangle_{\mathfrak b})_+\,|\,\mathfrak p_{\beta}\cap
\mathfrak t^{\perp}\not=\{0\}\}$ and $(\triangle_{\mathfrak b})^H_+:=\{\beta\in(\triangle_{\mathfrak b})_+\,|\,
\mathfrak p_{\beta}\cap\mathfrak t\not=\{0\}\}$.  
Since $\mathfrak t$ and $\mathfrak t^{\perp}$ are ${\rm ad}(\mathfrak b)$-invariant, we have 
$$\mathfrak t^{\perp}=\mathfrak b\oplus\left(\mathop{\oplus}_{\beta\in(\triangle_{\mathfrak b})^V_+}
(\mathfrak p_{\beta}\cap\mathfrak t^{\perp})\right)$$
and 
$$\mathfrak t=\mathfrak z_{\mathfrak t}(\mathfrak b)\oplus\left(
\mathop{\oplus}_{\beta\in(\triangle_{\mathfrak b})^H_+}(\mathfrak p_{\beta}\cap\mathfrak t)\right),$$
where $\mathfrak z_{\mathfrak t}(\mathfrak b)$ denotes the centralizer of $\mathfrak b$ in $\mathfrak t$.  
For the convenience, we denote the centralizer $\mathfrak z_{\mathfrak p}(\mathfrak b)$ of $\mathfrak b$ in 
$\mathfrak p$ by $\mathfrak p_0$.  It is clear that 
$\mathfrak z_{\mathfrak t}(\mathfrak b)=\mathfrak p_0\cap\mathfrak t$.  
Let $\widetilde B_{Z_0}(F)$ be the partial tube over $F$ through $x_0$ stated in the previous section.  
According to Theorem 5.2, $M=\widetilde B_{Z_0}(F)$ holds and each fibre of this tube is the image by the 
normal exponential map of a principal orbit of an orthogonal representation on the normal space of $F$ given as 
the direct sum representation of some s-representations and a trivial representation, which implies 
that each fibre of this tube is a principal orbit of the isotropy action of the symmetric space 
${\rm Exp}(\mathfrak t^{\perp})$.  Take any $v\in T^{\perp}_{{\rm Exp}\,Z_0}M$.  
Then we have 
\begin{align*}
&R(v)|_{(\exp\,Z_0)_{\ast o}(\mathfrak p_{\beta})}=-\beta(v)^2\,{\rm id}\quad(\beta\in(\triangle_{\mathfrak b})_+
\cup\{0\}).
\end{align*}
According to (ii) of Proposition 5.3, we can derive that the horizontal lift 
(which is denoted by $(\mathfrak p_{\beta}\cap\mathfrak t)^L_{Z_0}$ of $\mathfrak p_{\beta}\cap\mathfrak t$ 
($\beta\in(\triangle_{\mathfrak b})_+^H\cup\{0\}$) to ${\rm Exp}\,Z_0$ is included by an eigenspace of $A_v$.  
According to (i) of Proposition 5.3 and the fact that each fibre of $M=\widetilde B_{Z_0}(F)$ is the image 
by the normal exponential map of a principal orbit of an orthogonal representation on the normal space of $F$ 
given as the direct sum representation of some s-representations and a trivial representation, we can derive that 
$(\exp\,Z_0)_{\ast}(\mathfrak p_{\beta}\cap\mathfrak t^{\perp})$ ($\beta\in(\triangle_{\mathfrak b})_+^V$) is 
included by an eigenspace of $A_v$.  
Also we have 
$$T_{{\rm Exp}\,Z_0}M=\left(\mathop{\oplus}_{\beta\in(\triangle_{\mathfrak b})^H_+\cup\{0\}}
(\mathfrak p_{\beta}\cap\mathfrak t)^L_{Z_0}\right)\oplus
\left(\mathop{\oplus}_{\beta\in(\triangle_{\mathfrak b})^V_+}
(\exp\,Z_0)_{\ast o}(\mathfrak p_{\beta}\cap\mathfrak t^{\perp})\right).$$
From the above facts, it follows that this decomposition is the common eigenspace decomposition of $A_v$ and $R(v)$.  
Hence $A_v$ and $R(v)$ commute.  It is clear that the same fact holds at other points of 
$M$.  Hence $M$ is curvature-adapted.  \qed

\section{Proof of Theorem B} 
In this section, we shall prove Theorem B.  
Let $M$ be as in Theorem B and $M^{\mathbb C}$ be the complete extrinsic complexification of $M$.  
See the proof of Theorem B in \cite{Koi7} about the construction of the complete extrinsic complexification of $M$.  
Let $\pi$ be the natural projection of $G^{\mathbb C}$ onto $G^{\mathbb C}/K^{\mathbb C}$ 
and $\phi:H^0([0,1],\mathfrak g^{\mathbb C})\to G^{\mathbb C}$ the parallel transport map for $G^{\mathbb C}$.  
Set $\widehat M^{\mathbb C}:=\pi^{-1}(M^{\mathbb C})$ and 
$\widetilde M^{\mathbb C}:=(\pi\circ\phi)^{-1}(M^{\mathbb C})$.  
Without loss of generality, we may assume that $K^{\mathbb C}$ is connected and 
that $G^{\mathbb C}$ is simply connected.  Hence $\widetilde M^{\mathbb C}$ is connected.  Denote by $A$ 
the shape tensor of $M$ and $R$ the curvature tensor of $G/K$.  First we shall show the following fact.  

\vspace{0.4truecm}

\noindent
{\bf Proposition 7.2.} {\sl The lifted submanifold $\widetilde M^{\mathbb C}$ is a full irreducible complete 
anti-Kaehler isoparametric submanifold with $J$-diagonalizable shape operators.}

\vspace{0.4truecm}

\noindent
{\it Proof.} Fix $x\in M$ and a unit normal vector $v$ of $M$ at $x$.  
Denote by ${\rm Spec}\,A_v$ and ${\rm Spec}\,R(v)$ the spectrum of $A_v$ and $R(v)$, repsectively.  
For each $(\lambda,\mu)\in({\rm Spec}\,A_v)\times({\rm Spec}\,R(v))$, set 
$$D_{\lambda\mu}:={\rm Ker}(A_v-\lambda\,{\rm id})\cap{\rm Ker}(R(v)-\mu\,{\rm id}).$$
Also, set 
$$\begin{array}{c}
{\mathcal S}:=\{(\lambda,\mu)\in{\rm Spec}\,A_v\times{\rm Spec}\,R(v)\,|\,D_{\lambda\mu}\not=\{0\}\}\,\,\,\,
{\mathcal S}_+:=\{(\lambda,\mu)\in{\mathcal S},\,|\,|\lambda|>{\sqrt{-\mu}}\}\\
{\rm and}\,\,\,\,{\mathcal S}_-:=\{(\lambda,\mu)\in{\mathcal S}\,|\,|\lambda|<{\sqrt{-\mu}}\}.
\end{array}$$
Since $M$ is curvature-adapted, $\displaystyle{T_xM=\mathop{\oplus}_{(\lambda,\mu)\in{\mathcal S}}D_{\lambda\mu}}$ 
holds.  For the simplicity, set 
$$Q(z):=\cos\left(z\sqrt{R(v)^{\mathbb C}}\right)-\frac{\sin\left(z\sqrt{R(v)^{\mathbb C}}\right)}
{\sqrt{R(v)^{\mathbb C}}}\circ A^{\mathbb C}_v.$$
Clearly we have 
$$Q(z)|_{D_{\lambda\mu}}=\left(\cos(\sqrt{-1}z\sqrt{-\mu})
-\frac{\lambda\sin(\sqrt{-1}z\sqrt{-\mu})}{\sqrt{-1}\sqrt{-\mu}}\right)
{\rm id}.$$
Hence, if $(\lambda,\mu)\in{\mathcal S}_+$ and $\mu\not=0$, then 
$\frac{1}{\sqrt{-\mu}}\left({\rm arctanh}\frac{\sqrt{-\mu}}{\lambda}
+k\pi\sqrt{-1}\right)$ ($k\in{\Bbb Z}$) are complex focal radii along $\gamma_v$ including $D_{\lambda\mu}$ 
as its nullity space.  Also, if $(\lambda,\mu)\in{\mathcal S}_-$ and $\mu\not=0$, then 
$\frac{1}{\sqrt{-\mu}}\left({\rm arctanh}\frac{\lambda}{\sqrt{-\mu}}+(k+\frac12)\pi\sqrt{-1}\right)$ ($k\in{\Bbb Z}$) 
are complex focal radii along $\gamma_v$ including $D_{\lambda\mu}$ as its nullity space.  
Also, if $\lambda\in{\rm Spec}\,A_v\setminus\{0\}$ satisfying $(\lambda,0)\in{\mathcal S}$, 
then $\frac{1}{\lambda}$ is a focal radii along $\gamma_v$ including $D_{\lambda0}$ as its nullity space.  
Also, if $|\lambda|=\sqrt{-\mu}$, then there exists no 
complex focal radius along $\gamma_v$ including $D_{\lambda\mu}$ as its nullity space.  
Hence, since $M$ satisfies the condition ($\ast_{\mathbb C}$), there exists no 
$(\lambda,\mu)\in{\mathcal S}$ satisfying $|\lambda|=\sqrt{-\mu}\not=0$.  
Thus we have 
$$T_xM=D_{00}\oplus\left(\mathop{\oplus}_{(\lambda,\mu)\in{\mathcal S}_+\cup
{\mathcal S}_-}D_{\lambda\mu}\right).\leqno{(7.1)}$$
Denote by $\widetilde A$ the shape tensor of $\widetilde M^{\mathbb C}$.  
Fix $u\in(\pi\circ\phi)^{-1}(x)$.  Let $v^L_u$ be the horizontal lift of $v$ to $u$.  
Then it follows from the above facts and Proposition 4 of \cite{Koi3} that 
$$\begin{array}{l}
\displaystyle{{\rm Spec}\,\widetilde A_{v_u^L}=\{\lambda\,|\,
\lambda\in{\rm Spec}A_v\,{\rm s.t.}\,(\lambda,0)\in{\mathcal S}_+\}}\\
\hspace{2.2truecm}\displaystyle{\bigcup\left\{\left.
\frac{\sqrt{-\mu}}{{\rm arctanh}\frac{\sqrt{-\mu}}{\lambda}+\pi k\sqrt{-1}}
\,\right|\,(\lambda,\mu)\in{\mathcal S}_+\,{\rm s.t.}\,\mu\not=0,\,\,\,\,k\in{\mathbb Z}\right\}}\\
\hspace{2.2truecm}\displaystyle{\bigcup\left\{\left.\frac{\sqrt{-\mu}}{{\rm arctanh}
\frac{\lambda}{\sqrt{-\mu}}+(k+\frac12)\pi\sqrt{-1}}\,\right|\,(\lambda,\mu)\in{\mathcal S}_-,
\,\,\,\,k\in{\mathbb Z}\right\}.}
\end{array}$$
For the simplicity, set 
$$\Lambda^+_{\lambda,\mu,k}:=\frac{\sqrt{-\mu}}
{{\rm arctanh}\frac{\sqrt{-\mu}}{\lambda}+\pi k\sqrt{-1}}\quad
((\lambda,\mu)\in{\mathcal S}_+\,\,{\rm s.t.}\,\,\mu\not=0,\,\,k\in{\Bbb Z})$$
and 
$$\Lambda^-_{\lambda,\mu,k}:=\frac{\sqrt{-\mu}}
{{\rm arctanh}\frac{\lambda}{\sqrt{-\mu}}+(k+\frac12)\pi\sqrt{-1}}\quad
((\lambda,\mu)\in{\mathcal S}_-,\,\,k\in{\Bbb Z}).$$
Also, we set 
$$\begin{array}{l}
\hspace{0.8truecm}\displaystyle{
\widetilde D_{\lambda}:={\rm Ker}\left(\widetilde A_{v_u^L}-\lambda\,{\rm id}\right)\,\,\,\,
(\lambda\in{\rm Spec}\,A_v\,\,{\rm s.t.}\,\,(\lambda,0)\in{\mathcal S}_+),}\\
\hspace{0.8truecm}\displaystyle{\widetilde D_{\Lambda^+_{\lambda,\mu,k}}
:={\rm Ker}\left(\widetilde A_{v_u^L}-\Lambda^+_{\lambda,\mu,k}\,{\rm id}\right)\,\,\,\,
((\lambda,\mu)\in{\mathcal S}_+\,\,{\rm s.t.}\,\,\mu\not=0,\,\,\,k\in{\Bbb Z})}\\
{\rm and}\quad\,\,\displaystyle{\widetilde D_{\Lambda^-_{\lambda,\mu,k}}
:={\rm Ker}\left(\widetilde A_{v_u^L}-\Lambda^-_{\lambda,\mu,k}\,{\rm id}\right)\,\,\,\,
((\lambda,\mu)\in{\mathcal S}_-,\,\,k\in{\Bbb Z}).}
\end{array}$$
Furthermore, by using $(7.1)$ and Lemma 9 of \cite{Koi3} (see Lemma 7.3 of \cite{Koi2} also), 
we can derive that $T_u\widetilde M^{\mathbb C}$ is equal to 
$$
\overline{
\left(\mathop{\oplus}_{\lambda\in{\rm Spec}\,A_v\,{\rm s.t.}\,(\lambda,0)
\in{\mathcal S}_+}\widetilde D_{\lambda}\right)\oplus\left(
\mathop{\oplus}_{(\lambda,\mu)\in{\mathcal S}_+\,{\rm s.t.}\,\mu\not=0}
\mathop{\oplus}_{k\in{\Bbb Z}}\widetilde D_{\Lambda^+_{\lambda,\mu,k}}\right)
\oplus\left(\mathop{\oplus}_{(\lambda,\mu)\in{\mathcal S}_-}
\mathop{\oplus}_{k\in{\Bbb Z}}\widetilde D_{\Lambda^-_{\lambda,\mu,k}}\right).}
$$
This implies that $\widetilde A_{v_u^L}$ is diagonalized with respect to a $J$-orthonormal base of 
$T_u\widetilde M^{\mathbb C}$.  Therefore it follows that $\widetilde M^{\mathbb C}$ is an anti-Kaehler 
isoparametric submanifold with $J$-diagonalizable shape operators from the arbitrarinesses of $x,v$ and $u$.  
Since $M$ is irreducible, it follows from Theorem 2 of \cite{Koi5} that 
the complex Coxeter group associated with $M$ is not decomposable, 
where we note that the complex Coxeter groups associated with $M$ is equal to 
the complex Coxeter groups associated with $\widetilde M^{\mathbb C}$ (see Introduction of \cite{Koi5}).  
Hence, it follows from Theorem 1 of \cite{Koi5} that $\widetilde M^{\mathbb C}$ is irreducible.  
Also, since $M$ is full, it is shown that the $J$-curvature normals of 
$\widetilde M^{\mathbb C}$ span the normal space of $\widetilde M^{\mathbb C}$ at each 
point of $\widetilde M^{\mathbb C}$ (see the discussion in the proof of Theorem 2 of \cite{Koi5}).  
Furthermore, it follows from this fact that $\widetilde M^{\mathbb C}$ is full (see the discussion in the proof 
of Theorem 1 of \cite{Koi5}).  The completeness of $\widetilde M^{\mathbb C}$ follows from the completeness of 
$M^{\mathbb C}$ and the fact that the fibres of $\pi\circ\phi$ are isometric to the complete anti-Kaehler Hilbert 
manifold $P(G^{\mathbb C},\{e\}\times K^{\mathbb C})$, where 
$P(G^{\mathbb C},\{e\}\times K^{\mathbb C})$ denotes the Hilbert Lie group 
$\{g\in H^1([0,1],G^{\mathbb C})\,|\,(g(0),g(1))\in\{e\}\times K^{\mathbb C}\}$ equipped with the natural 
complete anti-Kaehler metric.  This completes the proof.  \qed

\vspace{0.4truecm}

\noindent
{\it Remark 7.1.} According to this proposition, $M$ is proper complex equifocal in the sense of \cite{Koi4}.  

\vspace{0.4truecm}

Without loss of generality, we may assume 
$\hat 0\in\widetilde M^{\mathbb C}$ and hence $e\in\widehat M^{\mathbb C}$.  
For the simplicity, set $V:=H^0([0,1],\mathfrak g^{\mathbb C})$, 
$\langle\,\,,\,\,\rangle:=\langle\,\,,\,\,\rangle^A_0$ and 
$\langle\,\,,\,\,\rangle^{\mathcal I}:=\langle\,\,,\,\,\rangle^{\mathcal I}_0$.  
Also, denote by $||\cdot||$ the norm associated with 
$\langle\,\,,\,\,\rangle^{\mathcal I}$.  
Let ${\mathcal K}^h$ be the Lie algebra of all holomorphic Killing vector fields defined entirely on $V$ and 
${\mathcal K}^h_{\widetilde M^{\mathbb C}}$ the Lie subalgebra of ${\mathcal K}^h$ consisting of elements of ${\mathcal K}^h$ 
which are tangent to $\widetilde M^{\mathbb C}$ along $\widetilde M^{\mathbb C}$.  
Also, denote by $\mathfrak o_{AK}(V)$ be the Lie algebra of all continuous skew-symmetric complex linear maps 
from $V$ to oneself.  
Any $X\in{\mathcal K}^h$ is described as $X_u=Au+b\,\,(u\in V)$ for some 
$A\in\mathfrak o_{AK}(V)$ and some $b\in V$.  
Hence ${\mathcal K}^h$ is identified with $\mathfrak o_{AK}(V)\ltimes V$.  
Give $\mathfrak o_{AK}(V)$ the operator norm (which we denote by $||\cdot||_{\rm op}$) 
associated with $\langle\,\,,\,\,\rangle^{\mathcal I}$ and ${\mathcal K}^h$ the product norm of this norm 
$||\cdot||_{\rm op}$ of $\mathfrak o_{AK}(V)$ and the norm $||\cdot||$ of $V$.  
Then the space ${\mathcal K}^h$ is a Banach Lie algebra with respect to this norm.  
Let $I_h(V)$ be the group of all holomorphic isometries of $V$ and $I_h^b(V)$ be the subgroup of $I_h(V)$ 
generated by one-parameter transformation groups induced by elements of ${\mathcal K}^h$.  Since ${\mathcal K}^h$ 
is a Banach Lie algebra, $I_h^b(V)$ is a Banach Lie group.  Note that, for a general holomorphic isometry $f$ of $V$, 
$\displaystyle{\left.\frac{d}{dt}\right|_{t=0}(f_t)_{\ast}}$ is defined on a dense linear subspace of $V$ 
but it is not necessarily defined entirely on $V$ (see Example in Appendix of \cite{Koi13}).  
It is clear that ${\mathcal K}^h$ is the Banach Lie algebra of this Banach Lie group $I_h^b(V)$.  
Let $H_b$ be the closed Banach Lie subgroup of $I_h^b(V)$ of all elements of $I_h^b(V)$ preserving 
$\widetilde M^{\mathbb C}$ invariantly.  
From Fact 1.1 stated in Introduction and Proposition 7.2, we can derive the following extrinsic homogeneity theorem.  

\vspace{0.4truecm}

\noindent
{\bf Lemma 7.2.1.} {\sl We have $\widetilde M^{\mathbb C}=H_b\cdot\hat 0$.}

\vspace{0.4truecm}

Denote by $\rho$ the homomorphism from $H^1([0,1],G^{\mathbb C})$ to $I_h(V)$ defined by assigning 
$g\ast\cdot$ to each $g\in H^1([0,1],G^{\mathbb C})$ (i.e., $\rho(g)(u):=g\ast u$ 
($g\in H^1([0,1],G^{\mathbb C}),\,\,u\in V$)), where $g\ast u$ is as stated in Section 4.  

\vspace{0.4truecm}

\noindent
{\bf Lemma 7.2.2.} {\sl The group $\rho(H^1([0,1],G^{\mathbb C}))$ is a closed subgroup of $I_h^b(V)$.}

\vspace{0.4truecm}

\noindent
{\it Proof.} Take an arbitrary 
$v\in H^1([0,1],\mathfrak g^{\mathbb C})$ and set 
$\psi_s:=\rho(\exp\circ sv)$, where $\exp$ is the exponential map of 
$G^{\mathbb C}$.  Note that $\exp\circ sv$ is equal to the image of 
$sv\in H^1([0,1],\mathfrak g^{\mathbb C})$ by the exponential map of 
$H^1([0,1],G^{\mathbb C})$.  
The group $\{\psi_s\,|\,s\in{\mathbb R}\}$ is a one-parameter transformation group consisting of 
holomorphic isometries of $V$.  The holomorphic Killing vector field $X$ associated with 
$\{\psi_s\,|\,s\in{\mathbb R}\}$ is given by 
$$X_u=\left.\frac{d}{ds}\right|_{s=0}\psi_s(u)
=\left.\frac{d}{ds}\right|_{s=0}(\exp\circ s v)\ast u
={\rm ad}(v)(u)-v'.$$
Set $I_c:=\{t\in[0,1]\,|\,\max\,{\rm Spec}(-{\rm ad}(v(t))^2)\geq c\}$ and 
$c_0:=\min\{c\,|\,I_c\,\,{\rm is}\,\,{\rm of}\,\,{\rm measure}\,\,
{\rm zero}\,\,{\rm in}$\newline
$[0,1]\}$, where ${\rm ad}$ is the adjoint operator of $\mathfrak g^{\mathbb C}$.  Then we have 
$$\begin{array}{l}
\displaystyle{||{\rm ad}(v)u||^2
=\int_0^1\langle{\rm ad}(v(t))u(t),{\rm ad}(v(t))u(t)
\rangle^{\mathcal I}\,dt}\\
\hspace{1.65truecm}\displaystyle{=-\int_0^1\langle{\rm ad}(v(t))^2u(t),u(t)\rangle^{\mathcal I}\,dt
\leq c_0|| u||^2,}
\end{array}$$
where $\langle\,\,,\,\,\rangle^{\mathcal I}$ is 
the inner product of $\mathfrak g^{\mathbb C}$ stated in Section 2.  
Thus ${\rm ad}(v)$ is bounded.  Hence we have $X\in{\mathcal K}^h$, that is, 
$\rho(\exp\circ v)\in I_h^b(V)$.  
Therefore, it follows from the arbitrariness of $v$ that 
$\rho(H^1([0,1],G^{\mathbb C}))$ is a subgroup of $I_h^b(V)$.  
The closedness of $\rho(H^1([0,1],G^{\mathbb C}))$ is trivial.  
\qed

\vspace{0.3truecm}

In the proof of Theorem 7.1, it is key to show the following fact.  

\vspace{0.3truecm}

\noindent
{\bf Proposition 7.3.} {\sl The above group $H_b$ is a subgroup of 
$\rho(H^1([0,1],G^{\mathbb C}))$.}

\vspace{0.3truecm}

To prove this proposition, we prepare some lemmas.  
For $X\in{\mathcal K}^h$, 
we define a map $F_X:\Omega_e(G^{\mathbb C})\to\mathfrak g^{\mathbb C}$ by $F_X(g)
:=\phi_{\ast\widehat 0}((\rho(g)_{\ast}X)_{\widehat 0})$.  
For the simplicity, denote by ${\rm Ad}$ the adjoint operator 
${\rm Ad}_{G^{\mathbb C}}$ of $G^{\mathbb C}$.  
For this map $F_X$, we have the following fact.  

\vspace{0.3truecm}

\noindent
{\bf Lemma 7.3.1.} {\sl {\rm (i)} For $g\in\Omega_e(G^{\mathbb C})$, $F_X(g)
=\int_0^1{\rm Ad}(g)(X_{\rho(g^{-1})(\hat 0)})dt$.  

{\rm (ii)} If $X\in{\mathcal K}^h_{\widetilde M^{\mathbb C}}$, then the image of 
$F_X$ is included by $T_e\widehat M^{\mathbb C}$.}

\vspace{0.3truecm}

\noindent
{\it Proof.} Let $\{\psi_s\}_{s\in{\mathbb R}}$ be the one-parameter 
transformation group associated with $X$.  For each $g\in\Omega_e(G^{\mathbb C})$, 
we have 
$$\begin{array}{l}
\displaystyle{(\rho(g)_{\ast}X)_{\hat 0}=\left.\frac{d}{ds}
\right|_{s=0}\rho(g)(\psi_s(g^{-1}\ast\hat 0))}\\
\displaystyle{=\left.\frac{d}{ds}\right|_{s=0}({\rm Ad}(g)
(\psi_s(\rho(g^{-1})(\hat0)))-g'g_{\ast}^{-1})
={\rm Ad}(g)(X_{\rho(g^{-1})(\hat 0)}).}
\end{array}$$
Also we have $\phi_{\ast\hat0}(u)=\int_0^1u(t)dt$ 
($u\in T_{\hat0}V(=V)$) (see Lemma 6 of \cite{Koi3}).  
Hence we obtain the relation in (i).  
Since $g\in\Omega_e(G^{\mathbb C})$, it maps each fibre of $\phi$ to 
oneself.  Hence, if $X\in{\mathcal K}^h_{\widetilde M^{\mathbb C}}$, then 
$\rho(g)_{\ast}X\in{\mathcal K}^h_{\widetilde M^{\mathbb C}}$.  In particular, 
we have $(\rho(g)_{\ast}X)_{\hat0}\in T_{\hat0}\widetilde M^{\mathbb C}$.  
Therefore we obtain $F_X(g)\in\phi_{\ast\hat0}(T_{\hat0}\widetilde M^{\mathbb C})
=T_e\widehat M^{\mathbb C}$.  \qed

\vspace{0.3truecm}

For $v\in H^1([0,1],\mathfrak g^{\mathbb C})$, we define a vector field $X^v$ on 
$V$ by $(X^v)_u:=[v,u]-v'$ ($u\in V$).  
Let $\{\exp\circ sv\,|\,s\in{\mathbb R}\}$ be the one-parameter 
subgroup of $H^1([0,1],G^{\mathbb C})$ associated with $v$.  Then the holomorphic 
Killing vector field associated with the one-parameter transformation group 
$\{\rho(\exp\circ sv)\,|\,s\in{\mathbb R}\}$ of $V$ is equal to $X^v$.  
Furthermore, we can show $X^v\in{\mathcal K}^h_b$ by the discussion in the proof 
of Lemma 7.3.  For $X^v$, we have the following fact.  

\vspace{0.4truecm}

\noindent
{\bf Lemma 7.3.2.} {\sl The map $F_{X^v}$ is a constant map.}

\vspace{0.4truecm}

\noindent
{\it Proof.} Take elements $g_1$ and $g_2$ of $\Omega_e(G^{\mathbb C})$.  
Since $\rho(g_i)$ maps each fibre of $\phi$ to oneself by the fact (iii) for 
$\phi$ (stated in Section 2), we have 
$\phi\circ\rho(g_i)=\phi$ ($i=1,2$) and hence 
$$F_{X^v}(g_i)=\phi_{\ast\hat0}((\rho(g_i)_{\ast}(X^v))_{\hat0})
=\phi_{\ast\rho(g_i^{-1})(\hat0)}((X^v)_{\rho(g_i^{-1})(\hat0)})\quad(i=1,2).
\leqno{(7.2)}$$
Since $\rho(\exp\circ sv)$ maps the fibres of $\phi$ to them by the fact (iii) for $\phi$ and 
$\phi(\rho(g_1^{-1})(\hat0))=\phi(\rho(g_2^{-1})(\hat0))$, we have 
$\phi(\rho(\exp\,sv)(\rho(g_1^{-1})(\hat0)))
=\phi(\rho(\exp\,sv)(\rho(g_2^{-1})(\hat0)))$ and hence 
$\phi_{\ast\rho(g_1^{-1})(\hat0)}(X^v_{\rho(g_1^{-1})(\hat0)})
=\phi_{\ast\rho(g_2^{-1})(\hat0)}(X^v_{\rho(g_2^{-1})(\hat0)})$.  
From this relation and $(7.2)$, we obtain $F_{X^v}(g_1)=F_{X^v}(g_2)$.  
Therefore it follows from the arbitrarinesses of $g_1$ and $g_2$ that 
$F_{X^v}$ is a constant map.  \qed

\vspace{0.3truecm}

For each $u\in V$, denote by $\widetilde u$ the element 
$t\mapsto\int_0^tu(t)dt\,\,\,(0\leq t\leq1)$ of 
$H^1([0,1],\mathfrak g^{\mathbb C})$.  

Also we have the following fact for $F_X$.  

\vspace{0.3truecm}

\noindent
{\bf Lemma 7.3.3.} {\sl {\rm (i)} The map $X\mapsto F_X$ is linear.

{\rm (ii)} $F_X(g_1g_2)=F_{\rho(g_2)_{\ast}X}(g_1)$ 
($g_1,g_2\in\Omega_e(G^{\mathbb C})$).  

{\rm (iii)} $(dF_X)_g\circ(dR_g)_{\hat e}=(dF_{\rho(g)_{\ast}X})_{\hat e}$ 
($g\in\Omega_e(G^{\mathbb C})$).  

{\rm (iv)} If $X_u=Au+b$ ($u\in V$) for some linear transformation $A$ of $V$ 
and some $b\in V$, then we have 
$(dF_X)_{\hat e}(u)=\int_0^1(A+{\rm ad}(\widetilde b))u'dt$ 
($u\in\Omega_0(\mathfrak g^{\mathbb C})$), where ${\rm ad}$ is the adjoint 
representation of $\mathfrak g^{\mathbb C}$ and 
$\Omega_0(\mathfrak g^{\mathbb C}):=\{u\in H^1([0,1],\mathfrak g^{\mathbb C})\,|\,
u(0)=u(1)=0\}$.  

{\rm (v)} If $X,\overline X\in{\mathcal K}^h$ and if $\overline X-X=X^v$ for some 
$v\in H^1([0,1],\mathfrak g^{\mathbb C})$, then $F_{\overline X}-F_X$ is 
a constant map.}

\vspace{0.3truecm}

\noindent
{\it Proof.} The statements ${\rm (i)}\sim{\rm (iii)}$ are trivial.  
The statement (iv) is shown by imitating the proof of Proposition 2.3 of 
\cite{Ch}.  The statement (v) follows from Lemma 7.3.2 and (i) directly.  
\qed

\vspace{0.3truecm}

By imitating the proof of Theorem 2.2 of \cite{Ch}, we can show the following fact 
in terms of Lemmas 7.3.1$\sim$7.3.3.  

\vspace{0.3truecm}

\noindent
{\bf Lemma 7.3.4.} {\sl Let $X$ be an element of ${\mathcal K}^h$ given by 
$X_u:=[v,u]-b$ ($u\in V$) for some $v,b\in V$.  If 
$X\in{\mathcal K}^h_{\widetilde M^{\mathbb C}}$, then we have $v\in 
H^1([0,1],\mathfrak g^{\mathbb C})$ and $b=v'$ (i.e., $X=X^v$).}

\vspace{0.3truecm}

\noindent
{\it Proof.} Set $\overline X:=X-X^{\widetilde b}$ and $w:=v-\widetilde b$.  
First we consider the case where $G^{\mathbb C}$ is simple.  
From $\overline X={\rm ad}(w)$, we have 
$$(\rho(g)_{\ast}\overline X)_u=\rho(g)_{\ast}(\overline X_{\rho(g^{-1})(u)})
={\rm Ad}(g)([w,\rho(g^{-1})(u)])=[{\rm Ad}(g)w,\,u-g\ast\hat0]\,\,\,\,
(u\in V).$$
From this relation and (i) of Lemma 7.3.1, we have 
$$\begin{array}{l}
\hspace{0.5truecm}\displaystyle{(dF_{\rho(g)_{\ast}\overline X})_{\hat e}(u)
=\left.\frac{d}{ds}\right|_{s=0}F_{\rho(g)_{\ast}\overline X}(\exp\,su)}\\
\displaystyle{=\left.\frac{d}{ds}\right|_{s=0}\int_0^1{\rm Ad}(\exp\,su)
((\rho(g)_{\ast}\overline X)_{\rho(\exp(-su))(\hat0)})dt}\\
\displaystyle{=\int_0^1\left([u,\,(\rho(g)_{\ast}\overline X)_{\hat 0}]
+\left.\frac{d}{ds}\right|_{s=0}(\rho(g)_{\ast}\overline X)
_{\rho(\exp(-su))(\hat 0)}\right)dt}\\
\displaystyle{=\int_0^1\left(-[u,\,[{\rm Ad}(g)w,g\ast\hat 0]]
+\left.\frac{d}{ds}\right|_{s=0}[{\rm Ad}(g)w,
\rho(\exp(-su))(\hat 0)-g\ast\hat 0]\right)dt}\\
\displaystyle{=\int_0^1\left([u,\,[{\rm Ad}(g)w,g'g_{\ast}^{-1}]]
-[{\rm Ad}(g)w,\,\left.\frac{d}{ds}\right|_{s=0}
((\exp(-su))'\exp(-su)_{\ast}^{-1})]\right)dt}\\
\displaystyle{=\int_0^1\left([u,\,[{\rm Ad}(g)w,g'g_{\ast}^{-1}]]
+[{\rm Ad}(g)w,\,u']\right)dt}
\\
\displaystyle{=\left.
[u(t),\,\widetilde{[{\rm Ad}(g)w,g'g_{\ast}^{-1}]}(t)]\right|_{t=1}
-\left.[u(t),\,\widetilde{[{\rm Ad}(g)w,g'g_{\ast}^{-1}]}(t)]\right|_{t=0}
}\\
\hspace{0.5truecm}\displaystyle{
-\int_0^1[u',\,\widetilde{[{\rm Ad}(g)w,g'g_{\ast}^{-1}]}]dt
+\int_0^1[{\rm Ad}(g)w,\,u']dt}\\
\displaystyle{=
\int_0^1[\widetilde{[{\rm Ad}(g)w,g'g_{\ast}^{-1}]}+{\rm Ad}(g)w,\,\,u']dt}
\end{array}
\leqno{(7.3)}
$$
for $u\in T_{\hat e}(\Omega_e(G^{\mathbb C}))
(=\Omega_0(\mathfrak g^{\mathbb C}))$, 
where each of the notation $\,{}'\,\,$ means the derivative 
with respect to $t$, $\hat e$ is the constant path at the identity element 
$e$ of $G^{\mathbb C}$ and 
$\Omega_0(\mathfrak g^{\mathbb C}):=\{u\in H^1([0,1],\mathfrak g^{\mathbb C})\,|\,
u(0)=u(1)=0\}$.  
According to (ii) of Lemma 7.3.1, 
we have ${\rm Im}\,F_X\subset T_e\widehat M^{\mathbb C}$ and hence 
${\rm dim}_{\mathbb C}({\rm Span}_{\mathbb C}{\rm Im}\,F_X)\leq{\rm dim}_{\mathbb C}
T_e\widehat M^{\mathbb C}\leq{\rm dim}_{\mathbb C}\mathfrak g^{\mathbb C}-2$, 
where ${\rm Span}_{\mathbb C}(\cdot)$ means the complex linear span of $(\cdot)$ 
and ${\rm dim}_{\mathbb C}(\cdot)$ means the complex dimension of $(\cdot)$.  
Since $F_{\overline X}-F_X$ is a constant map by (v) of Lemma 7.3.3, we have 
${\rm dim}_{\mathbb C}({\rm Span}_{\mathbb C}{\rm Im}\,F_{\overline X})\leq
{\rm dim}_{\mathbb C}\mathfrak g^{\mathbb C}-1$, that is, 
${\rm dim}_{\mathbb C}(\mathfrak g^{\mathbb C}\ominus{\rm Span}_{\mathbb C}
{\rm Im}\,F_{\overline X})\geq1$.  
Take $Y(\not=0)\in\mathfrak g^{\mathbb C}\ominus
{\rm Span}_{\mathbb C}\,{\rm Im}\,F_{\overline X}$.  Also, take 
$g\in\Omega_e(G^{\mathbb C})$ and $u\in T_{\hat e}(\Omega_e(G^{\mathbb C}))$.  
By using (iii) of Lemma 7.3.3 and $(7.3)$, we have 
$$\begin{array}{l}
\hspace{0.5truecm}\displaystyle{
\langle(dF_{\overline X})_g((dR_g)_{\hat e}(u)),Y\rangle^A
=\langle(dF_{\rho(g)_{\ast}\overline X})_{\hat e}(u),Y\rangle^A}\\
\displaystyle{=\langle
\int_0^1[\widetilde{[{\rm Ad}(g)w,g'g_{\ast}^{-1}]}+{\rm Ad}(g)w,\,\,u']dt,\,\,
Y\rangle^A}\\
\displaystyle{=\int_0^1
\langle[\widetilde{[{\rm Ad}(g)w,g'g_{\ast}^{-1}]}+{\rm Ad}(g)w,\,\,u'],\,\,
Y\rangle^A\,dt}\\
\displaystyle{=-\int_0^1
\langle u',\,\,[\widetilde{[{\rm Ad}(g)w,g'g_{\ast}^{-1}]}+{\rm Ad}(g)w,\,\,
Y]\rangle^A\,dt}\\
\displaystyle{=-\langle u',\,\,[\widetilde{[{\rm Ad}(g)w,g'g_{\ast}^{-1}]}
+{\rm Ad}(g)w,\,\,Y]\rangle,}
\end{array}$$
where $\langle\,\,\,\,\rangle^A$ is the non-degenerate symmetric bilinear form of $\mathfrak g^{\mathbb C}$ 
stated in Section 4.  
For the simplicity, we set $\eta:=\widetilde{[{\rm Ad}(g)w,g'g_{\ast}^{-1}]}+{\rm Ad}(g)w$.  
On the other hand, from $(dF_{\overline X})_g((dR_g)_{\hat e}(u))
\in{\rm Span}_{\mathbb C}{\rm Im}\,F_{\overline X}$, we have 
$\langle(dF_{\overline X})_g((dR_g)_{\hat e}(u)),Y\rangle^A=0$.  
Hence we have $\langle u',\,[\eta,Y]\rangle=0$.  
The space $\Omega_0(\mathfrak g^{\mathbb C})$ is identified with the vertical 
space (which is denoted by ${\mathcal V}_{\hat 0}$) at $\hat 0$ of $\phi$ under 
the correspondence $u\mapsto u'$ ($u\in\Omega_0(\mathfrak g^{\mathbb C})$), where 
we note that $\phi_{\ast\hat0}(u')=\int_0^1u'(t)\,dt=0$ by Lemma 6 of \cite{Koi3} 
(hence $u'\in{\mathcal V}_{\hat 0}$).  
Hence, from the arbitrariness of $u$, it follows that $[\eta,Y]$ belongs to 
the horizontal space (which is denoted by ${\mathcal H}_{\hat 0}$) at $\hat 0$ of 
$\phi$.  Since $G^{\mathbb C}$ has no center, 
there exists $Z\in\mathfrak g^{\mathbb C}$ with $[Y,Z]\not=0$.  Set $W:=[Y,Z]$.  
By using Lemma 6 of \cite{Koi3}, we can show that ${\mathcal H}_{\hat 0}$ is equal to 
the set of all constant paths in $\mathfrak g^{\mathbb C}$.  Hence it follows from 
$[\eta,Y]\in{\mathcal H}_{\hat 0}$ that $[\eta,Y]$ is a constant path.  
Furthermore it follows from $\langle\eta,W\rangle^A=\langle[\eta,Y],Z
\rangle^A$ that $\langle\eta,W\rangle^A$ is constant, that is, 
$$\langle\widetilde{[{\rm Ad}(g)w,g'g_{\ast}^{-1}]},W\rangle^A
+\langle{\rm Ad}(g)w,W\rangle^A={\rm const.}\leqno{(7.4)}$$
Since $\mathfrak g^{\mathbb C}$ has no center, there exists 
$\overline W\in\mathfrak g^{\mathbb C}$ with $[W,\overline W]\not=0$.  Since 
$G^{\mathbb C}$ is simple, ${\rm Ad}(G^{\mathbb C})[W,\overline W]$ is full in $\mathfrak g^{\mathbb C}$.  
Hence there exist $h_1,\cdots,h_{2m}\in G^{\mathbb C}$ such that $({\rm Ad}(h_1)[W,\overline W],\cdots,
{\rm Ad}(h_{2m})[W,\overline W])$ is a base of $\mathfrak g^{\mathbb C}$ (regarded as a real vector space), 
where $m:={\rm dim}_{\mathbb C}G^{\mathbb C}$.  For a sufficiently small $\varepsilon>0$, we take 
$g_i\in\Omega_e(G^{\mathbb C})$ with $g_i|_{(\varepsilon,1-\varepsilon)}=h_i$ ($i=1,\cdots,2m$).  
Since $g_i$ ($i=1,\cdots,2m$) are constant over $[\varepsilon,1-\varepsilon]$, 
it follows form $(7.4)$ ($g=g_i$-case) that 
$\langle w,{\rm Ad}(h_i^{-1})W\rangle^A$ ($i=1,\cdots,2m$) are constant over 
$[\varepsilon,1-\varepsilon]$.  
Hence $w$ is constant over $[\varepsilon,1-\varepsilon]$.  
Hence it follows from the arbitrariness of $\varepsilon$ that $w$ is constant 
over $[0,1]$.  That is, we obtain $b=v'$ and hence 
$v\in H^1([0,1],\mathfrak g^{\mathbb C})$.  

Next we consider the case where $G^{\mathbb C}$ is not simple.  
Let $G^{\mathbb C}=G_1^{\mathbb C}\times\cdots\times G_k^{\mathbb C}$ be the irreducible 
decomposition of $G^{\mathbb C}$ and $\mathfrak g_i^{\mathbb C}$ be the Lie algebra of 
$G_i^{\mathbb C}$ ($i=1,\cdots,k$).  Let $\mathfrak g_{\overline X}^{\mathbb C}$ be 
the maximal 
ideal of $\mathfrak g^{\mathbb C}$ such that the orthogonal projection of 
$w=v-\widetilde b$ 
onto the ideal is a constant path, where we note that any ideal of 
$\mathfrak g^{\mathbb C}$ is equal to the direct sum of some 
$\mathfrak g^{\mathbb C}_i$'s and hence it is a non-degenerate subspace with 
respect to $\langle\,\,,\,\,\rangle^A$.  
Now we shall show 
$$(\mathfrak g_{\overline X}^{\mathbb C})^{\perp}\subset 
T_e\widehat M^{\mathbb C},\leqno{(7.5)}$$
where $(\mathfrak g^{\mathbb C}_{\overline X})^{\perp}$ is the orthogonal 
complement of $\mathfrak g^{\mathbb C}_{\overline X}$ in $\mathfrak g^{\mathbb C}$ with respect to 
$\langle\,\,,\,\,\rangle^A$.  
Let $V_i:=H^0([0,1],\mathfrak g_i^{\mathbb C})$ ($i=1,\cdots,k$).  
It is clear that $V=V_1\oplus\cdots\oplus V_k$ (orthogonal direct sum).  
The holomorphic Killing vector field $\overline X$ is 
described as $\overline X=\overline X_1^L+\cdots+\overline X_k^L$ in terms of 
some holomorphic Killing vector field $\overline X_i$ on $V_i$ ($i=1,\cdots,k$), 
where $\overline X_i^L$ is the holomorphic Killing vector field on $V$ defined by 
$(\overline X_i^L)_u=(\overline X_i)_{u_i}$ ($u=(u_1,\cdots,u_k)\in V$).  
For $g=(g_1,\cdots,g_k)\in\Omega_e(G^{\mathbb C})
(=\Omega_e(G_1^{\mathbb C})\times\cdots\times\Omega_e(G_k^{\mathbb C}))$, we have 
${\rm Ad}(g)(\overline X_{\rho(g^{-1})(\hat0)})=\sum\limits_{i=1}^k
{\rm Ad}_i(g_i)((\overline X_i)_{\rho_i(g_i^{-1})(\hat0)})$, where 
${\rm Ad}_i$ denotes the adjoint representation of $G^{\mathbb C}_i$ and $\rho_i$ denotes 
the homomorphism from $H^1([0,1],G_i^{\mathbb C})$ to $I_h(V_i)$ defined 
in similar to $\rho$.  Hence, from (i) of Lemma 7.3.1, we have 
$F_{\overline X}(g)=\sum\limits_{i=1}^kF^i_{\overline X_i}(g_i)$, where 
$F^i_{\overline X_i}$ is the map from $\Omega_e(G_i^{\mathbb C})$ to 
$\mathfrak g_i^{\mathbb C}$ defined in similar to $F_{\overline X}$.  
Therefore we obtain 
$\displaystyle{{\rm Span}_{\mathbb C}{\rm Im}\,F_{\overline X}
=\mathop{\oplus}_{i=1}^k{\rm Span}_{\mathbb C}{\rm Im}\,F^i_{\overline X_i}}$.  
Let $v=\sum\limits_{i=1}^kv_i$ and $\widetilde b=\sum\limits_{i=1}^k
\widetilde b_i$, where $v_i,\widetilde b_i\in V_i$ ($i=1,\cdots,k$).  
Since $\mathfrak g^{\mathbb C}_{\overline X}$ is an ideal of 
$\mathfrak g^{\mathbb C}$, it is described as 
$\displaystyle{\mathfrak g^{\mathbb C}_{\overline X}=\mathop{\oplus}_{i\in I}
\mathfrak g^{\mathbb C}_i}$ ($I\subset\{1,\cdots,k\}$).  
Since $v_i-\widetilde b_i$ ($i\in I$) are constant paths by the definition of 
$\mathfrak g^{\mathbb C}_{\overline X}$, 
${\rm Ad}_i(g_i)[v_i-\widetilde b_i,\rho(g_i^{-1})(\hat 0)]$ ($i\in I$) are 
loops and hence 
$$F^i_{\overline X_i}(g_i)=\int_0^1{\rm Ad}_i(g_i)
[v_i-\widetilde b_i,\rho(g_i^{-1})(\hat 0)]dt=0\qquad(i\in I).$$
Hence we have 
$${\rm Span}_{\mathbb C}{\rm Im}\,F_{\overline X}\subset
(\mathfrak g^{\mathbb C}_{\overline X})^{\perp}
(=\mathop{\oplus}_{i\notin I}\mathfrak g^{\mathbb C}_i).$$
Also we can show 
${\rm Span}_{\mathbb C}{\rm Im}\,F^i_{\overline X_i}=\mathfrak g^{\mathbb C}_i\,\,\,\,
(i\notin I)$.  Therefore we obtain 
$${\rm Span}_{\mathbb C}{\rm Im}\,F_{\overline X}
=(\mathfrak g^{\mathbb C}_{\overline X})^{\perp}.\leqno{(7.6)}$$
Also, since $F_{\overline X}-F_X$ is a constant map by (v) of Lemma 7.3.3 
and $0\in{\rm Im}\,F_{\overline X}$, we have 
$${\rm Span}_{\mathbb C}{\rm Im}\,F_{\overline X}\subset{\rm Span}_{\mathbb C}
{\rm Im}\,F_X.\leqno{(7.7)}$$
From $(7.6),\,(7.7)$ and (ii) of Lemma 7.3.1, we obtain 
$(\mathfrak g_{\overline X}^{\mathbb C})^{\perp}\subset T_e\widehat M^{\mathbb C}$.  
Next we shall show that $(R_g)_{\ast}((\mathfrak g_{\overline X}^{\mathbb C})
^{\perp})\subset T_g\widehat M^{\mathbb C}$ for any $g\in\widehat M^{\mathbb C}$.  
Fix $g\in\widehat M^{\mathbb C}$.  Define $\widehat g\in H^1([0,1],G^{\mathbb C})$ 
with $\widehat g(0)=e$ and $\widehat g(1)=g$ by $\widehat g(t)
:=\exp tY$ for some $Y\in\mathfrak g^{\mathbb C}$.  
Since $\phi\circ\rho(\widehat g)=R_g^{-1}\circ\phi$, we have 
$\phi^{-1}(R_g^{-1}(\widehat M^{\mathbb C}))=\rho(\widehat g)
(\widetilde M^{\mathbb C})$.  
Also we have $\rho(\widehat g)_{\ast}X\in
{\mathcal K}^h_{\rho(\widehat g)(\widetilde M^{\mathbb C})}$.  
Hence, by imitating the above discussion, we can show 
$$(\mathfrak g^{\mathbb C}_{\overline{\rho(\widehat g)_{\ast}X}})^{\perp}\subset 
T_eR_g^{-1}(\widehat M^{\mathbb C})=(R_g)_{\ast}^{-1}(T_g\widehat M^{\mathbb C}).
\leqno{(7.8)}$$
Also, we have 
$$(\rho(\widehat g)_{\ast}X)_u=\rho(\hat g)_{\ast}
(X_{\rho(\hat g)^{-1}(u)})
=[{\rm Ad}(\widehat g)v,u]-[{\rm Ad}(\widehat g)v,\rho(\widehat g)(\hat 0)]
-{\rm Ad}(\widehat g)b.\leqno{(7.9)}$$
Set $\overline v:={\rm Ad}(\widehat g)v$ and 
$\overline b:=[{\rm Ad}(\widehat g)v,\rho(\widehat g)(\hat 0)]
+{\rm Ad}(\widehat g)b$.  
Denote by ${\rm pr}_{\mathfrak g^{\mathbb C}_{\overline X}}$ the orthogonal 
projection of $\mathfrak g^{\mathbb C}$ onto $\mathfrak g^{\mathbb C}_{\overline X}$.  
Since ${\rm Ad}(\widehat g)$ preserves each $\mathfrak g_i^{\mathbb C}$ 
invariantly, it preserves $\mathfrak g^{\mathbb C}_{\overline X}$ and 
$(\mathfrak g^{\mathbb C}_{\overline X})^{\perp}$ invariantly, respectively.  
Hence we have 
${\rm pr}_{\mathfrak g^{\mathbb C}_{\overline X}}\circ{\rm Ad}(\widehat g)={\rm Ad}
(\widehat g)\circ{\rm pr}_{\mathfrak g^{\mathbb C}_{\overline X}}$ and 
${\rm pr}_{\mathfrak g^{\mathbb C}_{\overline X}}\circ{\rm ad}(Y)={\rm ad}(Y)\circ
{\rm pr}_{\mathfrak g^{\mathbb C}_{\overline X}}$.  Also, we have 
$\rho(\widehat g)(\hat0)=-Y=-{\rm Ad}(\widehat g)Y$.  By using these facts and 
noticing that ${\rm pr}_{\mathfrak g^{\mathbb C}_{\overline X}}(v-\widetilde b)$ 
is a constant path, we have 
$$\begin{array}{l}
\hspace{0.7truecm}\displaystyle{\frac{d}{dt}
{\rm pr}_{\mathfrak g^{\mathbb C}_{\overline X}}
\left(\overline v-\widetilde{\overline b}\right)}\\
\displaystyle{=\frac{d}{dt}{\rm pr}_{\mathfrak g^{\mathbb C}_{\overline X}}\left(
{\rm Ad}(\widehat g)(v-\widetilde b)+{\rm Ad}(\widehat g)\widetilde b
-\widetilde{[{\rm Ad}(\widehat g)v,\rho(\widehat g)(\hat0)]}
-\widetilde{{\rm Ad}(\widehat g)b}\right)}\\
\displaystyle{={\rm Ad}(\widehat g)
[Y,{\rm pr}_{\mathfrak g^{\mathbb C}_{\overline X}}
(v-\widetilde b)]+{\rm Ad}(\widehat g)
[Y,{\rm pr}_{\mathfrak g^{\mathbb C}_{\overline X}}
(\widetilde b)]}\\
\hspace{0.7truecm}
\displaystyle{+{\rm Ad}(\widehat g){\rm pr}_{\mathfrak g^{\mathbb C}_{\overline X}}
(b)+{\rm pr}_{\mathfrak g^{\mathbb C}_{\overline X}}[{\rm Ad}(\widehat g)v,Y]-
{\rm Ad}(\widehat g){\rm pr}_{\mathfrak g^{\mathbb C}_{\overline X}}(b)}\\
\displaystyle{=({\rm pr}_{\mathfrak g^{\mathbb C}_{\overline X}}\circ
{\rm Ad}(\widehat g))\left([Y,v-\widetilde b]+[Y,\widetilde b]+[v,Y]\right)=0.}
\end{array}$$
Thus ${\rm pr}_{\mathfrak g^{\mathbb C}_{\overline X}}
(\overline v-\widetilde{\overline b})$ 
is a constant path.  This fact together with $(7.9)$ implies 
$\mathfrak g^{\mathbb C}_{\overline X}\subset
\mathfrak g^{\mathbb C}_{\overline{\rho(\widehat g)_{\ast}X}}$.  
By exchanging the roles 
of $X$ and $\rho(\widehat g)_{\ast}X$, we have 
$\mathfrak g^{\mathbb C}_{\overline{\rho(\widehat g)_{\ast}X}}\subset
\mathfrak g^{\mathbb C}_{\overline X}$.  Thus we obtain 
$\mathfrak g^{\mathbb C}_{\overline X}
=\mathfrak g^{\mathbb C}_{\overline{\rho(\widehat g)_{\ast}X}}$.  
Therefore 
the relation $(R_g)_{\ast}(\mathfrak g^{\mathbb C}_{\overline X})^{\perp}\subset 
T_g\widehat M^{\mathbb C}$ follows from $(7.8)$.  
Since this relation holds for any $g\in\widehat M^{\mathbb C}$ and 
$\mathfrak g^{\mathbb C}_{\overline X}$ is an ideal of $\mathfrak g^{\mathbb C}$, 
we have 
$\widehat M^{\mathbb C}=\widehat M^{{\mathbb C}'}\times G^{{\mathbb C}\perp}_{\overline X}
\subset G^{\mathbb C}_{\overline X}\times G_{\overline X}^{{\mathbb C}\perp}
(=G^{\mathbb C})$ for some submanifold $\widehat M^{\mathbb C'}$ in 
$G^{\mathbb C}_{\overline X}$, where 
$G^{\mathbb C}_{\overline X}:=\exp(\mathfrak g^{\mathbb C}_{\overline X})$ 
and 
$G_{\overline X}^{{\mathbb C}\perp}:=\exp
((\mathfrak g^{\mathbb C}_{\overline X})^{\perp})$.  
Since $\widehat M^{\mathbb C}$ is irreducible and 
${\rm dim}\,\widehat M^{\mathbb C}<{\rm dim}\,G^{\mathbb C}$, we have 
$(\mathfrak g^{\mathbb C}_{\overline X})^{\perp}=\{0\}$, that is, 
$\mathfrak g^{\mathbb C}_{\overline X}=\mathfrak g^{\mathbb C}$.  This implies that 
$v-\widetilde b$ is a constant path.  Therefore we obtain $b=v'$ and hence 
$v\in H^1([0,1],\mathfrak g^{\mathbb C})$.  
\qed

\vspace{0.5truecm}

Also we have the following fact.  

\vspace{0.5truecm}

\noindent
{\bf Lemma 7.3.5.} {\sl The set ${\mathcal K}^h_{\widetilde M^{\mathbb C}}$ is 
closed in ${\mathcal K}^h$.}

\vspace{0.5truecm}

\noindent
{\it Proof.} Denote by $\overline{{\mathcal K}^h_{\widetilde M^{\mathbb C}}}$ the 
closure of ${\mathcal K}^h_{\widetilde M^{\mathbb C}}$ in ${\mathcal K}^h$.  
Take $X\in\overline{{\mathcal K}^h_{\widetilde M^{\mathbb C}}}$.  Then there exists 
a sequence $\{X_n\}_{n=1}^{\infty}$ in ${\mathcal K}^h_{\widetilde M^{\mathbb C}}$ 
with $\lim\limits_{n\to\infty}X_n=X$ (in ${\mathcal K}^h$).  Let 
$(X_n)_u=A_nu+b_n$ ($A_n\in\mathfrak o_{AK}(V),\,\,b_n\in V$) and 
$X_u=Au+b$ ($A\in\mathfrak o_{AK}(V),\,\,b\in V$).  
From $\lim\limits_{n\to\infty}X_n=X$ 
(in ${\mathcal K}^h$), we have 
$\lim\limits_{n\to\infty}A_n=A$ 
(in ${\mathfrak o}_{AK}(V)$) and hence 
$\lim\limits_{n\to\infty}A_nu=Au$ ($u\in V$).  Also, we have 
$\lim\limits_{n\to\infty}b_n=b$.  Hence we have 
$\lim\limits_{n\to\infty}(X_n)_u=X_u$ ($u\in V$).  For each 
$u\in\widetilde M^{\mathbb C}$, denote by ${\rm pr}^{\perp}_u$ the orthogonal 
projection of $V$ onto $T^{\perp}_u\widetilde M^{\mathbb C}$.  Since 
${\rm dim}\,T^{\perp}_u\widetilde M^{\mathbb C}<\infty$, 
${\rm pr}^{\perp}_u$ is a compact operator.  Hence, since 
${\rm pr}^{\perp}_u((X_n)_u)=0$ for all $n$, we obtain 
${\rm pr}_u^{\perp}(X_u)=0$ and hence 
$X\in{\mathcal K}^h_{\widetilde M^{\mathbb C}}$.  Therefore we obtain 
$\overline{{\mathcal K}^h_{\widetilde M^{\mathbb C}}}
={\mathcal K}^h_{\widetilde M^{\mathbb C}}$.  
\qed

\vspace{0.5truecm}

Take $v\in V$ and $X\in{\mathcal K}^h$.  Also, define 
$g_n\in H^1([0,1],G^{\mathbb C})$ ($n\in{\bf N}$) by $g_n(t):=\exp(n\widetilde v(t))$ and 
a vector field $X_n^v$ ($n\in{\bf N}$) by $X_n^v:=\frac 1n\rho(g_n)_{\ast}X$.  
Since $\rho(g_n)\in I_h^b(V)$ by Lemma 7.2.2, we have $X_n^v\in{\mathcal K}^h$.  
Let $X_u=Au+b$ ($A\in{\mathfrak o}_{AK}(V),\,\,b\in V$), where $u\in V$, and 
$(X_n^v)_u=A_n^vu+b_n^v$ ($A_n^v\,:\,$ a skew-symmetric complex linear map from the domain of $X_n^v$ to $V$, 
$b_n^v\in V$), where $u$ is an arbitrary point of the domain of $X_n^v$.  
Then we have 
$$\begin{array}{l}
\displaystyle{(X_n^v)_u=\frac1n{\rm Ad}(g_n)(X_{\rho(g_n^{-1})(u)})
=\frac1n{\rm Ad}(g_n)(A\rho(g_n^{-1})(u)+b)}\\
\displaystyle{=\frac1n({\rm Ad}(g_n)\circ A\circ{\rm Ad}(g_n^{-1}))(u)
+\frac1n{\rm Ad}(g_n)(A\rho(g_n^{-1})(\hat0)+b)}
\end{array}$$
and hence 
$$A_n^v=\frac1n{\rm Ad}(g_n)\circ A\circ{\rm Ad}(g_n^{-1})
\,\,\,\,{\rm and}
\,\,\,\,b_n^v=\frac1n{\rm Ad}(g_n)A(\rho(g_n^{-1})(\hat0)+b).\leqno{(7.10)}$$
From the first relation in $(7.10)$, we have 
$A_n^v\in{\mathfrak o}_{AK}(V)$ and hence $X^v_n\in{\mathcal K}^h$.  

\vspace{0.5truecm}



For $\{X_n^v\}_{n=1}^{\infty}$, we have the following fact.  

\vspace{0.5truecm}

\noindent
{\bf Lemma 7.3.6.} {\sl If $X\in{\mathcal K}^h_{\widetilde M^{\mathbb C}}$ and 
$v$ is an element of $H^{0,{\mathbb C}}_-$ with 
$$\exp\left(n\int_0^1v(t)\,dt\right)=e\quad(n\in{\mathbb N}),$$
then there exists a subsequence of $\{X_n^v\}_{n=1}^{\infty}$ converging to the zero 
vector field.}

\vspace{0.5truecm}

\noindent
{\it Proof.} Take $u\in V$.  Let $u=u_-+u_+$ 
($u_-\in H^{0,{\mathbb C}}_-,\,u_+\in H^{0,{\mathbb C}}_+$).  Then we have 
$$({\rm Ad}(g_n)u_{\varepsilon})(t)
={\rm Ad}(\exp(n\widetilde v(t)))u_{\varepsilon}(t)
=\exp({\rm ad}(n\widetilde v(t)))u_{\varepsilon}(t)\in
\mathfrak g_{\varepsilon}^{\mathbb C}\,\,\,(\varepsilon=-\,\,{\rm or}\,\,+)$$
for each $t\in [0,1]$ because $\widetilde v(t)\in\mathfrak g^{\mathbb C}_-\,\,
(0\leq t\leq 1)$ by the assumption and 
$[\mathfrak g^{\mathbb C}_-,\mathfrak g^{\mathbb C}_{\varepsilon}]\subset
\mathfrak g^{\mathbb C}_{\varepsilon}$ ($\varepsilon=-$ or $+$).  
Hence we have 
$$\begin{array}{l}
\displaystyle{\langle{\rm Ad}(g_n)u,{\rm Ad}(g_n)u
\rangle^{\mathcal I}=-\langle{\rm Ad}(g_n)u_-{\rm Ad}(g_n)u_-\rangle
+\langle{\rm Ad}(g_n)u_+,{\rm Ad}(g_n)u_+\rangle}\\
\hspace{4.2truecm}\displaystyle{=-\langle u_-,u_-\rangle
+\langle u_+,u_+\rangle=\langle u,u\rangle^{\mathcal I}.}
\end{array}$$
Therefore, by using $(7.10)$, we can show $|| A_n^v||_{\rm op}
=\frac1n|| A||_{\rm op}\to0$ ($n\to\infty$).  
Also, since $v\in\mathfrak g^{\mathbb C}_-$ and $G^{\mathbb C}_-$ is a compact Lie group, we have 
$$||\rho(g_n^{-1})(\hat 0)||=||-(g_n^{-1})'(g_n^{-1})_{\ast}^{-1}||=||(g_n^{-1})'||
=||\exp_{\ast}(nv)||\leq n||v||.$$
and hence 
$$|| b_n^v||\leq\frac1n\left(|| A\rho(g_n^{-1})(\hat0)||+|| b||\right)
\leq|| A||_{\rm op}\cdot||v||+\frac1n|| b||\to
||A||_{\rm op}\cdot||v||\,\,\,(n\to\infty).$$
Since the sequence $\{X_n^v\,|\,n\in{\bf N}\}$ in ${\mathcal K}^h$ is bounded, 
there exists its convergent subsequence $\{X^v_{n_j}\}_{j=1}^{\infty}$.  
Set $X^v_{\infty}:=\lim\limits_{j\to\infty}X^v_{n_j}$.  From 
$\lim\limits_{n\to\infty}A_n^v=0$, $X^v_{\infty}$ is a parallel Killing vector field 
on $V$.  
From $\displaystyle{\exp\left(n\int_0^1v(t)\,dt\right)=e}$, we have 
$g_n\in\Omega_e(G^{\mathbb C})$ and hence 
$\rho(g_n)(\widetilde M^{\mathbb C})=\widetilde M^{\mathbb C}$.  This fact together 
with $X\in{\mathcal K}^h_{\widetilde M^{\mathbb C}}$ deduces $X^v_n\in
{\mathcal K}^h_{\widetilde M^{\mathbb C}}$.  
Also, from $|| A_n^v||_{\rm op}=\frac1n|| A
||_{\rm op}<\infty$, we have $X^v_n\in{\mathcal K}^h$.  
Hence we have $X^v_n\in{\mathcal K}^h_{\widetilde M^{\mathbb C}}$.  
Therefore we have 
$X^v_{\infty}\in\overline{{\mathcal K}^h_{\widetilde M^{\mathbb C}}}$.  
Furthermore, from Lemma 7.3.5, 
we have $X^v_{\infty}\in{\mathcal K}^h_{\widetilde M^{\mathbb C}}$.  Thus, since 
$X^v_{\infty}$ is parallel and 
$X^v_{\infty}\in{\mathcal K}^h_{\widetilde M^{\mathbb C}}$, it follows 
from Lemma 7.3.4 that $X^v_{\infty}=0$.  This completes the proof.  
\qed

\vspace{0.3truecm}

On the other hand, we have the following fact.  

\vspace{0.3truecm}

\noindent
{\bf Lemma 7.3.7.} {\sl Let $X$ be an element of 
${\mathcal K}^h_{\widetilde M^{\mathbb C}}$ given by $X_u=Au+b\,\,(u\in V)$ for some 
$A\in\mathfrak o_{AK}(V)$ and some $b\in V$, 
$Y$ an element of $\mathfrak g^{\mathbb C}_-$ and $f$ an element of 
$H^0([0,1],{\mathbb C})(=H^0([0,1],{\mathbb R}^2))$ satisfying 
$\int_0^1f(t)dt=0$ or $f={\rm const}$.  
Then we have $A(fY)=[Y,w]$ for some $w\in V$.}

\vspace{0.3truecm}

\noindent
{\it Proof.} Set $v:=fY$.  
Define $\widetilde f\in H^1([0,1],{\mathbb C})$ by $\widetilde f(t):=
\int_0^tf(t)dt\,\,(0\leq t\leq 1)$.  
Let $A(fY)(t)=u_1(t)+u_2(t)$ ($u_1(t)\in{\rm Ker}\,
{\rm ad}(Y)$ and $u_2(t)\in{\rm Im}\,{\rm ad}(Y)$), and 
$u_i(t)=u_i^-(t)+u_i^+(t)$ 
($u_i^-(t)\in\mathfrak g_-^{\mathbb C},\,u_i^+(t)\in\mathfrak g_+^{\mathbb C}$) 
($i=1,2$) and 
$b(t)=b^-(t)+b^+(t)$ ($b^-(t)\in\mathfrak g_-^{\mathbb C},\,b^+(t)\in
\mathfrak g_+^{\mathbb C}$).  Let $g_n(t):=\exp(n\widetilde v(t))=\exp(n\widetilde f(t)Y)$.  
From $(7.10)$ 
and ${\rm Ad}(g_n)|_{{\rm Ker}\,{\rm ad}(Y)}={\rm id}$, we have 
$$\begin{array}{l}
\displaystyle{b_n^v=\frac1n{\rm Ad}(g_n)(A\rho(g_n^{-1})(\hat0)+b)
={\rm Ad}(g_n)\left(A(fY)+\frac bn\right)}\\
\hspace{0.45truecm}\displaystyle{=u_1+{\rm Ad}(g_n)(u_2+\frac bn).}
\end{array}$$
Since ${\rm Ad}(g_n)$ preserves 
$\mathfrak g_-^{\mathbb C}$ and $\mathfrak g_+^{\mathbb C}$ invariantly, respectively, 
and ${\rm Ad}(g_n)|_{{\rm Ker}\,{\rm ad}(Y)}={\rm id}$, we have 
$$\begin{array}{l}
\displaystyle{\langle b_n^v,u_1\rangle^{\mathcal I}
=\langle u_1,u_1\rangle^{\mathcal I}
+\langle{\rm Ad}(g_n)(u_2+\frac bn),
{\rm Ad}(g_n)u_1\rangle^{\mathcal I}}\\
\hspace{1.85truecm}\displaystyle{=\langle u_1,u_1
\rangle^{\mathcal I}+\frac1n\langle b,u_1\rangle^{\mathcal I}\to
\langle u_1,u_1\rangle^{\mathcal I}\,\,\,\,(n\to\infty).}
\end{array}$$
First we consider the case where ``$\int_0^1f(t)dt=0\,$'' or "$f={\rm const}$ and $Y$ is the initial 
vector of a closed geodesic in $G^{\mathbb C}_-$ of period $f\,\,$".  
Then we have $\exp\left(n\int_0^1v(t)dt\right)=e$ ($n\in{\mathbb N})$.  
Also we have $v\in H^{0,{\mathbb C}}_-$ because of $Y\in\mathfrak g^{\mathbb C}_-$.  
Hence, according to Lemma 7.3.6, there exists a subsequence 
$\{X^v_{n_i}\}_{i=1}^{\infty}$ of $\{X_n\}_{n=1}^{\infty}$ converging to the zero vector field.  
Clearly we have $\lim\limits_{i\to\infty}b^v_{n_i}=0$ and hence $u_1=0$.  
Thus we see that $A(fY)(t)\in{\rm Im}\,{\rm ad}(Y)$ holds for all 
$t\in[0,1]$.  That is, we have $A(fY)=[Y,w]$ for some $w\in V$.  
Next we consider the case where $f={\rm const}$ and $Y$ is the initial vector of a closed geodesic in 
$G^{\mathbb C}_-$ (not necessarily of period $f$).  
Let $a$ be the period of the closed geodesic.  Since $aY$ is the initial vector of a closed geodesic in 
$G^{\mathbb C}_-$ of period one, it follows from the above discussion that 
$A(aY)=[Y,\bar w]$ holds for some $\bar w\in V$.  
Hence we have 
$$A(fY)=\frac{f}{a}A(aY)=\frac{f}{a}[Y,\bar w]=\left[Y,\frac{f}{a}\bar w\right].$$
Next we consider the case where $f={\rm const}$ and $Y$ is the initial vector 
of non-closed geodesic in $G^{\mathbb C}_-$.   Set 
$$B:=\{Z\,|\,Z\,:\,{\rm the}\,\,{\rm initial}\,\,{\rm vector}\,\,{\rm of}\,\,{\rm a}\,\,{\rm closed}\,\,
{\rm geodesic}\,\,{\rm in}\,\,G^{\mathbb C}_-\}.$$
Since $\mathfrak g^{\mathbb C}_-$ is  the compact real of $\mathfrak g^{\mathbb C}$, 
$B$ is dense in $\mathfrak g^{\mathbb C}_-$.  Take a sequence 
$\{Z_i\}_{i=1}^{\infty}$ in $B$ with $\lim\limits_{i\to\infty}Z_i=fY$.  
As showed in the above, there exists $w_i\in V$ with $A(Z_i)=[Z_i,w_i]$ for each $i$.  
We can show that the sequence $\{w_i\}_{i=1}^{\infty}$ is a convergent sequence and that 
$$A(fY)=\lim_{i\to\infty}[Z_i,w_i]=[Y,f\lim_{i\to\infty}w_i].$$
This completes the proof.  \qed


\vspace{0.5truecm}

Since $w$ in this lemma depends on $X,f$ and $Y$, we denote it by 
$w_{X,f,Y}$.  
According to Lemma 2.10 of \cite{Ch}, we have the following fact.  

\vspace{0.5truecm}

\noindent
{\bf Lemma 7.3.8.} {\sl Let $B$ be a map from $\mathfrak g_-^{\mathbb C}$ to 
oneself defined by $B(Y)=[\mu(Y),Y]$ {\rm(}$Y\in\mathfrak g_-^{\mathbb C}${\rm)} 
in terms of a map $\mu:\mathfrak g_-^{\mathbb C}\to\mathfrak g_-^{\mathbb C}$.  
If $B$ is linear, then $\mu$ is a constant map.}

\vspace{0.5truecm}

By using Lemmas 7.3.7 and 7.3.8, we can show the following fact.  

\vspace{0.5truecm}

\noindent
{\bf Lemma 7.3.9.} {\sl Fix $X\in{\mathcal K}^h_{\widetilde M^{\mathbb C}}$ and 
$f\in H^0([0,1],{\mathbb C})$ satisfying $\int_0^1f(t)dt=0$ or $f={\rm const}$.  
Then $w_{X,f,Y}$ is independent of the choice of $Y\in\mathfrak g^{\mathbb C}_-$.}

\vspace{0.5truecm}

\noindent
{\it Proof.} For the simplicity, set $w_Y:=w_{X,f,Y}$.  Define a linear map 
$B_1^t:\mathfrak g_-^{\mathbb C}\to\mathfrak g_-^{\mathbb C}$ by $B_1^t(Y):=
A(fY)(t)_{\mathfrak g_-^{\mathbb C}}$ ($Y\in\mathfrak g_-^{\mathbb C}$) and 
a linear map $B_2^t:\mathfrak g_-^{\mathbb C}\to\mathfrak g_-^{\mathbb C}$ by 
$B_2^t(Y):=\sqrt{-1}(A(fY)(t)_{\mathfrak g_+^{\mathbb C}})$ 
($Y\in\mathfrak g_-^{\mathbb C}$), where 
$(\cdot)_{\mathfrak g_{\varepsilon}^{\mathbb C}}$ ($\varepsilon=-$ or $+$) is 
the $\mathfrak g_{\varepsilon}^{\mathbb C}$-component of $(\cdot)$.  
Since $A(fY)=[Y,w_Y]$, we have 
$B_1^t(Y)=[Y,w_Y(t)_{\mathfrak g_-^{\mathbb C}}]$ and 
$B_2^t(Y)=[Y,\sqrt{-1}w_Y(t)_{\mathfrak g_+^{\mathbb C}}]$, 
it follows from Lemma 7.3.8 that, for each $t\in[0,1]$, 
$w_Y(t)_{\mathfrak g^{\mathbb C}_-}$ and $w_Y(t)_{\mathfrak g^{\mathbb C}_+}$ 
are independent of the choice of $Y\in\mathfrak g_-^{\mathbb C}$.  
Hence $w_Y$ is independent of the choice of $Y\in\mathfrak g^{\mathbb C}_-$.  
\qed

\vspace{0.5truecm}

According to this lemma, $w_{X,f,Y}$ is independent of the choice of 
$Y\in\mathfrak g^{\mathbb C}_-$, we denote it by $w_{X,f}$.  
Define $\psi_n\in H^0([0,1],{\mathbb C})$ by $\psi_n(t)=\exp(2n\pi\sqrt{-1}t)$ 
($0\leq t\leq1$), where $n\in{\mathbb Z}$.  

\vspace{0.5truecm}

\noindent
{\bf Lemma 7.3.10.} {\sl For each $X\in{\mathcal K}^h_{\widetilde M^{\mathbb C}}$ 
and each $f\in H^0([0,1],{\mathbb C})$ satisfying 
$\int_0^1f(t)dt=0$ or $f={\rm const}$, 
we have $w_{X,f}=fw_{X,1}$, where the subscript $1$ in $w_{X,1}$ means 
$1\in H^0([0,1],{\mathbb C})$.}

\vspace{0.5truecm}

\noindent
{\it Proof.} Let $\langle\,\,,\,\,\rangle^{\mathbb C}$ be the complexification of 
the ${\rm Ad}(G)$-invariant non-degenerate symmetric bilinear form 
$\langle\,\,,\,\,\rangle$ of $\mathfrak g$ inducing the metric of $G/K$.  
Let $\mathfrak a$ be a maximal abelian subspace of $\sqrt{-1}\mathfrak p$ and 
$\mathfrak g^{\mathbb C}_-=\mathfrak z_{\mathfrak g^{\mathbb C}_-}(\mathfrak a)
+\sum_{\alpha\in\triangle}(\mathfrak g^{\mathbb C}_-)_{\alpha}$ the root space 
decomposition of $\mathfrak g^{\mathbb C}_-$ with respect to $\mathfrak a$, where 
$\mathfrak z_{\mathfrak g^{\mathbb C}_-}(\mathfrak a)$ is the centralizer of 
$\mathfrak a$ in $\mathfrak g^{\mathbb C}_-$ and 
$\triangle:=\{\alpha\in\mathfrak a^{\ast}\,|\,
(\mathfrak g^{\mathbb C}_-)_{\alpha}\not=\{0\}\}$ 
($(\mathfrak g^{\mathbb C}_-)_{\alpha}:=\{Z\in\mathfrak g^{\mathbb C}_-\,|\,
{\rm ad}(a)Z=\sqrt{-1}\alpha(a)Z\,\,(\forall\,a\in\mathfrak a)\}$).  
For any $\alpha\in\triangle$ and any $n\in{\bf N}\cup\{0\}$, define 
$H_{\alpha}\in\mathfrak a$ by $\langle H_{\alpha},\cdot
\rangle=\alpha(\cdot)$ and $c_{\alpha,n}:=\frac{2n\pi\sqrt{-1}}
{\alpha(H_{\alpha})}$.  
Define $g_{\alpha,n}\in H^1([0,1],G^{\mathbb C})$ by 
$g_{\alpha,n}(t):=\exp(tc_{\alpha,n}H_{\alpha})$ ($0\leq t\leq1$).  It is 
clear that $g_{\alpha,n}\in\Omega_e(G^{\mathbb C})$.  
Let $\overline X_{\alpha,n}:=\rho(g_{\alpha,n})_{\ast}^{-1}X$.  
Since $\rho(g_{\alpha,n})(\widetilde M^{\mathbb C})=\widetilde M^{\mathbb C}$, 
$\overline X_{\alpha,n}$ is tangent to $\widetilde M^{\mathbb C}$ along 
$\widetilde M^{\mathbb C}$.  Also, we can show 
$\overline X_{\alpha,n}\in{\mathcal K}^h$.  
Hence we have $\overline X_{\alpha,n}\in{\mathcal K}^h_{\widetilde M^{\mathbb C}}$.  
Let 
$(\overline X_{\alpha,n})_u=\overline A_{\alpha,n}u+\overline b_{\alpha,n}$ 
($\overline A_{\alpha,n}\in\mathfrak o_{AK}(V),
\,\overline b_{\alpha,n}\in V$).  
We can show that $\overline A_{\alpha,n}={\rm Ad}(g_{\alpha,n})^{-1}\circ A
\circ{\rm Ad}(g_{\alpha,n})$ in similar to the first relation in (7.10).  
Take any $Y_0\in\mathfrak z_{\mathfrak g^{\mathbb C}_-}(\mathfrak a)$ and any 
$Y_{\alpha}\in(\mathfrak g^{\mathbb C}_-)_{\alpha}$.  
Then, from ${\rm Ad}(g_{\alpha,n})Y_0=Y_0$, we have 
$$\begin{array}{l}
\displaystyle{[{\rm Ad}(g_{\alpha,n})w_{\overline X_{\alpha,n},1},Y_0]
=[{\rm Ad}(g_{\alpha,n})w_{\overline X_{\alpha,n},1},{\rm Ad}(g_{\alpha,n})Y_0]
}\\
\displaystyle{=-{\rm Ad}(g_{\alpha,n})(\overline A_{\alpha,n}Y_0)
=-A({\rm Ad}(g_{\alpha,n})Y_0)=-AY_0=[w_{X,1},Y_0].}
\end{array}$$
It follows from the arbitrariness of 
$Y_0(\in\mathfrak z_{\mathfrak g^{\mathbb C}_-}(\mathfrak a)$) that 
$${\rm Im}({\rm Ad}(g_{\alpha,n})
w_{\overline X_{\alpha,n},1}-w_{X,1})\subset\mathfrak a.\leqno{(7.11)}$$
Also, from ${\rm Ad}(g_{\alpha,n})Y_{\alpha}=\psi_nY_{\alpha}$, we have 
$$\begin{array}{l}
\hspace{0.5truecm}\displaystyle{[{\rm Ad}(g_{\alpha,n})
w_{\overline X_{\alpha,n},1},Y_{\alpha}]
=\psi_{-n}[{\rm Ad}(g_{\alpha,n})w_{\overline X_{\alpha,n},1},
{\rm Ad}(g_{\alpha,n})Y_{\alpha}]}\\
\displaystyle{=-\psi_{-n}{\rm Ad}(g_{\alpha,n})
(\overline A_{\alpha,n}Y_{\alpha})=-\psi_{-n}A({\rm Ad}(g_{\alpha,n})
Y_{\alpha})}\\
\displaystyle{=-\psi_{-n}A(\psi_nY_{\alpha})
=\psi_{-n}[w_{X,\psi_n},Y_{\alpha}]}
\end{array}$$
and hence 
$$[{\rm Ad}(g_{\alpha,n})w_{\overline X_{\alpha,n},1}
-\psi_{-n}w_{X,\psi_n},Y_{\alpha}]=0.$$
It follows from the arbitrariness of 
$Y_{\alpha}(\in(\mathfrak g^{\mathbb C}_-)_{\alpha})$ that 
$${\rm Im}\left({\rm Ad}(g_{\alpha,n})w_{\overline X_{\alpha,n},1}
-\psi_{-n}w_{X,\psi_n}\right)\subset\mathfrak z_{\mathfrak g^{\mathbb C}_-}
((\mathfrak g^{\mathbb C}_-)_{\alpha}).$$
This together with $(7.11)$ implies 
$${\rm Im}\left(\psi_nw_{X,1}-w_{X,\psi_n}\right)
\subset\mathfrak a\oplus\mathfrak z_{\mathfrak g^{\mathbb C}_-}
((\mathfrak g^{\mathbb C}_-)_{\alpha}).$$
From the arbitrariness of $\alpha$, we obtain 
$${\rm Im}\left(\psi_nw_{X,1}-w_{X,\psi_n}\right)
\subset\mathfrak a\oplus\left(\mathop{\cap}_{\alpha\in\triangle}
\mathfrak z_{\mathfrak g^{\mathbb C}_-}((\mathfrak g^{\mathbb C}_-)_{\alpha})\right)
=\mathfrak a.$$
Take another maximal abelian subspace $\mathfrak a'$ of 
$\sqrt{-1}\mathfrak p$ with $\mathfrak a'\cap\mathfrak a=\{0\}$.  Similarly 
we can show 
$${\rm Im}\left(\psi_nw_{X,1}-w_{X,\psi_n}\right)
\subset\mathfrak a'$$
and hence 
$$w_{X,\psi_n}=\psi_nw_{X,1}.\leqno{(7.12)}$$
Take any $f\in H^0([0,1],{\mathbb C})$ satisfying $\int_0^1f(t)dt=0$ or 
$f={\rm const}$.  Let $f=\sum\limits_{n=-\infty}^{\infty}c_n\psi_n$ be the 
Fourier's expansion of $f$, where $c_n$ is constant for each $n$.  Then, 
since $A$ is continuous and linear, we have 
$$A(fY)=\sum\limits_{n=-\infty}^{\infty}
c_nA(\psi_nY)\,\,\,\,(Y\in\mathfrak g^{\mathbb C}_-).\leqno{(7.13)}$$
From $(7.12)$ and $(7.13)$, we obtain 
$$[Y,w_{X,f}]=A(fY)=\sum_{n=-\infty}^{\infty}c_n
[Y,w_{X,\psi_n}]=[Y,fw_{X,1}]\,\,\,\,(Y\in\mathfrak g^{\mathbb C}_-).$$
Thus $w_{X,f}-fw_{X,1}$ belongs to the center of $\mathfrak g^{\mathbb C}_-$.  
Therefore, since $\mathfrak g^{\mathbb C}_-$ has no center, we obtain 
$w_{X,f}=fw_{X,1}$.  \qed

\vspace{0.7truecm}

From Lemmas 7.3.7 and 7.3.10, we have the following fact.  

\vspace{0.7truecm}

\noindent
{\bf Lemma 7.3.11.} {\sl Let $X$ be an element of 
${\mathcal K}^h_{\widetilde M^{\mathbb C}}$ given by 
$X_u=Au+b$ ($u\in V$) for some $A\in\mathfrak o_{AK}(V)$ and $b\in V$.  
Then we have $A={\rm ad}(v)$ for some $v\in V$.}

\vspace{0.7truecm}

\noindent
{\it Proof.} Take any $u\in V$ and a base $\{e_1,\cdots,e_m\}$ of 
$\mathfrak g^{\mathbb C}_-$.  Let $u=\sum\limits_{i=1}^mu_ie_i$ and 
$u_i=\sum\limits_{n=-\infty}^{\infty}c_{i,n}\psi_n$ be the Fourier expansion of $u_i$.  
Then, since $A$ is continuous and linear, we have 
$Au=\sum\limits_{n=-\infty}^{\infty}\sum\limits_{i=1}^mc_{i,n}A(\psi_ne_i)$.  
According to Lemmas 7.3.7 and 7.3.10, we have 
$A(fY)=[w_{X,1},fY]$ for any $Y\in\mathfrak g^{\mathbb C}_-$ and any 
$f\in H^0([0,1],{\mathbb C})$ satisfying $\int_0^1f(t)dt=0$ or $f={\rm const}$.  
Hence we have 
$$Au=\sum\limits_{n=-\infty}^{\infty}\sum_{i=1}^mc_{i,n}[w_{X,1},\psi_ne_i]=[w_{X,1},u].$$
Thus we obtain $A={\rm ad}(w_{X,1})$.  
\qed

\vspace{0.7truecm}

By using Lemmas 7.3.4 and 7.3.11, we shall prove Proposition 7.3.  

\vspace{0.7truecm}

\noindent
{\it Proof of Proposition 7.3.} 
Take any $X\in{\rm Lie}\,H_b$.  Since ${\rm Lie}\,H_b\subset{\mathcal K}^h_{\widetilde M^{\mathbb C}}$, 
it follows from Lemmas 7.3.4 and 7.3.11 that $X=X^v$ for some $v\in V$.  
Since $X^v$ is the holomorphic Killing vector field associated with an one-parameter 
subgroup $\{\rho(\exp\circ sv)\,|\,s\in{\mathbb R}\}$ of 
$\rho(H^1([0,1],G^{\mathbb C}))$, we have $X\in{\rm Lie}\,\rho(H^1([0,1],
G^{\mathbb C}))$.  Hence we obtain ${\rm Lie}\,H_b\subset{\rm Lie}\,
\rho(H^1([0,1],G^{\mathbb C})$, that is, $H_b\subset\rho(H^1([0,1],G^{\mathbb C}))$.  
\qed

\vspace{0.7truecm}

By using Proposition 7.3, we shall prove Theorem B.  

\vspace{0.7truecm}

\noindent
{\it Proof of Theorem B.} 
Since $H_b$ is a subgroup of $\rho(H^1([0,1],G^{\mathbb C}))$ by 
Proposition 7.3, we have $H_b=\rho(Q)$ for some subgroup $Q$ of $H^1([0,1],G^{\mathbb C})$.  
Let $Q'$ be a closed connected subgroup of $G^{\mathbb C}\times G^{\mathbb C}$ 
generated by $\{(h(0),h(1))\,|\,h\in Q\}$.  Since $\phi\circ\rho(h)=(L_{h(0)}\circ R_{h(1)}^{-1})\circ\phi$ 
for each $h\in H$, we have $\widehat M^{\mathbb C}=Q'\cdot e$, where $e$ is the identity element of 
$G^{\mathbb C}$.  Here we note that $G^{\mathbb C}\times G^{\mathbb C}$ acts on $G^{\mathbb C}$ by 
$(g_1,g_2)\cdot g:=(L_{g_1}\circ R_{g_2}^{-1})(g)\,\,\,(g_1,g_2,g\in G^{\mathbb C})$.  
Set $\widehat M:=\pi_{\mathbb R}^{-1}(M)$, where $\pi_{\mathbb R}$ is the natural projection of $G$ onto $G/K$.   
Since $\widehat M$ is a component of $\widehat M^{\mathbb C}\cap G$ containing $e$ and $(Q'\cap(G\times G))\cdot e$ 
is a complete open submanifold of $\widehat M^{\mathbb C}\cap G$, $\widehat M$ is a component of 
$(Q'\cap(G\times G))\cdot e$.  Therefore we have $\widehat M=(Q'\cap(G\times G))_0\cdot e$, where 
$(Q'\cap(G\times G))_0$ is the identity component of $Q'\cap(G\times G)$.  
Set $Q'_{\mathbb R}:=(Q'\cap(G\times G))_0$.  Since $\widehat M$ consists of fibres of $\pi_{\mathbb R}$, we have 
$\langle Q'_{\mathbb R}\cup(e\times K)\rangle\cdot e=\widehat M$, where 
$\langle Q'_{\mathbb R}\cup(e\times K)\rangle$ is the group generated by 
$Q'_{\mathbb R}\cup(e\times K)$.  Denote by the same symbol $Q'_{\mathbb R}$ the group 
$\langle Q'_{\mathbb R}\cup(e\times K)\rangle$ under abuse of the notation.  
Set $(Q'_{\mathbb R})_1:=\{g_1\in G\,|\,\exists\,g_2\in G\,\,{\rm s.t.}\,\,
(g_1,g_2)\in Q'_{\mathbb R}\}$ and $(Q'_{\mathbb R})_2:=\{g_2\in G\,|\,\exists\,
g_1\in G\,\,{\rm s.t.}\,\,(g_1,g_2)\in Q'_{\mathbb R}\}$.  Also, set 
$(Q'_{\mathbb R})^{\bullet}_1:=\{g\in G\,|\,(g,e)\in Q'_{\mathbb R}\}$ and 
$(Q'_{\mathbb R})^{\bullet}_2:=\{g\in G\,|\,(e,g)\in Q'_{\mathbb R}\}$.  
It is clear that $(Q'_{\mathbb R})^{\bullet}_i$ is a normal subgroup of 
$(Q'_{\mathbb R})_i$ ($i=1,2$).  From $e\times K\subset Q'_{\mathbb R}$, we have 
$K\subset(Q'_{\mathbb R})^{\bullet}_2$.  Since $K\subset(Q'_{\mathbb R})^{\bullet}_2
\subset(Q'_{\mathbb R})_2\subset G$ and $K$ is a maximal subgroup of $G$, we have 
$(Q'_{\mathbb R})_2=K$ or $G$ and $(Q'_{\mathbb R})^{\bullet}_2=K$ or $G$.  
Suppose that $(Q'_{\mathbb R})^{\bullet}_2=G$.  Then we have $\widehat M=G$ 
and hence $M=G/K$.  Thus a contradiction arises.  
Hence we have $(Q'_{\mathbb R})^{\bullet}_2=K$.  Since $K$ is not a normal 
subgroup of $G$ and it is a normal subgroup of $(Q'_{\mathbb R})_2$, 
we have $(Q'_{\mathbb R})_2\not=G$.  Therefore we have 
$(Q'_{\mathbb R})_2=K$ and hence $Q'_{\mathbb R}\subset G\times K$.  
Set $Q''_{\mathbb R}:=\{g\in G\,|\,(\{g\}\times K)\cap Q'_{\mathbb R}\not=\emptyset
\}$.  Then, since $\widehat M=Q'_{\mathbb R}\cdot e$ and $M=\pi(\widehat M)$, 
we have $M=Q''_{\mathbb R}(eK)$.  Thus $M$ is extrinsically homogeneous.  \qed

\section{Proof of Theorem C} 
In this section, we prove Theorem C (Main theorem) by using Theorems A and B.  
Let $M$ be as in Theorem C and $F$ be its reflective focal submanifold.  
Without loss of genereality, we may assume that $o:=eK\in F$.  Denote by $A$ the shape tensor of $M$ and $R$ 
the curvature tensor of $G/K$.  

First we prove the following fact by using Theorem A.  

\vspace{0.5truecm}

\noindent
{\bf Proposition 8.1.} {\sl The submanifold $M$ satisfies the condition $(\ast_{\mathbb C})$.}


\vspace{0.5truecm}

\noindent
{\it Proof.} We prove this statement in the case where $G/K$ is of non-compact type (this statement is proved 
similarly in the case where $G/K$ is of compact type).  Take $Z_0\in\mathfrak p$ with ${\rm Exp}\,Z_0\in M$.  
Set $x_0:={\rm Exp}\,Z_0,\,\mathfrak t:=T_oF,\,\mathfrak t^{\perp}:=T^{\perp}_oF$ and 
$\mathfrak b:=(\exp\,Z_0)_{\ast o}^{-1}(T^{\perp}_{x_0}M)$.  
We use the notations in the proof of Theorem A (in Section 6).  
Take any $v\in T^{\perp}_{x_0}M$.  As stated in the proof of Theorem A, the decomposition 
$$T_{{\rm Exp}\,Z_0}M=\left(\mathop{\oplus}_{\beta\in(\triangle_{\mathfrak b})^H_+\cup\{0\}}
(\mathfrak p_{\beta}\cap\mathfrak t)^L_{Z_0}\right)\oplus
\left(\mathop{\oplus}_{\beta\in(\triangle_{\mathfrak b})^V_+}
(\exp\,Z_0)_{\ast o}(\mathfrak p_{\beta}\cap\mathfrak t^{\perp})\right)$$
is the common eigenspace decomposition of $A_v$ and $R(v)$.  
Also, we have $R(v)|_{(\exp\,Z_0)_{\ast o}(\mathfrak p_{\beta})}=\beta(v)^2\,{\rm id}$ 
($\beta\in(\triangle_{\mathfrak b})_+\cup\{0\}$).  
From (ii) of Proposition 5.3 that 
$$(\mathfrak p_{\beta}\cap\mathfrak t)^L_{Z_0}\subset{\rm Ker}(A_v+\beta(\bar v)\tanh(\beta(Z_0)){\rm id})
\quad(\beta\in(\triangle_{\mathfrak b})_+^H).$$
Also, since $F$ is reflective and the fibre $M\cap{\rm Exp}(\mathfrak t^{\perp})$ is a principal orbit of the 
isotropy action of the symmetric space ${\rm Exp}(\mathfrak t^{\perp})$, 
it follows from (i) of Proposition 5.3 that 
$$(\exp\,Z_0)_{\ast}(\mathfrak p_{\beta}\cap\mathfrak t^{\perp})
\subset{\rm Ker}\left(A_v+\frac{\beta(\bar v)}{\tanh(\beta(Z_0))}{\rm id}\right).$$
From these facts, it follows that are not equal the absolute values of the eigenvalues $A_v$ and $R(v)$ 
on each of the common eigenspaces $(\mathfrak p_{\beta}\cap\mathfrak t)^L_{Z_0}$'s 
($\beta\in(\triangle_{\mathfrak b})^H_+\cup\{0\}$) and 
$(\exp\,Z_0)_{\ast o}(\mathfrak p_{\beta}\cap\mathfrak t^{\perp})$'s ($(\triangle_{\mathfrak b})^V_+$) of $A_v$ 
and $R(v)$.  This implies that $M$ satisfies the condition $(\ast_{\mathbb C})$.  \qed

\vspace{0.5truecm}

From Theorem B and this proposition, we can derive Theorem C.  

\vspace{0.5truecm}

\noindent
{\it Proof of Theorem C.} 
Since $M$ satisfies the condition $(\ast_{\mathbb C})$ by Proposition 8.1, it follows from Theorem B that 
$M$ is extrinsically homogeneous.  Hence it follows from Theorem A of \cite{Koi6} that $M$ is a principal orbit of 
a (complex) hyperpolar action on $G/K$.  See \cite{Koi3} (or \cite{Koi6}) about the definition of a (complex) 
hyperpolar action.  Furthermore, since this action admits a reflective (hence totally geodesic) singular orbit 
and it is of cohomogeneity greater than one, it follows from Theorem C and Remark 1.1 of \cite{Koi6} that this action 
is orbit equivalent to a Hermann type action.  Therefore $M$ is a principal orbit of a Hermann type action.  \qed

\section{Classifications}
From Theorem C and the list of Hermann type actions in \cite{Koi6}, 
we can classify isoparametric submanifolds as in Theorem C as follows.  

\vspace{0.7truecm}

\noindent
{\bf Theorem 9.1.} {\sl Let $M$ be a full irreducible isoparametric $C^{\omega}$-submanifold of codimension greater 
than one in an irreducible symmetric space $G/K$ of non-compact type.  If $M$ admits a reflective focal submanifold, 
then it is a principal orbit of the action of one of symmetric subgroups $H$'s of $G$ as in Tables $1$-$3$.}

\vspace{1truecm}

$$\begin{tabular}{|c|c|}
\hline
$G/K$ & $H$\\
\hline
\scriptsize{$SL(n,{\mathbb R})/SO(n)$} & 
\scriptsize{$SO(n),\,\,\,\,SO_0(p,n-p)\,\,(1\leq p\leq n-1),\,\,\,\,
Sp(\frac n2,{\mathbb R}),\,\,\,\,SL(\frac n2,{\mathbb C})\cdot U(1)$}\\
\scriptsize{($n\geq 6,\,\,n:$ even)} & 
\scriptsize{$(SL(p,{\mathbb R})\times SL(n-p,{\mathbb R}))\cdot{\mathbb R}_{\ast}
\,\,\,\,(2\leq p\leq n-2)$}\\
\hline
\scriptsize{$SL(4,{\mathbb R})/SO(4)$} & 
\scriptsize{$SO(4),\,\,\,\,SO_0(1,3),\,\,\,\,SO_0(2,2),\,\,\,\,
SL(2,{\mathbb C})\cdot U(1),\,\,\,\,
(SL(2,{\mathbb R})\times SL(2,{\mathbb R}))\cdot{\mathbb R}_{\ast}$}\\
\hline
\scriptsize{$SL(n,{\mathbb R})/SO(n)$} & 
\scriptsize{$SO(n),\,\,\,\,SO_0(p,n-p)\,\,(1\leq p\leq n-1),$}\\
\scriptsize{($n\geq 5,\,\,n:$ odd)} & 
\scriptsize{$(SL(p,{\mathbb R})\times SL(n-p,{\mathbb R}))\cdot{\mathbb R}_{\ast}\,\,
(2\leq p\leq n-2)$}\\
\hline
\scriptsize{$SL(3,{\mathbb R})/SO(3)$} & 
\scriptsize{$SO(3),\,\,\,\,SO_0(1,2)$}\\
\hline
\scriptsize{$SU^{\ast}(2n)/Sp(n)\,\,(n\geq 4)$} & 
\scriptsize{$Sp(n),\,\,\,\,SO^{\ast}(2n),\,\,\,\,Sp(p,n-p)\,\,
(1\leq p\leq n-1),\,\,\,\,SL(n,{\mathbb C})\cdot U(1)$}\\
\scriptsize{} & \scriptsize{$SU^{\ast}(2p)\times SU^{\ast}(2n-2p)\times U(1)
\,\,(2\leq p\leq n-2)$}\\
\hline
\scriptsize{$SU^{\ast}(6)/Sp(3)$} & 
\scriptsize{$Sp(3),\,\,\,\,SO^{\ast}(6),\,\,\,\,Sp(1,2)$}\\
\hline
\scriptsize{$SU(p,q)/S(U(p)\times U(q))$} & 
\scriptsize{$S(U(p)\times U(q)),\,\,\,\,SO_0(p,q),\,\,\,\,\,Sp(\frac p2,
\frac q2),$}\\
\scriptsize{$(4\leq p<q,\,p,q:{\rm even})$} & 
\scriptsize{$S(U(i,j)\times U(p-i,q-j))\,\,\,\,(1\leq i\leq p-1,\,
1\leq j\leq q-1)$}\\
\hline
\scriptsize{$SU(p,q)/S(U(p)\times U(q))$} & 
\scriptsize{$S(U(p)\times U(q)),\,\,\,\,SO_0(p,q),$}\\
\scriptsize{$(3\leq p<q,\,\,\,p\,\,{\rm or}\,\,q:{\rm odd})$} & 
\scriptsize{$S(U(i,j)\times U(p-i,q-j))\,\,(1\leq i\leq p-1,\,
1\leq j\leq q-1)$}\\
\hline
\scriptsize{$SU(2,q)/S(U(2)\times U(q))$} & 
\scriptsize{$S(U(2)\times U(q)),\,\,\,\,SO_0(2,q),\,\,\,\,\,S(U(1,j)\times 
U(1,q-j))\,\,(1\leq j\leq q-1)$}\\
\scriptsize{$(q\geq3)$} & \scriptsize{}\\
\hline
\scriptsize{$SU(p,p)/S(U(p)\times U(p))$} & 
\scriptsize{$S(U(p)\times U(p)),\,\,\,\,SO_0(p,p),\,\,\,\,SO^{\ast}(2p),
\,\,\,\,Sp(\frac p2,\frac p2),\,\,\,\,Sp(p,{\mathbb R}),\,\,\,\,
SL(p,{\mathbb C})\cdot U(1)$}\\
\scriptsize{$(p\geq4,\,\,\,p:{\rm even})$} & 
\scriptsize{$S(U(i,j)\times U(p-i,p-j))\,\,\,\,(1\leq i\leq p-1,\,
1\leq j\leq p-1)$}\\
\hline
\scriptsize{$SU(2,2)/S(U(2)\times U(2))$} & 
\scriptsize{$S(U(2)\times U(2)),\,\,\,\,SO_0(2,2),\,\,\,\,SO^{\ast}(4),
\,\,\,\,SL(2,{\mathbb C})\cdot U(1),\,\,\,\,S(U(1,1)\times U(1,1))$}\\
\hline
\scriptsize{$SU(p,p)/S(U(p)\times U(p))$} & 
\scriptsize{$S(U(p)\times U(p)),\,\,\,\,SO_0(p,p),\,\,\,\,SO^{\ast}(2p),
\,\,\,\,Sp(p,{\mathbb R}),\,\,\,\,SL(p,{\mathbb C})\cdot U(1)$}\\
\scriptsize{$(p\geq5,\,\,\,p:{\rm odd})$} & 
\scriptsize{$S(U(i,j)\times U(p-i,p-j))\,\,\,\,(1\leq i\leq p-1,\,1\leq j\leq 
p-1)$}\\
\hline
\scriptsize{$SU(3,3)/S(U(3)\times U(3))$} & 
\scriptsize{$S(U(3)\times U(3)),\,\,\,\,SO_0(3,3),\,\,\,\,SO^{\ast}(6),\,\,\,\,
SL(3,{\mathbb C})\cdot U(1),$}\\
\scriptsize{} & \scriptsize{$S(U(1,1)\times U(2,2)),\,\,\,\,
S(U(1,2)\times U(2,1))$}\\
\hline
\scriptsize{$SL(n,{\mathbb C})/SU(n)$} & 
\scriptsize{$SU(n),\,\,\,\,SO(n,{\mathbb C}),\,\,\,\,SL(n,{\mathbb R}),\,\,\,\,
SU(i,n-i)\,\,\,(1\leq i\leq n-1),\,\,\,\,
Sp(\frac n2,{\mathbb C}),\,\,\,\,SU^{\ast}(n)$}\\
\scriptsize{$(n\geq 6,\,\,\,n:{\rm even})$} & 
\scriptsize{$SL(i,{\mathbb C})\times SL(n-i,{\mathbb C})\times U(1)\,\,\,\,
(2\leq i\leq n-2)$}\\
\hline
\scriptsize{$SL(4,{\mathbb C})/SU(4)$} & 
\scriptsize{$SU(4),\,\,\,\,SO(4,{\mathbb C}),\,\,\,\,SL(4,{\mathbb R}),\,\,\,\,
SU(i,4-i)\,\,\,(1\leq i\leq 3),\,\,\,\,SU^{\ast}(4)$}\\
\scriptsize{} & 
\scriptsize{$SL(2,{\mathbb C})\times SL(2,{\mathbb C})\times U(1)$}\\
\hline
\scriptsize{$SL(n,{\mathbb C})/SU(n)$} & 
\scriptsize{$SU(n),\,\,\,\,SO(n,{\mathbb C}),\,\,\,\,SL(n,{\mathbb R}),\,\,\,\,
SU(i,n-i)\,\,\,(1\leq i\leq n-1)$}\\
\scriptsize{$(n\geq 5,\,\,\,n:{\rm odd})$} & 
\scriptsize{$SL(i,{\mathbb C})\times SL(n-i,{\mathbb C})\times U(1)\,\,\,\,
(2\leq i\leq n-2)$}\\
\hline
\scriptsize{$SL(3,{\mathbb C})/SU(3)$} & 
\scriptsize{$SU(3),\,\,\,\,SO(3,{\mathbb C})$}\\
\hline
\end{tabular}$$

\vspace{1truecm}

\centerline{{\bf Table 1.}}

\vspace{0.5truecm}

$$\begin{tabular}{|c|c|}
\hline
$G/K$ & $H$\\
\hline
\scriptsize{$SO_0(p,q)/SO(p)\times SO(q)$} & 
\scriptsize{$SO(p)\times SO(q),\,\,\,\,SU(\frac p2,\frac q2)\cdot U(1),$}\\
\scriptsize{$(4\leq p<q,\,\,\,p,q:{\rm even})$} & 
\scriptsize{$SO_0(i,j)\times SO_0(p-i,q-j)\,\,\,\,(1\leq i\leq p-1,\,
1\leq j\leq q-1)$}\\
\hline
\scriptsize{$SO_0(2,q)/SO(2)\times SO(q)$} & 
\scriptsize{$SO(2)\times SO(q),\,\,\,\,SO_0(1,j)\times SO_0(1,q-j)\,\,\,\,
(1\leq j\leq q-1)$}\\
\scriptsize{$(4\leq q,\,\,\,q:{\rm even})$} & 
\scriptsize{}\\
\hline
\scriptsize{$SO_0(p,q)/SO(p)\times SO(q)$} & 
\scriptsize{$SO(p)\times SO(q),\,\,\,\,SO_0(i,j)\times SO_0(p-i,q-j)\,\,\,\,
(1\leq i\leq p-1,\,\,1\leq j\leq q-1)$}\\
\scriptsize{$(2\leq p<q,\,\,\,p\,\,{\rm or}\,\,q:{\rm odd})$} & 
\scriptsize{}\\
\hline
\scriptsize{$SO_0(p,p)/SO(p)\times SO(p)$} & 
\scriptsize{$SO(p)\times SO(p),\,\,\,\,SO(p,{\mathbb C}),\,\,\,\,
SU(\frac p2,\frac p2)\cdot U(1),\,\,\,\,SL(p,{\mathbb R})\cdot U(1)$}\\
\scriptsize{$(p\geq4,\,\,\,p:{\rm even})$} & 
\scriptsize{$SO_0(i,j)\times SO_0(p-i,p-j)\,\,\,\,(1\leq i\leq p-1,\,
1\leq j\leq p-1)$}\\
\hline
\scriptsize{$SO_0(2,2)/SO(2)\times SO(2)$} & 
\scriptsize{$SO(2)\times SO(2),\,\,\,\,SO(2,{\mathbb C}),\,\,\,\,
SO_0(1,1)\times SO_0(1,1)$}\\
\hline
\scriptsize{$SO_0(p,p)/SO(p)\times SO(p)$} & 
\scriptsize{$SO(p)\times SO(p),\,\,\,\,SO(p,{\mathbb C}),\,\,\,\,
SL(p,{\mathbb R})\cdot U(1),$}\\
\scriptsize{$(p\geq5,\,\,\,p:{\rm odd})$} & 
\scriptsize{$SO_0(i,j)\times SO_0(p-i,p-j)\,\,
(1\leq i\leq p-1,\,1\leq j\leq p-1)$}\\
\hline
\scriptsize{$SO_0(3,3)/SO(3)\times SO(3)$} & 
\scriptsize{$SO(3)\times SO(3),\,\,\,\,SO(3,{\mathbb C}),\,\,\,\,
SO_0(1,1)\times SO_0(2,2)$}\\
\scriptsize{} & \scriptsize{$SO_0(1,2)\times SO_0(2,1)$}\\
\hline
\scriptsize{$SO^{\ast}(2n)/U(n)$} & 
\scriptsize{$U(n),\,\,\,\,SO(n,{\mathbb C}),\,\,\,\,SU^{\ast}(n)\cdot U(1)$}\\
\scriptsize{$(n\geq6,\,\,\,n:{\rm even})$} & 
\scriptsize{$SO^{\ast}(2i)\times SO^{\ast}(2n-2i)\,\,\,(2\leq i\leq n-2),$}\\
\scriptsize{} & \scriptsize{$SU(i,n-i)\cdot U(1)\,\,\,\,
(\left[\frac i2\right]+\left[\frac{n-i}{2}\right]\geq2)$}\\
\hline
\scriptsize{$SO^{\ast}(8)/U(4)$} & 
\scriptsize{$U(4),\,\,\,\,SO(4,{\mathbb C}),\,\,\,\,
SO^{\ast}(4)\times SO^{\ast}(4),\,\,\,\,SU(2,2)\cdot U(1)$}\\
\hline
\scriptsize{$SO^{\ast}(2n)/U(n)$} & 
\scriptsize{$U(n),\,\,\,\,SO(n,{\mathbb C}),\,\,\,\,
SO^{\ast}(2i)\times SO^{\ast}(2n-2i)\,\,\,\,(2\leq i\leq n-2),$}\\
\scriptsize{$(n\geq5,\,\,\,n:{\rm odd})$} & 
\scriptsize{$SU(i,n-i)\cdot U(1)\,\,\,\,
(\left[\frac i2\right]+\left[\frac{n-i}{2}\right]\geq2)$}\\
\hline
\scriptsize{$SO(n,{\mathbb C})/SO(n)$} & 
\scriptsize{$SO(n),\,\,\,\,SO(i,{\mathbb C})\times SO(n-i,{\mathbb C})\,\,\,\,
(2\leq i\leq n-2),$}\\
\scriptsize{$(n\geq 8,\,\,\,n:{\rm even})$} & 
\scriptsize{$SO_0(i,n-i)\,\,\,\,(\left[\frac i2\right]+\left[\frac{n-i}{2}
\right]\geq2),\,\,\,\,\,SL(\frac n2,{\mathbb C})\cdot SO(2,{\mathbb C}),\,\,\,\,
SO^{\ast}(n)$}\\
\hline
\scriptsize{$SO(6,{\mathbb C})/SO(6)$} & 
\scriptsize{$SO(6),\,\,\,\,SO(i,{\mathbb C})\times SO(6-i,{\mathbb C})\,\,\,\,
(2\leq i\leq 4),$}\\
\scriptsize{} & 
\scriptsize{$SO_0(2,4),\,\,\,\,SO_0(3,3),\,\,\,\,SO^{\ast}(6)$}\\
\hline
\scriptsize{$SO(4,{\mathbb C})/SO(4)$} & 
\scriptsize{$SO(4),\,\,\,\,SO(2,{\mathbb C})\times SO(2,{\mathbb C}),\,\,\,\,
SO_0(2,2),\,\,\,\,SO^{\ast}(4)$}\\
\hline
\scriptsize{$SO(n,{\mathbb C})/SO(n)$} & 
\scriptsize{$SO(n),\,\,\,\,SO(i,{\mathbb C})\times SO(n-i,{\mathbb C})\,\,\,
(2\leq i\leq n-2),$}\\
\scriptsize{$(n\geq 5,\,\,\,n:{\rm odd})$} & 
\scriptsize{$SO_0(i,n-i)\,\,\,\,
(\left[\frac i2\right]+\left[\frac{n-i}{2}\right]\geq2)$}\\
\hline
\scriptsize{$Sp(n,{\mathbb R})/U(n)$} & 
\scriptsize{$U(n),\,\,\,\,SU(i,n-i)\cdot U(1)\,\,\,\,(1\leq i\leq n-1),\,\,\,\,
SL(n,{\mathbb R})\cdot U(1),$}\\
\scriptsize{$(n\geq 4,\,\,n:{\rm even})$} & 
\scriptsize{$Sp(\frac n2,{\mathbb C}),\,\,\,\,
Sp(i,{\mathbb R})\times Sp(n-i,{\mathbb R})\,\,\,\,(2\leq i\leq n-2)$}\\
\hline
\scriptsize{$Sp(2,{\mathbb R})/U(2)$} & 
\scriptsize{$U(2),\,\,\,\,SU(1,1)\cdot U(1)$}\\
\hline
\scriptsize{$Sp(n,{\mathbb R})/U(n)$} & 
\scriptsize{$U(n),\,\,\,\,SU(i,n-i)\cdot U(1)\,\,\,\,(1\leq i\leq n-1),\,\,\,\,
SL(n,{\mathbb R})\cdot U(1),$}\\
\scriptsize{$(n\geq 5,\,\,n:{\rm odd})$} & 
\scriptsize{$Sp(i,{\mathbb R})\times Sp(n-i,{\mathbb R})\,\,\,\,
(2\leq i\leq n-2)$}\\
\hline
\scriptsize{$Sp(3,{\mathbb R})/U(3)$} & 
\scriptsize{$U(3),\,\,\,\,SU(1,2)\cdot U(1),\,\,\,\,
SL(3,{\mathbb R})\cdot U(1)$}\\
\hline
\end{tabular}$$

\vspace{0.5truecm}

\centerline{{\bf Table 2.}}

\vspace{0.5truecm}

$$\begin{tabular}{|c|c|}
\hline
$G/K$ & $H$\\
\hline
\scriptsize{$Sp(p,q)/Sp(p)\times Sp(q)$} & 
\scriptsize{$Sp(p)\times Sp(q),\,\,\,\,SU(p,q)\cdot U(1),$}\\
\scriptsize{$(2\leq p<q)$} & 
\scriptsize{$Sp(i,j)\times Sp(p-i,q-j)\,\,\,\,(1\leq i\leq p-1,\,1\leq j\leq 
q-1)$}\\
\hline
\scriptsize{$Sp(p,p)/Sp(p)\times Sp(p)$} & 
\scriptsize{$Sp(p)\times Sp(p),\,\,\,\,SU(p,p)\cdot U(1),\,\,\,\,
SU^{\ast}(2p)\cdot U(1),\,\,\,\,Sp(p,{\mathbb C})$}\\
\scriptsize{$(p\geq3)$} & 
\scriptsize{$Sp(i,j)\times Sp(p-i,p-j)\,\,\,\,
(1\leq i\leq p-1,\,1\leq j\leq p-1)$}\\
\hline
\scriptsize{$Sp(2,2)/Sp(2)\times Sp(2)$} & 
\scriptsize{$Sp(2)\times Sp(2),\,\,\,\,SU(2,2)\cdot U(1),\,\,\,\,
SU^{\ast}(4)\cdot U(1),\,\,\,\,Sp(1,1)\times Sp(1,1)$}\\
\hline
\scriptsize{$Sp(n,{\mathbb C})/Sp(n)$} & 
\scriptsize{$Sp(n),\,\,\,\,SL(n,{\mathbb C})\cdot SO(2,{\mathbb C}),\,\,\,\,
Sp(n,{\mathbb R}),\,\,\,\,Sp(i,n-i)\,\,\,(1\leq i\leq n-1),$}\\
\scriptsize{$(n\geq 4)$} & 
\scriptsize{$Sp(i,{\mathbb C})\times Sp(n-i,{\mathbb C})\,\,\,(2\leq i\leq n-2)$}\\
\hline
\scriptsize{$Sp(n,{\mathbb C})/Sp(n)$} & 
\scriptsize{$Sp(n),\,\,\,\,SL(n,{\mathbb C})\cdot SO(2,{\mathbb C}),\,\,\,\,
Sp(n,{\mathbb R}),\,\,\,\,Sp(i,n-i)\,\,\,(1\leq i\leq n-1)$}\\
\scriptsize{$(n=2,3)$} & \scriptsize{}\\
\hline
\scriptsize{$E_6^6/(Sp(4)/\{\pm1\})$} & 
\scriptsize{$Sp(4)/\{\pm1\},\,\,\,\,Sp(4,{\mathbb R}),\,\,\,\,Sp(2,2),\,\,\,\,
SU^{\ast}(6)\cdot SU(2),$}\\
\scriptsize{} & \scriptsize{$SL(6,{\mathbb R})\times SL(2,{\mathbb R}),\,\,\,\,
SO_0(5,5)\cdot{\mathbb R},\,\,\,\,F_4^4$}\\
\hline
\scriptsize{$E_6^2/SU(6)\cdot SU(2)$} & 
\scriptsize{$SU(6)\cdot SU(2),\,\,\,\,Sp(1,3),\,\,\,\,Sp(4,{\mathbb R}),\,\,\,\,
SU(2,4)\cdot SU(2),\,\,\,\,SU(3,3)\cdot SL(2,{\mathbb R}),$}\\
\scriptsize{} & \scriptsize{$SO^{\ast}(10)\cdot U(1),\,\,\,\,
SO_0(4,6)\cdot U(1)$}\\
\hline
\scriptsize{$E_6^{-14}/Spin(10)\cdot U(1)$} & 
\scriptsize{$Spin(10)\cdot U(1),\,\,\,\,Sp(2,2),\,\,\,\,SU(2,4)\cdot SU(2),
\,\,\,\,SU(1,5)\cdot SL(2,{\mathbb R}),$}\\
\scriptsize{} & \scriptsize{$SO^{\ast}(10)\cdot U(1),\,\,\,\,
SO_0(2,8)\cdot U(1)$}\\
\hline
\scriptsize{$E_6^{-26}/F_4$} & 
\scriptsize{$F_4,\,\,\,\,F_4^{-20},\,\,\,\,Sp(1,3)$}\\
\hline
\scriptsize{$E_6^{\mathbb C}/E_6$} & 
\scriptsize{$E_6,\,\,\,\,E_6^6,\,\,\,\,E_6^2,\,\,\,\,E_6^{-14},\,\,\,\,
Sp(4,{\mathbb C}),\,\,\,\,SL(6,{\mathbb C})\cdot SL(2,{\mathbb C}),\,\,\,\,
SO(10,{\mathbb C})\cdot Sp(1),\,\,\,\,F_4^{\mathbb C}.\,\,\,\,E_6^{-26}$}\\
\hline
\scriptsize{$E_7^7/(SU(8)/\{\pm1\})$} & 
\scriptsize{$SU(8)/\{\pm1\},\,\,\,\,SL(8,{\mathbb R}),\,\,\,\,SU^{\ast}(8),
\,\,\,\,SU(4,4),\,\,\,\,SO^{\ast}(12)\cdot SU(2),$}\\
\scriptsize{} & \scriptsize{$SO_0(6,6)\cdot SL(2,{\mathbb R}),\,\,\,\,
E_6^6\cdot U(1),\,\,\,\,E_6^2\cdot U(1)$}\\
\hline
\scriptsize{$E_7^{-5}/SO'(12)\cdot SU(2)$} & 
\scriptsize{$SO'(12)\cdot SU(2),\,\,\,\,SU(4,4),\,\,\,\,SU(2,6),\,\,\,\,
SO^{\ast}(12)\cdot SL(2,{\mathbb R}),$}\\
\scriptsize{} & \scriptsize{$SO_0(4,8)\cdot SU(2),\,\,\,\,E_6^2\cdot U(1),
\,\,\,\,E_6^{-14}\cdot U(1)$}\\
\hline
\scriptsize{$E_7^{-25}/E_6\cdot U(1)$} & 
\scriptsize{$E_6\cdot U(1),\,\,\,\,SU^{\ast}(8),\,\,\,\,SU(2,6),\,\,\,\,
SO^{\ast}(12)\cdot SU(2),$}\\
\scriptsize{} & \scriptsize{$SO_0(2,10)\cdot SL(2,{\mathbb R}),\,\,\,\,
E_6^{-14}\cdot U(1),\,\,\,\,E_6^{-26}\cdot U(1)$}\\
\hline
\scriptsize{$E_7^{\mathbb C}/E_7$} & 
\scriptsize{$E_7,\,\,\,\,E_7^7,\,\,\,\,E_7^{-5},\,\,\,\,E_7^{-25},\,\,\,\,
SL(8,{\mathbb C}),\,\,\,\,SO(12,{\mathbb C})\cdot SL(2,{\mathbb C}),\,\,\,\,
E_6^{\mathbb C}\cdot{\mathbb C}^{\ast}$}\\
\hline
\scriptsize{$E_8^8/SO'(16)$} & 
\scriptsize{$SO'(16),\,\,\,\,SO^{\ast}(16),\,\,\,\,SO_0(8,8),\,\,\,\,
E_7^{-5}\cdot Sp(1),\,\,\,\,E_7^7\cdot SL(2,{\mathbb R})$}\\
\hline
\scriptsize{$E_8^{-24}/E_7\cdot Sp(1)$} & 
\scriptsize{$E_7\cdot Sp(1),\,\,\,\,E_7^{-5}\cdot Sp(1),\,\,\,\,
E_7^{-25}\cdot SL(2,{\mathbb R}),\,\,\,\,SO^{\ast}(16),\,\,\,\,SO_0(4,12)$}\\
\hline
\scriptsize{$E_8^{\mathbb C}/E_8$} & 
\scriptsize{$E_8,\,\,\,\,E_8^8,\,\,\,\,E_8^{-24},\,\,\,\,SO(16,{\mathbb C}),
\,\,\,\,E_7^{\mathbb C}\times SL(2,{\mathbb C})$}\\
\hline
\scriptsize{$F_4^4/Sp(3)\cdot Sp(1)$} & 
\scriptsize{$Sp(3)\cdot Sp(1),\,\,\,\,Sp(1,2)\cdot Sp(1),\,\,\,\,
Sp(3,{\mathbb R})\cdot SL(2,{\mathbb R})$}\\
\hline
\scriptsize{$F_4^{\mathbb C}/F_4$} & 
\scriptsize{$F_4,\,\,\,\,F_4^4,\,\,\,\,F_4^{-20},\,\,\,\,
Sp(3,{\mathbb C})\cdot SL(2,{\mathbb C})$}\\
\hline
\scriptsize{$G_2^2/SO(4)$} & 
\scriptsize{$SO(4),\,\,\,\,SL(2,{\mathbb R})\times SL(2,{\mathbb R}),\,\,\,\,
\alpha(SO(4))$}\\
\scriptsize{} & \scriptsize{($\alpha:$an outer automorphism of $G_2^2$)}\\
\hline
\scriptsize{$G_2^{\mathbb C}/G_2$} & 
\scriptsize{$G_2,\,\,\,\,G_2^2,\,\,\,\,SL(2,{\mathbb C})\times SL(2,{\mathbb C})$}\\
\hline
\end{tabular}$$

\vspace{0.5truecm}

\centerline{{\bf Table 3.}}

\section{Proof of Theorem D} 
In 1991, G. Thorbergsson (\cite{Th}) proved that any full irreducible isoparametric submanifold of codimension 
greater than two in a Euclidean space is extrinsically homogeneous by using the building theory.  
In this section, we shall prove Theorem D by defining the topological Tits building of spherical type 
associated to an isoparametric submanifold as in Theorem D and using it, where we refer the proof in \cite{Th}.  
First we recall the notion of a topological Tits building.  
Let $\Delta=({\mathcal V},{\mathcal S})$ be an $r$-dimensional simplicial complex, where 
${\mathcal V}$ denotes the set of all vertices and ${\mathcal S}$ denotes the set of all simplices.  
Each $r$-simplex of $\Delta$ is called a {\it chamber} of $\Delta$.  
Let ${\mathcal A}:=\{{\mathcal A}_{\lambda}\}_{\lambda\in\Lambda}$ be a family of subcomplexes of $\Delta$.  
The pair ${\mathcal B}:=(\Delta,{\mathcal A})$ is called a {\it Tits building} if the following conditions hold:

\vspace{0.25truecm}

\noindent
(B1) Each $(r-1)$-dimensional simplex of $\Delta$ is contained in at least three chambers.  

\noindent
(B2) Each $(r-1)$-dimensional simplex in a subcomplex ${\mathcal A}_{\lambda}$ are contained in exactly two chambers 
of ${\mathcal A}_{\lambda}$.  

\noindent
(B3) Any two simplices of $\Delta$ are contained in some ${\mathcal A}_{\lambda}$.  

\noindent
(B4) If two subcomplexes ${\mathcal A}_{\lambda_1}$ and ${\mathcal A}_{\lambda_2}$ share a chamber, then there is 
an isomorphism of ${\mathcal A}_{\lambda_1}$ onto ${\mathcal A}_{\lambda_2}$ fixing 
${\mathcal A}_{\lambda_1}\cap{\mathcal A}_{\lambda_2}$ pointwisely.  

\vspace{0.25truecm}

\noindent 
Each subcomplex belonging to ${\mathcal A}$ is called an {\it apartment} of ${\mathcal B}$.  
In this appendix, we assume that all Tits building futhermore satisfies the following condition: 

\vspace{0.25truecm}

\noindent
(B5) Each apartment ${\mathcal A}_{\lambda}$ is a Coxeter complex.  

\vspace{0.25truecm}

\noindent
If ${\mathcal A}_{\lambda}$ is finite (resp. infinite), then the building ${\mathcal B}$ is said to be {\it spherical type} 
(resp. {\it affine type}).  Let ${\mathcal O}$ be a Hausdorff topology of ${\mathcal V}$.  The pair 
$({\mathcal B},{\mathcal O})$ is called a {\it topological Tits building} if the following conditions hold:

\vspace{0.25truecm}

\noindent
(TB1) $({\mathcal B},{\mathcal A})$ is a Tits building.  

\noindent
(TB2) For $k\in\{1,\cdots,r\}$, $\widehat{\mathcal S}_k:=\{(x_1,\cdots,x_{k+1})\in{\mathcal V}^{k+1}\,\vert\,\,
\vert x_1\cdots x_{k+1}\vert\in{\mathcal S}_k\}$ is closed in the product topological space 
$({\mathcal V}^{k+1},{\mathcal O}^{k+1})$, where ${\mathcal S}_k$ denotes the set of all $k$-simplices of 
${\mathcal S}$ and $\vert x_1\cdots x_{k+1}\vert$ denotes the $k$-simplex with vertices $x_1,\cdots,x_{k+1}$.  

\vspace{0.25truecm}

A homeomorphism $\phi$ of $({\mathcal V},{\mathcal O})$ is called a {\it topological automorphism} of 
the topological Tits building $(\Delta,{\mathcal A},{\mathcal O})$ if the following conditions hold:

\vspace{0.25truecm}

\noindent
(TA1) $\phi$ preserves ${\mathcal S}$ (i.e., 
``$\sigma=\vert x_1\cdots x_{k+1}\vert\in{\mathcal S}\Rightarrow\phi(\sigma)
:=\vert\phi(x_1)\cdots\phi(x_{k+1})\vert\in{\mathcal S}$.  

\noindent
(TA2) $\phi$ preserves ${\mathcal A}$ (i.e., for each $\lambda\in\Lambda$, 
$\phi({\mathcal A}_{\lambda}):=\{\phi(\sigma)\,\vert\,\sigma\in{\mathcal A}_{\lambda}\}\in{\mathcal A}$.)

\noindent
(TA3) For each $k\in\{1,\cdots,r\}$, $\phi$ gives a homeomorphism of $\widehat{\mathcal S}_k$ onto oneself.  

\vspace{0.25truecm}

\noindent
According to (TA1) (resp. (TA2)), $\phi$ gives a bijection of ${\mathcal S}$ onto oneself 
(resp. ${\mathcal A}$ onto oneself).  

Let $M$ be a full irreducible curvature-adapted isoparametric submanifold of codimension $r(\geq2)$ in 
an irreducible symmetric space $G/K$ of non-compact type.  Assume that $M$ satisfies the condition 
$(\ast'_{\mathbb R})$.  
Set $\mathfrak p:=T_{eK}(G/K)$ and $\mathfrak b:=T^{\perp}_{eK}M$.  
Let $\mathfrak a$ be a maximal abelian subspace of $\mathfrak p\,(\subset\mathfrak g)$ containing $\mathfrak b$ 
and $\displaystyle{\mathfrak p=\mathfrak a\oplus
\left(\mathop{\oplus}_{\alpha\in\triangle_+}\mathfrak p_{\alpha}\right)}$ be the root space 
decomposition with respect to $\mathfrak a$, that is, $\mathfrak p_{\alpha}:=\{X\in\mathfrak p\,|\,
{\rm ad}(a)^2(X)=\alpha(a)^2X\,\,(\forall\,\,a\in\mathfrak a)\}$ and 
$\triangle_+$ is the positive root system of the root system 
$\triangle:=\{\alpha\in{\mathfrak a}^{\ast}\setminus\{0\}\,|\,
\mathfrak p_{\alpha}\not=\{0\}\}$ under a lexicographic ordering of 
${\mathfrak a}^{\ast}$.  
Set $\triangle_{\mathfrak b}:=\{\alpha|_{\mathfrak b}\,|\,\alpha\in
\triangle\,\,{\rm s.t.}\,\,\alpha|_{\mathfrak b}\not=0\}$ and 
let $\displaystyle{\mathfrak p=\mathfrak z_{\mathfrak p}(\mathfrak b)\oplus
\left(\mathop{\oplus}_{\beta\in(\triangle_{\mathfrak b})_+}\mathfrak p_{\beta}\right)}$ be the root space 
decomposition with respect to $\mathfrak b$, where $\mathfrak z_{\mathfrak p}(\mathfrak b)$ is 
the centralizer of $\mathfrak b$ in $\mathfrak p$, 
$\displaystyle{\mathfrak p_{\beta}=\mathop{\oplus}_{\alpha\in\triangle_+\,{\rm s.t.}\,\alpha|_{\mathfrak b}=\pm\beta}
\mathfrak p_{\alpha}}$ and $(\triangle_{\mathfrak b})_+$ is the positive root system of the root system 
$\triangle_{\mathfrak b}$ under a lexicographic ordering of 
$\mathfrak b^{\ast}$.  
For convenience, we denote $\mathfrak z_{\mathfrak p}(\mathfrak b)$ by $\mathfrak p_0$.  
Denote by $A$ the shape tensor of $M$ and $R$ the curvature tensor of $G/K$.  
Let $m_A:=\displaystyle{\mathop{\max}_{v\in\mathfrak b\setminus\{0\}}
\sharp{\rm Spec}\,A_v}$ and 
$m_R:=\displaystyle{\mathop{\max}_{v\in\mathfrak b\setminus\{0\}}\sharp{\rm Spec}\,R(v)}$, where $\sharp(\cdot)$ 
is the cardinal number of $(\cdot)$.  Note that $m_R=\sharp(\triangle_{\mathfrak b})_+$.  
Let $U:=\{v\in\mathfrak b\setminus\{0\}\,|\,\sharp{\rm Spec}A_v=m_A,\,\,\sharp{\rm Spec}\,R(v)=m_R\}$, 
which is an open dense subset of $\mathfrak b\setminus\{0\}$.  Fix $v\in U$.  Note that 
${\rm Spec}\,R(v)=\{-\beta(v)^2\,|\,\beta\in(\triangle_{\mathfrak b})_+\}$.  From $v\in U$, $\beta(v)^2$'s 
($\beta\in(\triangle_{\mathfrak b})_+$) are mutually distinct.  
Let ${\rm Spec}A_v=\{\lambda^v_1,\cdots,\lambda^v_{m_A}\}$ ($\lambda^v_1>\cdots>\lambda^v_{m_A}$).  
Set 
$$\begin{array}{l}
\displaystyle{I_0^v:=\{i\,|\,\mathfrak p_0\cap{\rm Ker}
(A_v-\lambda_i^v{\rm id})\not=\{0\}\},}\\
\displaystyle{I_{\beta}^v:=\{i\,|\,\mathfrak p_{\beta}\cap{\rm Ker}
(A_v-\lambda^v_i{\rm id})\not=\{0\}\},}\\
\displaystyle{(I_{\beta}^v)^+:=\{i\in I_{\beta}^v\,|\,
|\lambda^v_i|\,>\,|\beta(v)|\},}\\
\displaystyle{(I_{\beta}^v)^-:=\{i\in I_{\beta}^v\,|\,
|\lambda^v_i|\,<\,|\beta(v)|\},}\\
\displaystyle{(I_{\beta}^v)^0:=\{i\in I_{\beta}^v\,|\,
|\lambda^v_i|=|\beta(v)|\}.}
\end{array}$$
Since $M$ is curvature-adapted and satisfies the condition $(\ast'_{\mathbb R})$, 
we have $I^0_v=\emptyset$, $(I_{\beta}^v)^-=\emptyset$ (i.e., $I_{\beta}^v=(I_{\beta}^v)^+$) 
($\beta\in(\triangle_{\mathfrak b})_+$) and $\mathfrak a=\mathfrak b$ (hence $\triangle=\triangle_{\mathfrak b}$).  
In similar to the fact $(2.2)$ stated in Section 2, we have 
$${\mathcal FR}^{\mathbb R}_{M,v}=\left\{\left.\frac{1}{\beta(v)}{\rm arctanh}\frac{\beta(v)}{\lambda_i^v}\,\right|\,
\beta\in\triangle_+,\,\,i\in I_{\beta}^v\right\}.\leqno{(10.1)}$$
From the arbitrarinees of $v$ and the fact that $U$ is open and dense in $\mathfrak b$, the relation $(10.1)$ holds 
for any $v\in \mathfrak b$.  
Hence the tangential focal set ${\mathcal F}^{\mathbb R}_{M,eK}$ of $M$ at $eK$ is given by 
$${\mathcal F}^{\mathbb R}_{M,eK}=\mathop{\bigcup}_{v\in T^{\perp}_xM\,\,{\rm s.t.}\,\,||v||=1}
\left\{\left.\frac{1}{\beta(v)}{\rm arctanh}\frac{\beta(v)}{\lambda_i^v}\cdot v\,\right|\,
\beta\in\triangle_+,\,\,i\in I_{\beta}^v\right\}.\leqno{(10.2)}$$
On the other hand, H. Ewert (\cite{E}) showed that the tangential focal set of an isoparametric submanifold 
in a symmetric spaces of non-compact type at any point consists of finitely many (real) hyperplanes 
(which are called {\it focal hyperplanes}) in the normal space at the point and the reflections with respect to 
the hyperplanes generates a Weyl group (see \cite{E} for example), 
where we note that he (\cite{E}) treated not only an isoparametric submanifold(=equifocal submanifold) but also 
a submanifold with parallel focal structure (whose sections are not necessarily flat).  
Denote by ${\mathcal W}$ this Weyl group.  Note that the focal hyperplanes are not parallel pairwisely because 
the Weyl gorup is a finite Coxeter group.  
From this fact and $(10.2)$, we see that, for any $\beta\in\triangle_+$, $\sharp\,I_{\beta}^v=1$ and 
$\frac{\beta(v)}{\lambda_i^v}$ is independent of the choice of $v$ and furthermore 
$\frac{\beta(v)}{\lambda_i^v}=\frac{2\beta(v)}{\lambda_j^v}$ holds when $\beta,2\beta\in\triangle_+$, where 
$\{i\}=I_{\beta}$ and $\{j\}=I_{2\beta}$.  So we set $c_{\beta}:=\frac{\beta(v)}{\lambda_i^v}$ 
($\beta\in\triangle'_+$) and furthemore $\hat c_{\beta}:={\rm arctanh}\,c_{\beta}$.  Also, 
set $\triangle'_+:=\{\beta\in\triangle_+\,|\,2\beta\notin\triangle_+\}$ and $k:=\sharp\,\triangle'_+$.  
Then ${\mathcal F}^{\mathbb R}_{M,eK}$ is given by 
$${\mathcal F}^{\mathbb R}_{M,eK}=\mathop{\bigcup}_{\beta\in\triangle'_+}\beta^{-1}(\hat c_{\beta}).\leqno{(7.3)}$$
This fact implies that ${\mathcal W}$ is isomorphic to the Weyl group of $G/K$ (that is, the Coxeter group of 
the principal orbits (which are isoparametric submanifolds) of the s-representation of $G/K$).  
Since $M$ is full and irreducible, we can show that ${\mathcal W}$ is of rank $r$ and irreducible.  
Therefore $G/K$ is irreducible and its rank is equal to $r$.  
For the simplicity, set ${\it l}_{\beta}:=\beta^{-1}(\hat c_{\beta})$.  
It is clear that $\displaystyle{\mathop{\cap}_{\beta\in\triangle_+}{\it l}_{\beta}}$ is a one-point set.  
Denote by $v_0$ this point and set $p_0:=\exp^{\perp}(v_0)$ and $r_0:=||v_0||$.  
It is clear that the section $\Sigma_x$ of $M$ through any $x\in M$ passes through $p_0$.  
Let $S(r_0)$ be the unit sphere of radius $r_0$ centered at $0$ in $T_{p_0}(G/K)$.  
It is easy to show that $M$ is included by the geodesic sphere $\exp_{p_0}(S(r_0))$ in $G/K$.  
Let $\{{\it l}_i^x\,|\,i=1,\cdots,k\}$ be the set of all focal hyperplanes of $M$ at $x(\in M)$, that is, 
$\displaystyle{\mathop{\cup}_{i=1}^k{\it l}_i^x={\mathcal F}^{\mathbb R}_{M,x}}$.  
Set $\bar{\it l}_i^x:=\exp^{\perp}({\it l}_i^x),\,\widetilde{\it l}_i^x:=\exp_{p_0}^{-1}(\bar{\it l}_i^x)$ 
and $\widetilde{\Sigma}_x:=\exp_{p_0}^{-1}(\Sigma_x)$, 
where we note that $\widetilde{\Sigma}_x$ is an $r$-dimensional affine subspace in $T_{p_0}(G/K)$ through $0$ 
beacuse $\Sigma_x$ is a flat totally geodesic submanifold in $G/K$, and that $\widetilde{\it l}_i^x$ is an (affine) 
hyperplane in $\widetilde{\Sigma}_x$ through $0$.  It is clear that $\widetilde{\it l}_i^x\cap S(r_0)$'s ($i\in I_x$) 
and their intersections give a Coxeter complex in $\widetilde{\Sigma}_x\cap S(r_0)$.  
Denote by ${\mathcal A}_x$ this Coxeter complex.  
Let ${\mathcal V}_x$ (resp. ${\mathcal S}_x$) be the set of all vertices (resp. simplices) of ${\mathcal A}_x$.  
Set ${\mathcal V}_M:=\mathop{\cup}_{x\in M}{\mathcal V}_x$, ${\mathcal S}_M:=\mathop{\cup}_{x\in M}{\mathcal S}_x$ 
and ${\mathcal A}_M:=\{{\mathcal A}_x\,\vert\,x\in M\}$.  Also, set $\Delta_M:=({\mathcal V}_M,{\mathcal S}_M)$.  
Give ${\mathcal V}_M$ the relative topology (which we denote by ${\mathcal O}$) of $T_{p_0}(G/K)$.  
Note that $\exp_{p_0}({\mathcal V}_M)$ is equal to the sum of some lower dimensional submanifold 
$F_1,\cdots,F_{\it l}$.  It is shown that $F_1,\cdots,F_{\it l}$ are focal submanifolds of $M$.  
For example, see Figure 6 about the case where ${\mathcal A}_x$ is a Coxeter complex of type $(A_2)$.  
We have the following fact:

\vspace{0.225truecm}

{\sl $(\sharp)\quad$ 
$\displaystyle{{\mathcal F}^{\mathbb R}_{M,eK}=\mathop{\cup}_{\beta\in\triangle'_+}{\it l}_{\beta}
=\mathop{\cup}_{\beta\in\triangle'_+}\beta^{-1}(\hat c_{\beta})}$, the nullity space corresponding to the focal 
hyperplane 

\hspace{0.85truecm}${\it l}_{\beta}$ is equal to $\mathfrak p_{\beta}$ and $A_v=\lambda_i^v\,{\rm id}
=\frac{\beta(v)}{c_{\beta}}\,{\rm id}$ on $\mathfrak p_{\beta}$.}

\vspace{0.225truecm}

\noindent
Set $M':=\exp_{p_0}^{-1}(M)(\subset T_{p_0}(G/K))$.  
It is clear that $M'$ is included by $S(r_0)$.  
Also, we can show that $M'$ meets $\widetilde{\Sigma}_x$'s ($x\in M$) orthogonally 
by calculating the Jacobi vector fields along each radial geodesic starting from $p_0$ and reaching $M$, 
where we use the fact that the sections $\Sigma_x$'s ($x\in M$) are flat.  
Assume that $r\geq3$.  
Let $L$ be a principal orbit of the s-representation of $G/K$, which is a full irreducible isoparametric submanifold 
of codimension $r$ in $T_{eK}(G/K)$.  It is clear that the same fact as $(\sharp)$ holds at any point of $M$ 
(other than $eK$).  Hence it is shown that the above ${\mathcal B}_M:=(\Delta_M,{\mathcal A}_M,{\mathcal O})$ 
essentially coincides with the topological Tits building of spherical type associated to the full irreducible 
isoparametric submanifold $L$ constructed in \cite{Th} by comparing their constructions.  
Thus ${\mathcal B}_M$ is a topological Tits building of spherical type.  

\vspace{10pt}

\begin{figure}[h]
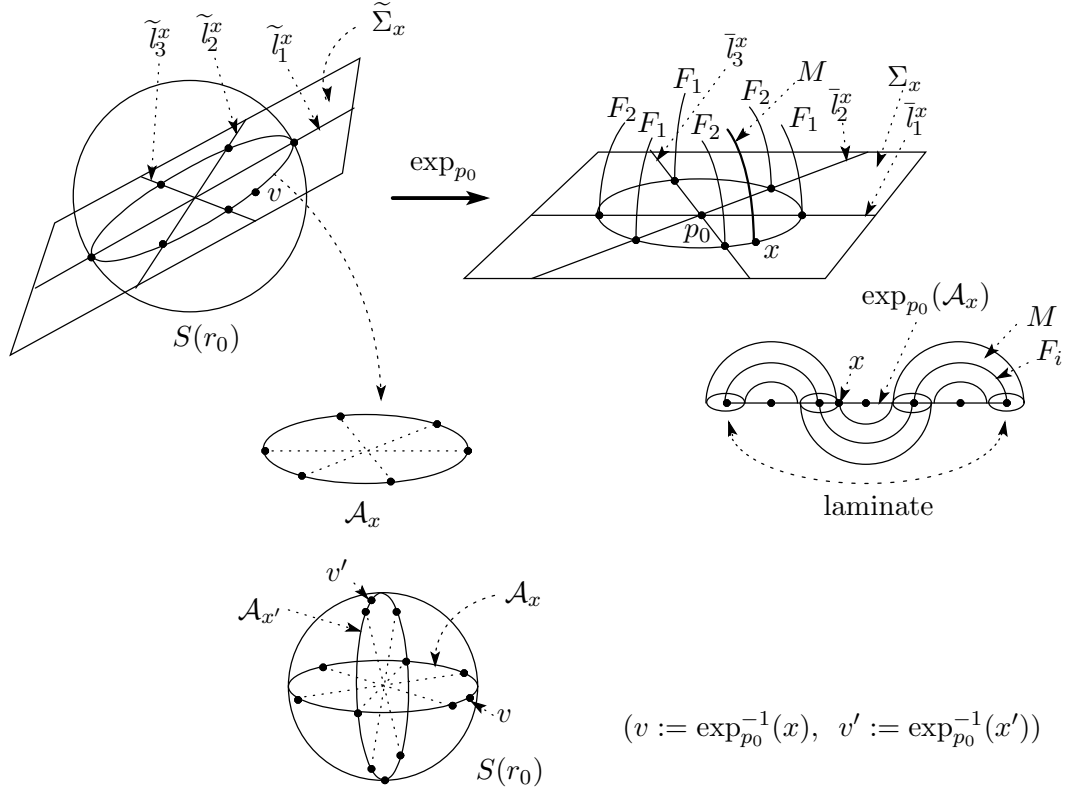

\centerline{
\unitlength 0.1in
%
\hspace{1truecm}}
\caption{The topological Tits building of an isoparametric submanifold as in Theorem D}
\label{figure6}
\end{figure}

\vspace{10pt}

Now we prove Theorem D by using this topological Tits building ${\mathcal B}_M$ of spherical type.  

\vspace{0.4truecm}

\noindent
{\it Proof of Theorem D.} Let $G'$ be the topological automorphism group of ${\mathcal B}_M$ and $G'_0$ be 
its identity component.  Then, by the result in \cite{BS}, it is shown that $G'_0$ is a semi-simple Lie group.  
Define an involution $s$ of ${\mathcal S}_M$ by $s(\sigma):=\{-p\,\vert\,p\in\sigma\}$ ($\sigma\in{\mathcal S}_M$).  
Let $K'$ be the subgroup consisting of all elements of $G'_0$ commuting with $s$.  
It is shown that $K'$ is a maximal compact subgroup of $G'_0$.  
We identify $T_{eK'}(G'_0/K')$ with $T_{p_0}(G/K)$ and denote these by the same symbol $\mathfrak p'$.  
We consider the action of $K'$ on $\mathfrak p'$ constructed as in the second paragraph of Section 4 (Page 444) 
of \cite{Th}.  That is, we consider the action of $K'$ on $\mathfrak p'$ constructed as follows.  
Take $k'\in K'$ and $v\in\mathfrak p'$.  Let $\sigma$ be the element of ${\mathcal S}_M$ including 
$\frac{r_0}{||v||}v$.  Let $w(k',v)$ be the element of $k'(\sigma)$ having the same barycentric coordinate as 
the barycentric coordinate of $\frac{r_0}{||v||}v$ with respect to $\sigma$.  
We define the $K'$-action on $\mathfrak p'$ by 
$$k'\cdot v:=\frac{||v||}{r_0}w(k',v)\quad\,\,(k'\in K',\,\,v\in\mathfrak p')$$
(see Figure 7).  From this construction, it is clear that this action $K'\curvearrowright\mathfrak p'$ has 
$M'$ as its orbit.  
It is shown that this action is a polar action on $\mathfrak p'$ by using the discussion in Page 444-445 of 
\cite{Th}, where we use also the fact that $M'$ meets $\widetilde{\Sigma}_x$'s ($x\in M$) orthogonally.  
Hence it follows that this action is orbit equivalent to the $s$-representation of $G'_0/K'$.  
Furthermore, since the same fact as $(\sharp)$ holds at any $x\in M$ other than $eK$, 
this action $K'\curvearrowright\mathfrak p'$ is orbit equivalent to the $s$-representation of $G/K$.  
Therefore $M'$ is a principal orbit of the $s$-representation of $G/K$ and hence 
$M$ is a principal orbit of the isotropy action $K\curvearrowright G/K$.  \qed

\vspace{10pt}

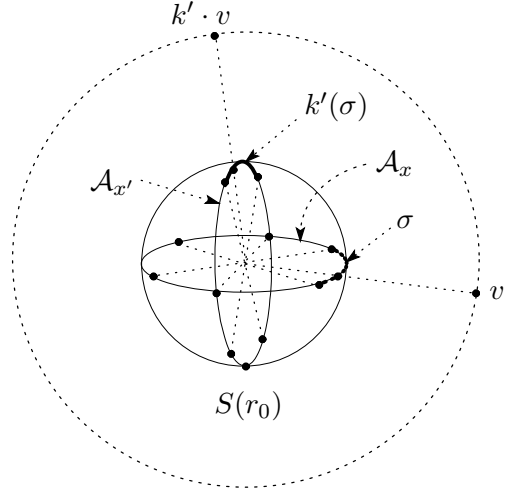
\begin{figure}[h]
\centerline{
\unitlength 0.1in
\begin{picture}( 35.4700, 25.9000)(  3.5300,-26.3000)
%
\special{pn 4}%
\special{ar 2628 1464 530 534  0.0000000 6.2831853}%
%
\special{pn 4}%
\special{ar 2628 1464 532 148  0.0000000 6.2831853}%
%
\special{pn 4}%
\special{pa 2616 930}%
\special{pa 2646 940}%
\special{pa 2668 962}%
\special{pa 2686 990}%
\special{pa 2700 1018}%
\special{pa 2710 1048}%
\special{pa 2720 1080}%
\special{pa 2728 1110}%
\special{pa 2736 1142}%
\special{pa 2742 1172}%
\special{pa 2748 1204}%
\special{pa 2752 1236}%
\special{pa 2758 1268}%
\special{pa 2760 1300}%
\special{pa 2764 1332}%
\special{pa 2766 1364}%
\special{pa 2768 1396}%
\special{pa 2770 1428}%
\special{pa 2770 1460}%
\special{pa 2770 1492}%
\special{pa 2770 1524}%
\special{pa 2770 1556}%
\special{pa 2768 1588}%
\special{pa 2766 1620}%
\special{pa 2762 1650}%
\special{pa 2760 1682}%
\special{pa 2756 1714}%
\special{pa 2752 1746}%
\special{pa 2746 1778}%
\special{pa 2740 1810}%
\special{pa 2732 1840}%
\special{pa 2724 1872}%
\special{pa 2714 1902}%
\special{pa 2702 1930}%
\special{pa 2686 1960}%
\special{pa 2666 1984}%
\special{pa 2638 1998}%
\special{pa 2608 1992}%
\special{pa 2584 1970}%
\special{pa 2566 1944}%
\special{pa 2552 1916}%
\special{pa 2540 1886}%
\special{pa 2530 1854}%
\special{pa 2522 1824}%
\special{pa 2514 1792}%
\special{pa 2506 1762}%
\special{pa 2502 1730}%
\special{pa 2496 1698}%
\special{pa 2492 1666}%
\special{pa 2488 1636}%
\special{pa 2484 1604}%
\special{pa 2482 1572}%
\special{pa 2480 1540}%
\special{pa 2478 1508}%
\special{pa 2478 1476}%
\special{pa 2478 1444}%
\special{pa 2478 1412}%
\special{pa 2478 1380}%
\special{pa 2480 1348}%
\special{pa 2482 1316}%
\special{pa 2484 1284}%
\special{pa 2488 1252}%
\special{pa 2492 1220}%
\special{pa 2496 1188}%
\special{pa 2500 1156}%
\special{pa 2506 1126}%
\special{pa 2514 1094}%
\special{pa 2522 1064}%
\special{pa 2532 1032}%
\special{pa 2544 1004}%
\special{pa 2560 974}%
\special{pa 2578 948}%
\special{pa 2604 932}%
\special{pa 2616 930}%
\special{sp}%
%
\special{pn 8}%
\special{pa 2486 1606}%
\special{pa 2758 1322}%
\special{dt 0.045}%
%
\special{pn 8}%
\special{pa 2290 1360}%
\special{pa 3006 1578}%
\special{dt 0.045}%
%
\special{pn 20}%
\special{sh 1}%
\special{ar 3018 1574 10 10 0  6.28318530717959E+0000}%
\special{sh 1}%
\special{ar 3018 1574 10 10 0  6.28318530717959E+0000}%
%
\special{pn 20}%
\special{sh 1}%
\special{ar 2574 974 10 10 0  6.28318530717959E+0000}%
\special{sh 1}%
\special{ar 2574 974 10 10 0  6.28318530717959E+0000}%
%
\special{pn 8}%
\special{pa 2140 1530}%
\special{pa 3094 1388}%
\special{dt 0.045}%
%
\special{pn 20}%
\special{sh 1}%
\special{ar 2160 1530 10 10 0  6.28318530717959E+0000}%
\special{sh 1}%
\special{ar 2160 1530 10 10 0  6.28318530717959E+0000}%
%
\special{pn 20}%
\special{sh 1}%
\special{ar 3082 1388 10 10 0  6.28318530717959E+0000}%
\special{sh 1}%
\special{ar 3082 1388 10 10 0  6.28318530717959E+0000}%
%
\special{pn 20}%
\special{sh 1}%
\special{ar 2290 1350 10 10 0  6.28318530717959E+0000}%
\special{sh 1}%
\special{ar 2290 1350 10 10 0  6.28318530717959E+0000}%
%
\special{pn 20}%
\special{sh 1}%
\special{ar 2758 1322 10 10 0  6.28318530717959E+0000}%
\special{sh 1}%
\special{ar 2758 1322 10 10 0  6.28318530717959E+0000}%
%
\special{pn 20}%
\special{sh 1}%
\special{ar 2486 1618 10 10 0  6.28318530717959E+0000}%
\special{sh 1}%
\special{ar 2486 1618 10 10 0  6.28318530717959E+0000}%
%
\special{pn 20}%
\special{sh 1}%
\special{ar 2638 2000 10 10 0  6.28318530717959E+0000}%
\special{sh 1}%
\special{ar 2638 2000 10 10 0  6.28318530717959E+0000}%
%
\special{pn 8}%
\special{pa 2562 1922}%
\special{pa 2704 1016}%
\special{dt 0.045}%
%
\special{pn 8}%
\special{pa 2724 1846}%
\special{pa 2530 1028}%
\special{dt 0.045}%
%
\special{pn 20}%
\special{sh 1}%
\special{ar 2700 1010 10 10 0  6.28318530717959E+0000}%
\special{sh 1}%
\special{ar 2700 1010 10 10 0  6.28318530717959E+0000}%
%
\special{pn 20}%
\special{sh 1}%
\special{ar 2530 1038 10 10 0  6.28318530717959E+0000}%
\special{sh 1}%
\special{ar 2530 1038 10 10 0  6.28318530717959E+0000}%
%
\special{pn 20}%
\special{sh 1}%
\special{ar 2562 1934 10 10 0  6.28318530717959E+0000}%
\special{sh 1}%
\special{ar 2562 1934 10 10 0  6.28318530717959E+0000}%
%
\special{pn 20}%
\special{sh 1}%
\special{ar 2724 1858 10 10 0  6.28318530717959E+0000}%
\special{sh 1}%
\special{ar 2724 1858 10 10 0  6.28318530717959E+0000}%
%
\special{pn 20}%
\special{sh 1}%
\special{ar 3114 1530 10 10 0  6.28318530717959E+0000}%
\special{sh 1}%
\special{ar 3114 1530 10 10 0  6.28318530717959E+0000}%
\put(39.0000,-15.8000){\makebox(0,0)[lt]{$v$}}%
%
\special{pn 8}%
\special{ar 3386 1312 466 382  3.1615900 3.1898919}%
\special{ar 3386 1312 466 382  3.2747975 3.3030994}%
\special{ar 3386 1312 466 382  3.3880051 3.4163070}%
\special{ar 3386 1312 466 382  3.5012126 3.5295145}%
\special{ar 3386 1312 466 382  3.6144202 3.6427221}%
\special{ar 3386 1312 466 382  3.7276277 3.7559296}%
\special{ar 3386 1312 466 382  3.8408353 3.8691372}%
\special{ar 3386 1312 466 382  3.9540428 3.9823447}%
\special{ar 3386 1312 466 382  4.0672504 4.0955523}%
\special{ar 3386 1312 466 382  4.1804579 4.2087598}%
\special{ar 3386 1312 466 382  4.2936655 4.3219673}%
\special{ar 3386 1312 466 382  4.4068730 4.4351749}%
%
\special{pn 8}%
\special{pa 2920 1290}%
\special{pa 2920 1344}%
\special{dt 0.045}%
\special{sh 1}%
\special{pa 2920 1344}%
\special{pa 2940 1278}%
\special{pa 2920 1292}%
\special{pa 2900 1278}%
\special{pa 2920 1344}%
\special{fp}%
\put(33.0900,-10.0500){\makebox(0,0)[lb]{${\mathcal A}_x$}}%
\put(20.6300,-11.0400){\makebox(0,0)[rb]{${\mathcal A}_{x'}$}}%
\put(35.0000,-25.8000){\makebox(0,0)[lt]{$\,$}}%
\put(24.7500,-21.3000){\makebox(0,0)[lt]{$S(r_0)$}}%
%
\special{pn 8}%
\special{pa 2096 1028}%
\special{pa 2498 1148}%
\special{dt 0.045}%
\special{sh 1}%
\special{pa 2498 1148}%
\special{pa 2440 1110}%
\special{pa 2446 1132}%
\special{pa 2428 1148}%
\special{pa 2498 1148}%
\special{fp}%
%
\special{pn 8}%
\special{pa 2638 1464}%
\special{pa 3830 1618}%
\special{dt 0.045}%
%
\special{pn 8}%
\special{pa 2638 1464}%
\special{pa 2480 262}%
\special{dt 0.045}%
%
\special{pn 20}%
\special{sh 1}%
\special{ar 3830 1618 10 10 0  6.28318530717959E+0000}%
\special{sh 1}%
\special{ar 3830 1618 10 10 0  6.28318530717959E+0000}%
%
\special{pn 20}%
\special{sh 1}%
\special{ar 2476 274 10 10 0  6.28318530717959E+0000}%
\special{sh 1}%
\special{ar 2476 274 10 10 0  6.28318530717959E+0000}%
%
\special{pn 20}%
\special{ar 2628 1464 532 154  5.8097341 5.8448218}%
\special{ar 2628 1464 532 154  5.9500850 5.9851727}%
\special{ar 2628 1464 532 154  6.0904358 6.1255236}%
\special{ar 2628 1464 532 154  6.2307867 6.2658744}%
\special{ar 2628 1464 532 154  6.3711376 6.4062253}%
\special{ar 2628 1464 532 154  6.5114885 6.5465762}%
\special{ar 2628 1464 532 154  6.6518394 6.6869271}%
\special{ar 2628 1464 532 154  6.7921902 6.8272780}%
\special{ar 2628 1464 532 154  6.9325411 6.9676288}%
%
\special{pn 8}%
\special{pa 3386 1278}%
\special{pa 3158 1454}%
\special{dt 0.045}%
\special{sh 1}%
\special{pa 3158 1454}%
\special{pa 3224 1428}%
\special{pa 3200 1420}%
\special{pa 3200 1396}%
\special{pa 3158 1454}%
\special{fp}%
%
\special{pn 8}%
\special{pa 2898 712}%
\special{pa 2638 930}%
\special{dt 0.045}%
\special{sh 1}%
\special{pa 2638 930}%
\special{pa 2702 902}%
\special{pa 2680 896}%
\special{pa 2676 872}%
\special{pa 2638 930}%
\special{fp}%
\put(34.1800,-12.7800){\makebox(0,0)[lb]{$\sigma$}}%
\put(29.4100,-7.2100){\makebox(0,0)[lb]{$k'(\sigma)$}}%
\put(25.6000,-2.1000){\makebox(0,0)[rb]{$k'\cdot v$}}%
%
\special{pn 8}%
\special{ar 2638 1442 1208 1188  0.0000000 0.0100209}%
\special{ar 2638 1442 1208 1188  0.0400835 0.0501044}%
\special{ar 2638 1442 1208 1188  0.0801670 0.0901879}%
\special{ar 2638 1442 1208 1188  0.1202505 0.1302714}%
\special{ar 2638 1442 1208 1188  0.1603340 0.1703549}%
\special{ar 2638 1442 1208 1188  0.2004175 0.2104384}%
\special{ar 2638 1442 1208 1188  0.2405010 0.2505219}%
\special{ar 2638 1442 1208 1188  0.2805846 0.2906054}%
\special{ar 2638 1442 1208 1188  0.3206681 0.3306889}%
\special{ar 2638 1442 1208 1188  0.3607516 0.3707724}%
\special{ar 2638 1442 1208 1188  0.4008351 0.4108559}%
\special{ar 2638 1442 1208 1188  0.4409186 0.4509395}%
\special{ar 2638 1442 1208 1188  0.4810021 0.4910230}%
\special{ar 2638 1442 1208 1188  0.5210856 0.5311065}%
\special{ar 2638 1442 1208 1188  0.5611691 0.5711900}%
\special{ar 2638 1442 1208 1188  0.6012526 0.6112735}%
\special{ar 2638 1442 1208 1188  0.6413361 0.6513570}%
\special{ar 2638 1442 1208 1188  0.6814196 0.6914405}%
\special{ar 2638 1442 1208 1188  0.7215031 0.7315240}%
\special{ar 2638 1442 1208 1188  0.7615866 0.7716075}%
\special{ar 2638 1442 1208 1188  0.8016701 0.8116910}%
\special{ar 2638 1442 1208 1188  0.8417537 0.8517745}%
\special{ar 2638 1442 1208 1188  0.8818372 0.8918580}%
\special{ar 2638 1442 1208 1188  0.9219207 0.9319415}%
\special{ar 2638 1442 1208 1188  0.9620042 0.9720251}%
\special{ar 2638 1442 1208 1188  1.0020877 1.0121086}%
\special{ar 2638 1442 1208 1188  1.0421712 1.0521921}%
\special{ar 2638 1442 1208 1188  1.0822547 1.0922756}%
\special{ar 2638 1442 1208 1188  1.1223382 1.1323591}%
\special{ar 2638 1442 1208 1188  1.1624217 1.1724426}%
\special{ar 2638 1442 1208 1188  1.2025052 1.2125261}%
\special{ar 2638 1442 1208 1188  1.2425887 1.2526096}%
\special{ar 2638 1442 1208 1188  1.2826722 1.2926931}%
\special{ar 2638 1442 1208 1188  1.3227557 1.3327766}%
\special{ar 2638 1442 1208 1188  1.3628392 1.3728601}%
\special{ar 2638 1442 1208 1188  1.4029228 1.4129436}%
\special{ar 2638 1442 1208 1188  1.4430063 1.4530271}%
\special{ar 2638 1442 1208 1188  1.4830898 1.4931106}%
\special{ar 2638 1442 1208 1188  1.5231733 1.5331942}%
\special{ar 2638 1442 1208 1188  1.5632568 1.5732777}%
\special{ar 2638 1442 1208 1188  1.6033403 1.6133612}%
\special{ar 2638 1442 1208 1188  1.6434238 1.6534447}%
\special{ar 2638 1442 1208 1188  1.6835073 1.6935282}%
\special{ar 2638 1442 1208 1188  1.7235908 1.7336117}%
\special{ar 2638 1442 1208 1188  1.7636743 1.7736952}%
\special{ar 2638 1442 1208 1188  1.8037578 1.8137787}%
\special{ar 2638 1442 1208 1188  1.8438413 1.8538622}%
\special{ar 2638 1442 1208 1188  1.8839248 1.8939457}%
\special{ar 2638 1442 1208 1188  1.9240084 1.9340292}%
\special{ar 2638 1442 1208 1188  1.9640919 1.9741127}%
\special{ar 2638 1442 1208 1188  2.0041754 2.0141962}%
\special{ar 2638 1442 1208 1188  2.0442589 2.0542797}%
\special{ar 2638 1442 1208 1188  2.0843424 2.0943633}%
\special{ar 2638 1442 1208 1188  2.1244259 2.1344468}%
\special{ar 2638 1442 1208 1188  2.1645094 2.1745303}%
\special{ar 2638 1442 1208 1188  2.2045929 2.2146138}%
\special{ar 2638 1442 1208 1188  2.2446764 2.2546973}%
\special{ar 2638 1442 1208 1188  2.2847599 2.2947808}%
\special{ar 2638 1442 1208 1188  2.3248434 2.3348643}%
\special{ar 2638 1442 1208 1188  2.3649269 2.3749478}%
\special{ar 2638 1442 1208 1188  2.4050104 2.4150313}%
\special{ar 2638 1442 1208 1188  2.4450939 2.4551148}%
\special{ar 2638 1442 1208 1188  2.4851775 2.4951983}%
\special{ar 2638 1442 1208 1188  2.5252610 2.5352818}%
\special{ar 2638 1442 1208 1188  2.5653445 2.5753653}%
\special{ar 2638 1442 1208 1188  2.6054280 2.6154489}%
\special{ar 2638 1442 1208 1188  2.6455115 2.6555324}%
\special{ar 2638 1442 1208 1188  2.6855950 2.6956159}%
\special{ar 2638 1442 1208 1188  2.7256785 2.7356994}%
\special{ar 2638 1442 1208 1188  2.7657620 2.7757829}%
\special{ar 2638 1442 1208 1188  2.8058455 2.8158664}%
\special{ar 2638 1442 1208 1188  2.8459290 2.8559499}%
\special{ar 2638 1442 1208 1188  2.8860125 2.8960334}%
\special{ar 2638 1442 1208 1188  2.9260960 2.9361169}%
\special{ar 2638 1442 1208 1188  2.9661795 2.9762004}%
\special{ar 2638 1442 1208 1188  3.0062630 3.0162839}%
\special{ar 2638 1442 1208 1188  3.0463466 3.0563674}%
\special{ar 2638 1442 1208 1188  3.0864301 3.0964509}%
\special{ar 2638 1442 1208 1188  3.1265136 3.1365344}%
\special{ar 2638 1442 1208 1188  3.1665971 3.1766180}%
\special{ar 2638 1442 1208 1188  3.2066806 3.2167015}%
\special{ar 2638 1442 1208 1188  3.2467641 3.2567850}%
\special{ar 2638 1442 1208 1188  3.2868476 3.2968685}%
\special{ar 2638 1442 1208 1188  3.3269311 3.3369520}%
\special{ar 2638 1442 1208 1188  3.3670146 3.3770355}%
\special{ar 2638 1442 1208 1188  3.4070981 3.4171190}%
\special{ar 2638 1442 1208 1188  3.4471816 3.4572025}%
\special{ar 2638 1442 1208 1188  3.4872651 3.4972860}%
\special{ar 2638 1442 1208 1188  3.5273486 3.5373695}%
\special{ar 2638 1442 1208 1188  3.5674322 3.5774530}%
\special{ar 2638 1442 1208 1188  3.6075157 3.6175365}%
\special{ar 2638 1442 1208 1188  3.6475992 3.6576200}%
\special{ar 2638 1442 1208 1188  3.6876827 3.6977035}%
\special{ar 2638 1442 1208 1188  3.7277662 3.7377871}%
\special{ar 2638 1442 1208 1188  3.7678497 3.7778706}%
\special{ar 2638 1442 1208 1188  3.8079332 3.8179541}%
\special{ar 2638 1442 1208 1188  3.8480167 3.8580376}%
\special{ar 2638 1442 1208 1188  3.8881002 3.8981211}%
\special{ar 2638 1442 1208 1188  3.9281837 3.9382046}%
\special{ar 2638 1442 1208 1188  3.9682672 3.9782881}%
\special{ar 2638 1442 1208 1188  4.0083507 4.0183716}%
\special{ar 2638 1442 1208 1188  4.0484342 4.0584551}%
\special{ar 2638 1442 1208 1188  4.0885177 4.0985386}%
\special{ar 2638 1442 1208 1188  4.1286013 4.1386221}%
\special{ar 2638 1442 1208 1188  4.1686848 4.1787056}%
\special{ar 2638 1442 1208 1188  4.2087683 4.2187891}%
\special{ar 2638 1442 1208 1188  4.2488518 4.2588727}%
\special{ar 2638 1442 1208 1188  4.2889353 4.2989562}%
\special{ar 2638 1442 1208 1188  4.3290188 4.3390397}%
\special{ar 2638 1442 1208 1188  4.3691023 4.3791232}%
\special{ar 2638 1442 1208 1188  4.4091858 4.4192067}%
\special{ar 2638 1442 1208 1188  4.4492693 4.4592902}%
\special{ar 2638 1442 1208 1188  4.4893528 4.4993737}%
\special{ar 2638 1442 1208 1188  4.5294363 4.5394572}%
\special{ar 2638 1442 1208 1188  4.5695198 4.5795407}%
\special{ar 2638 1442 1208 1188  4.6096033 4.6196242}%
\special{ar 2638 1442 1208 1188  4.6496868 4.6597077}%
\special{ar 2638 1442 1208 1188  4.6897704 4.6997912}%
\special{ar 2638 1442 1208 1188  4.7298539 4.7398747}%
\special{ar 2638 1442 1208 1188  4.7699374 4.7799582}%
\special{ar 2638 1442 1208 1188  4.8100209 4.8200418}%
\special{ar 2638 1442 1208 1188  4.8501044 4.8601253}%
\special{ar 2638 1442 1208 1188  4.8901879 4.9002088}%
\special{ar 2638 1442 1208 1188  4.9302714 4.9402923}%
\special{ar 2638 1442 1208 1188  4.9703549 4.9803758}%
\special{ar 2638 1442 1208 1188  5.0104384 5.0204593}%
\special{ar 2638 1442 1208 1188  5.0505219 5.0605428}%
\special{ar 2638 1442 1208 1188  5.0906054 5.1006263}%
\special{ar 2638 1442 1208 1188  5.1306889 5.1407098}%
\special{ar 2638 1442 1208 1188  5.1707724 5.1807933}%
\special{ar 2638 1442 1208 1188  5.2108559 5.2208768}%
\special{ar 2638 1442 1208 1188  5.2509395 5.2609603}%
\special{ar 2638 1442 1208 1188  5.2910230 5.3010438}%
\special{ar 2638 1442 1208 1188  5.3311065 5.3411273}%
\special{ar 2638 1442 1208 1188  5.3711900 5.3812109}%
\special{ar 2638 1442 1208 1188  5.4112735 5.4212944}%
\special{ar 2638 1442 1208 1188  5.4513570 5.4613779}%
\special{ar 2638 1442 1208 1188  5.4914405 5.5014614}%
\special{ar 2638 1442 1208 1188  5.5315240 5.5415449}%
\special{ar 2638 1442 1208 1188  5.5716075 5.5816284}%
\special{ar 2638 1442 1208 1188  5.6116910 5.6217119}%
\special{ar 2638 1442 1208 1188  5.6517745 5.6617954}%
\special{ar 2638 1442 1208 1188  5.6918580 5.7018789}%
\special{ar 2638 1442 1208 1188  5.7319415 5.7419624}%
\special{ar 2638 1442 1208 1188  5.7720251 5.7820459}%
\special{ar 2638 1442 1208 1188  5.8121086 5.8221294}%
\special{ar 2638 1442 1208 1188  5.8521921 5.8622129}%
\special{ar 2638 1442 1208 1188  5.8922756 5.9022965}%
\special{ar 2638 1442 1208 1188  5.9323591 5.9423800}%
\special{ar 2638 1442 1208 1188  5.9724426 5.9824635}%
\special{ar 2638 1442 1208 1188  6.0125261 6.0225470}%
\special{ar 2638 1442 1208 1188  6.0526096 6.0626305}%
\special{ar 2638 1442 1208 1188  6.0926931 6.1027140}%
\special{ar 2638 1442 1208 1188  6.1327766 6.1427975}%
\special{ar 2638 1442 1208 1188  6.1728601 6.1828810}%
\special{ar 2638 1442 1208 1188  6.2129436 6.2229645}%
\special{ar 2638 1442 1208 1188  6.2530271 6.2630480}%
%
\special{pn 20}%
\special{ar 2620 1460 140 530  4.0130720 5.3173802}%
\end{picture}%
\hspace{2.5truecm}}
\caption{The action defined by the topological Tits building of an isoparametric submanifold}
\label{figure7}
\end{figure}

\vspace{1truecm}

{\small\textit{Department of Mathematics, Faculty of Science, 
Tokyo University of Science,}}

{\small\textit{1-3 Kagurazaka Shinjuku-ku, Tokyo 162-8601 Japan}}

{\small\textit{E-mail address}: koike@rs.kagu.tus.ac.jp}

\end{document}